\documentclass[times]{nmeauth}
\usepackage{moreverb}
\usepackage[colorlinks,bookmarksopen,bookmarksnumbered,citecolor=red,urlcolor=red]{hyperref}
\usepackage{float}
\usepackage{appendix}
\usepackage{subcaption}
\usepackage{amsmath}
\usepackage{subcaption}
\usepackage{sansmath}
\usepackage{scalerel}
\usepackage{soul}
\usepackage{multirow}
\usepackage{bm}
\usepackage{mathrsfs}

\graphicspath{ {./figures/} }

\DeclareMathOperator*{\assemble}{\scalerel*{\mathsf{A}}{\sum}}

\newcommand\BibTeX{{\rmfamily B\kern-.05em \textsc{i\kern-.025em b}\kern-.08em
T\kern-.1667em\lower.7ex\hbox{E}\kern-.125emX}}

\newcommand{\narrowtimes}{{\mkern-2mu\times\mkern-2mu}}

\renewcommand{\vec}[1]{\ensuremath{\bm{#1}}}
\newcommand{\ary}[1]{\boldsymbol { \mathsf{#1}}}

\begin{document}

\runningheads{Z Liu et al.}{Coupled isogeometric FEM/BEM for structural-acoustic analysis}

\title{Isogeometric FEM-BEM coupled structural-acoustic analysis of shells using subdivision surfaces}

\author{Zhaowei Liu\affil{1}, Musabbir Majeed\affil{2}, Fehmi Cirak\affil{2}, Robert N. Simpson\affil{1}\corrauth}
\address{\affilnum{1}School of Engineering, University of Glasgow, Glasgow, G12 8QQ, UK \break \affilnum{2}Department of Engineering, University of Cambridge, Trumpington Street, Cambridge CB2 1PZ, UK
}

\corraddr{School of Engineering, University of Glasgow, Glasgow, G12 8QQ, UK. E-mail: robert.simpson.2@glasgow.ac.uk}

\begin{abstract}
We introduce a coupled finite and boundary element formulation for acoustic scattering analysis over thin shell structures. A  triangular Loop subdivision surface discretisation is used for both geometry and analysis fields. The Kirchhoff-Love shell equation is discretised with the finite element method and the Helmholtz equation for the acoustic field with the boundary element method. The use of the boundary element formulation allows the elegant handling of infinite domains and precludes the need for volumetric meshing. In the present work the subdivision control meshes for the shell displacements and the acoustic pressures have the same resolution. The corresponding smooth subdivision basis functions have the $C^1$ continuity property required for the Kirchhoff-Love formulation and are highly efficient for the acoustic field computations. We verify the proposed isogeometric formulation through a closed-form solution of acoustic scattering over a thin shell sphere. Furthermore, we demonstrate the ability of the proposed approach to handle complex geometries with arbitrary topology that provides an integrated isogeometric design and analysis workflow for coupled structural-acoustic analysis of shells.
\end{abstract}


\maketitle

\vspace{-6pt}

%  *****************************************************
% **************** INTRODUCTION *******************
% *****************************************************

\section{Introduction}
Structural-acoustic interaction plays a key role in the design of components and structures in a wide variety of applications where vibrational noise minimisation is a major design requirement, such as in  aerospace and automotive engineering, in the design of materials used to control acoustic absorption, or sound radiation from submerged structures, such as submarines and ships. Numerical methods like the finite element method (FEM), finite difference method (FDM) and boundary element method (BEM) play a key role in the design of such structures through simulations of prototype designs. In recent years where advances in manufacturing allow component designs with ever increasing geometric complexity the importance of numerical methods is becoming more apparent.  Industrial designers are continually striving for designs that can efficiently deliver improved performance. This necessitates design workflows where computer-aided design (CAD) and analysis methods are tightly integrated. In addition, recent developments into novel materials that exhibit unique physical properties, i.e.  metamaterials~\cite{cummer2016controlling},  coupled with the increased fidelity of modern manufacturing processes has opened up new design possibilities, but the lack of design tools that integrate geometric design, analysis and optimisation technologies is widely acknowledged as a key challenge that must be addressed before such materials can be used for practical applications~\cite{cirakscott2002, cottrell2006isogeometric, SchDeScEvBoRaHu12}.

The engineering design workflows  possible with most commercial software are based on disparate geometry and analysis models where expensive and error-prone model conversion processes are required. Considering the iterative nature of design, such conversion processes can dominate and impede the design process so that alternative solutions are sought. Much research has been presented on the idea of adopting a single common model to represent geometry and analysis fields, with such methods often falling under the umbrella of `isogeometric analysis' (a term initially coined by Hughes et al.~\cite{hughes2005isogeometric}). A number of isogeometric approaches have, and are, being developed based on non-uniform rational B-splines (NURBS)~\cite{hughes2005isogeometric}, subdivision surfaces~\cite{cirakortiz2000} and other geometry representations~\cite{nguyen2016c,majeedCirak:2016}. Locally refinable variants of NURBS  include  the T-splines~\cite{bazilevs2010isogeometric} and PHT-splines~\cite{nguyen2011rotation}, and  the corresponding locally refinable subdivision surfaces include THCCS~\cite{wei2015truncated} and CHARMS~\cite{Grinspun:2002aa}. Almost all the mentioned geometry representations are intended for surfaces so that their application to shell and BEM analysis is straightfoward. However, they need to be suitably extended for isogeometric FEM approaches that require volumetric discretisations.  In general, CAD software does not provide such volumetric discretisations but we note that recent research is progressing towards automatic generation of analysis-suitable trivariate splines, e.g.~\cite{zhang2016geometric, hu2016centroidal}, and alternative immersed/embedded methods that do not require a boundary-fitted volume mesh~\cite{Sanches:2011aa, SchDeScEvBoRaHu12}.

For structural-acoustic analysis problems, where the structure can be modelled as a thin-shell and the acoustic pressure field is governed by the time-harmonic Helmholtz equation, a sensible approach is to adopt a coupled FEM/BEM formulation \cite{ fritze2005fem, chen2014fem}.  In this way the infinite fluid domain is modelled through a BEM  discretisation which exhibits significant advantages over a traditional FEM approach.  For the latter case, a na\"{i}ve truncation of meshes that represent  infinite domains will lead to unphysical wave reflections that require appropriate absorbing boundary conditions at truncated boundaries,~\cite{webb1989absorbing,safjan1998highly,djellouli2000finite,  harari2004analytical} and an insufficient mesh resolution can lead to significant dispersion error for high frequency problems. In contrast, BEM formulations do not suffer from unphysical reflected waves and dispersion error but more crucially,  only a surface discretisation is required to model infinite acoustic domains. This makes an isogeometric approach for  coupled FEM/BEM analysis with shells a particularly attractive approach where high order discretisations generated through CAD software can be used directly without any requirement for model conversion, meshing or generation of volumetric splines. In addition, high order BEM discretisations are generally accepted as superior over equivalent low order discretisations for acoustic problems \cite{marburg2002six, marburg2003influence}. Several isogeometric BEM approaches have already been developed based on NURBS~\cite{simpsonbordas2012, simpsonscott2014, peaketrevelyan2013}, T-splines~\cite{scottsimpson2012, ginnis2014isogeometric} and subdivision surfaces~\cite{bandara2015boundary}. 

In applications where the structure can be approximated as a thin-shell  the Kirchhoff-Love model leads to particularly robust and simple finite element formulations with only the nodal displacements as the degrees of freedom. The most challenging aspect of the discretisation of the Kirchhoff-Love formulation is the need for $C^1$ continuity prompting the development of high order smooth discretisations. Several NURBS-based isogeometric shell formulations have been proposed in which a high order CAD discretisation is used as a basis for both geometry and analysis thus simultaneously satisfying continuity requirements and side-stepping mesh generation~\cite{kiendl2009isogeometric, benson2010isogeometric, benson2011large}.  An  earlier approach made use of subdivision surfaces~\cite{cirakortiz2000,Cirak:2001aa} which illustrated automatic satisfaction of the $C^1$ continuity requirement and an ability to handle geometries of arbitrary topology. Since this seminal work, the approach has been extended to shape optimisation~\cite{bandara2015boundary}, non-manifold shell geometries~\cite{cirak2011subdivision} and thick shells~\cite{Long:2012aa}. The present study proposes a method for performing coupled acoustic-structural analysis for shell structures using a common discretisation for geometry and analysis through subdivision basis functions. This facilitates the use of identical basis functions for geometric modelling, structural displacements and acoustic pressures. We adopt the Loop subdivision scheme~\cite{loop1987smooth} in which global basis functions are constructed from quartic box-splines defined over a triangular surface control mesh. We note that a somewhat similar coupled isogeometric BEM/FEM approach has been presented for simulating the Stokes flow~\cite{heltai2016interaction}. However, here we distinguish the novelty of our work through the use of subdivision basis functions. 

We organise the paper as follows: a brief overview of Loop subdivision surfaces and their use for numerical analysis is given; the formulation of a BEM discretisation for Helmholtz analysis with subdivision basis functions is introduced; an outline of the coupled system of equations that govern acoustic-structural analysis with Kirchhoff-Love shells is presented;  implementation details on how to efficiently handle the large dense system of equations generated through the present approach are outlined; verification of the method through a closed-form solution for acoustic scattering over a thin-shell sphere is shown; and finally, the ability of the approach to analyse arbitrarily complex 3D geometries is demonstrated.  We note that the approach is restricted to time-harmonic problems and in the present work we consider  medium-frequency problems with a normalised wavenumber up to 80. In all problems it can be assumed that the fluid domain resides on the outside of the shell structure with the position vector $\mathbf{x} \in \mathbb{R}^3$.

% Analytic solutions of simple objects, such as sphere and plate, were generated by Junger et al. \cite{junger1986sound}. However, since the real underwater structures have complex geometries, analytic solutions are hard to obtain. Numerical methods are required for analysing this problem. The Finite Element Method is proved as an widely-used numerical method for analysing structural dynamics but it is not an efficient method for acoustic analysis. A number of works which use the Boundary Element Method have been done to solve acoustic problems. They have shown that the Boundary Element Method has advantages in analysing acoustic problems  \cite{ghosh1986new, wrobel2002boundary}. A FEM-BEM coupling method was presented by Fritze et al. to solve this coupled problem. The concept of Isogeometric analysis has also been adopted for structure-fluid interaction problem by Heltai et al.\cite{heltai2016interaction} using NURBS-based isogeometric discretisations.

\section{Loop Subdivision Surfaces}

Subdivision surfaces were introduced in the 1970s and are now widely used in computer graphics and animation~\cite{Zorin:2000aa,Peters:2008aa}. They are also available in most industrial CAD solid modelling packages, including Autodesk Fusion 360, PTC Creo and CATIA. In computer graphics literature, subdivision surfaces are usually viewed as a process for generating smooth limit surfaces through repeated refinement and smoothing of a control mesh. Alternatively, they can be viewed as the generalisation of splines to arbitrary connectivity meshes which is a viewpoint more suitable to finite and boundary element analysis. The subdivision surfaces inherit from the splines their refinability property so that all control meshes generated during subdivision refinement describe exactly the same spline surface.

In the present paper we consider the triangular subdivision surfaces proposed by Loop~\cite{loop1987smooth}, which generalise the quartic box-splines to arbitrary connectivity meshes. Quartic box-splines are defined on shift-invariant three-direction meshes in which each vertex is attached to six triangles. Note that, for the sake of brevity, the treatment of vertices located on the domain boundaries is omitted in this paper. In common with all subdivision schemes, in Loop subdivision each step consists of a refinement and an averaging step. In the refinement step the mesh is refined by splitting each triangle into four triangles, after bisecting the triangle edges. Subsequently, the coordinates of the vertices on the refined mesh are determined as the weighted average of the coordinates of their neighbouring vertices on the coarse mesh. For mesh regions consisting only of regular vertices with six adjacent triangles, the averaging weights are derived from the knot insertion rules for box-splines. For remaining irregular, also referred to as extraordinary or star, vertices the weights are derived through a spectral analysis of the subdivision matrix underlying the refinement process, see~\cite{Zorin:2000aa, Peters:2008aa}. The averaging weights depend only on the connectivity of the mesh but not the actual vertex coordinates. The resulting subdivision surface  is  $C^2$ continuous almost everywhere except at the irregular vertices where it is only $C^1$.

\begin{figure}[]
\centering
\begin{subfigure}{.6\textwidth}
 \centering
\includegraphics[width=\textwidth]{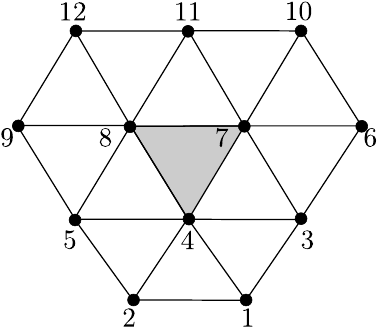}
\caption{A patch of elements required for interpolation over the centre triangle (shaded). The three vertices of the centre triangle are regular and the subdivision basis functions are identical to quartic box splines. The twelve non-zero quartic box-spline basis functions over the centre triangle correspond to the shown twelve vertices.} 
\label{fig:regular-patch}
\end{subfigure}
\\ \vspace{2ex}
\begin{subfigure}{.8\textwidth}
 \centering
\includegraphics[width=\textwidth]{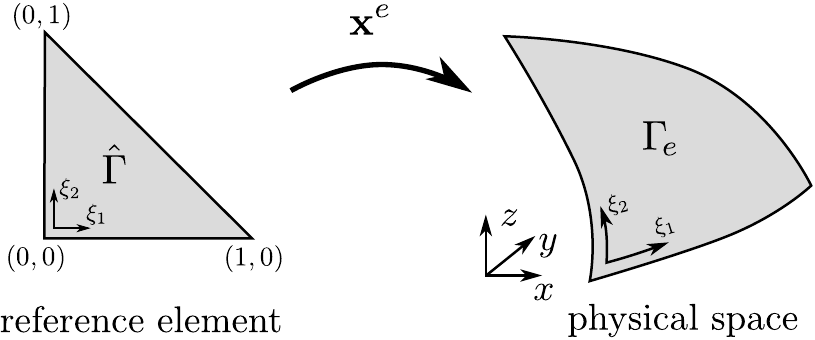}
\caption{Mapping of the reference triangle to the triangle on the subdivision surface. }
\label{fig:reference-element}
\end{subfigure}
\caption{Interpolation with Loop subdivision basis functions.}
\label{fig:test}
\end{figure}

In triangles with three regular vertices, there are twelve quartic non-zero basis functions associated with its and neighbouring triangle vertices, as can be seen in  Figure~\ref{fig:regular-patch}. Hence, the surface coordinates $\mathbf x^e$ within the triangle~$e$ can be determined through  
\begin{equation}
\label{eq:boxSpline_physicalcoord}
	\vec x^e(\xi_1, \xi_2) = \sum_{{a=1}}^{12} B_{a}(\xi_1, \xi_2) \vec P_a
\end{equation}
where $B_a$ are the box-spline basis functions,  $\vec \xi = (\xi_1, \xi_2)$ are two of the barycentric local coordinates in the triangle $e$ and $ \vec{P}_{a}$ are the coordinates of the twelve control vertices. In triangles with irregular vertices it is necessary to apply first a few steps of subdivision refinement until the considered point lies in a patch of elements so that~\eqref{eq:boxSpline_physicalcoord} can again be applied.  This is possible since  during subdivision refinement all newly created vertices are regular and repeated refinement leads to more and more patches with only regular vertices. As proposed by Stam~\cite{Stam1998a} this can be used to devise an algorithm for obtaining subdivision basis functions  $N_a (\vec \xi)$ for triangle patches that contain irregular vertices allowing the geometry to be interpolated as
\begin{equation}
\label{eq:loopSpline_physicalcoord}
	\vec x^e (\boldsymbol \xi) =  \sum_{{a=1}}^{n_v}   N_a (\boldsymbol \xi) \vec P_a  \, , 
\end{equation}
where $n_v$ are the number of vertices in the patch containing triangle $e$ and its neighbouring triangles (which share a vertex with $e$) and $\vec P_a$ are the coordinates of the vertices on the coarse control mesh. There is no closed-form expression for $N_a (\boldsymbol \xi) $ available, only an algorithm for evaluating it for almost any given $\boldsymbol \xi $. See~\cite{cirakortiz2000} for details. It is clear that on regular patches $N_a \equiv B_a$ so that in the following basis functions are always denoted with $N_a$.

The interpolation equations \eqref{eq:boxSpline_physicalcoord} and \eqref{eq:loopSpline_physicalcoord} can be adapted to provide a discretisation of analysis fields such as displacements in thin-shells and pressure in the  Helmholtz equation. The use of a common high-order basis for geometry and analysis makes such an approach inherently isogeometric.  It should also be noted that in comparison to traditional finite element interpolations the present Loop subdivision surfaces offer the advantage of providing unique normals and normal derivatives at vertices which both simplifies implementation and allows for superior accuracy.

\section{A coupled BEM/FEM with Loop subdivision surfaces basis functions}
\label{sec:coupled-formulation}

%
% Problem formulation
%
\subsection{Problem setup}
\label{sec:problem-formulation}

\begin{figure}[]
\centering
\begin{subfigure}{.6\textwidth}
	\centering
	\includegraphics[width=\textwidth]{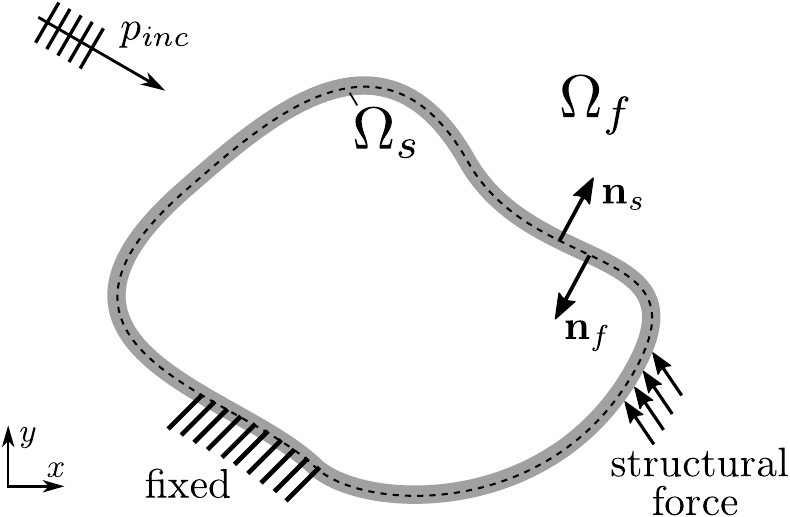}
	\caption{Definition of the structural thin-shell and fluid domains for the coupled structural acoustic problem with an impinged plane wave $p_{inc}$.   The mid-surface of the shell is denoted by a dashed line with $\mathbf{n}_s$ and $\mathbf{n}_f$ representing the normal vectors for the structural and fluid domains respectively.}
\end{subfigure} 
\\ \vspace{2ex}
\begin{subfigure}{.6\textwidth}
	\centering
	\includegraphics[width=\textwidth]{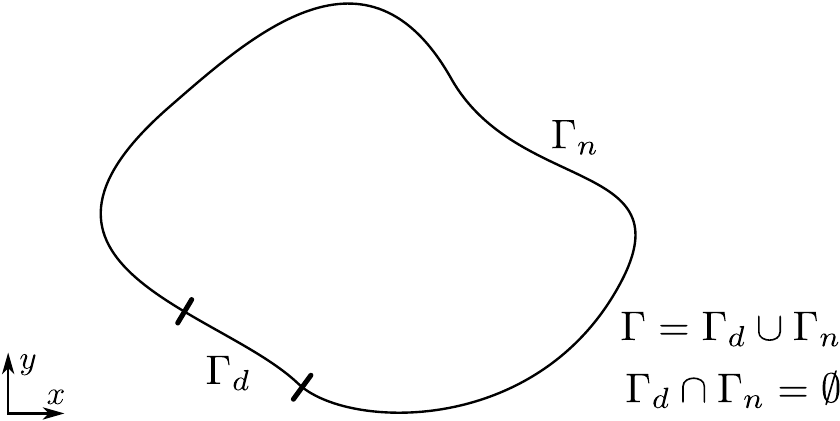}
	\caption{The interface boundary $\Gamma$ composed of Dirichlet and Neumann boundaries denoted by $\Gamma_d$ and $\Gamma_n$ respectively.}
\end{subfigure}
	\caption{Coupled problem setup.}
	\label{fig:domain-definitions}
\end{figure}

We first consider a domain $\Omega_{s}$ that represents an elastic thin-shell structure immersed in an infinite fluid domain $\Omega_{f}$. Furthermore, we assume that the behaviour of the thin-shell structure is governed by Kirchhoff-Love shell theory and all variables are time-harmonic with a time dependence of $e^{-i \omega t}$ with  $\omega$ denoting the angular frequency.  An acoustic pressure field exists in the fluid domain governed by the Helmholtz equation and the system is excited by an incident plane wave that is impinged on the shell structure with the entire setup depicted in Figure~\ref{fig:domain-definitions}. The acoustic pressure field at the fluid-structure interface induces displacements in the shell structure and, likewise, normal velocity components of the shell surface induce acoustic pressure gradients on the fluid domain forming a coupled system.   Our approach represents the fluid-structure interface with a Loop subdivision surface, and in order to construct the system of equations that governs the behaviour of the coupled system, first the discretisation of each domain is considered separately.  We then introduce the final discrete coupled system of equations in Section~\ref{sec:bem-fem-coupling}.

%
% BEM with subdivision surfaces
%
\subsection{Infinite fluid domain: collocation boundary element method with Loop subdivision surfaces}

In the infinite fluid domain $\Omega_f$ we wish to determine the complex-valued acoustic pressure field $p(\mathbf{x})$ given a wavenumber $k$ that is governed by the Helmholtz equation
\begin{equation}
\nabla^2 p(\mathbf{x}) + k^2 p(\mathbf{x}) = 0 \qquad  \text{in } \mathbf{x} \in \Omega_f \, ,
\label{eqn:helmholtz}
\end{equation} 
where $\nabla^2$ denotes the Laplacian operator. The (total) pressure field is composed of the incident and reflected acoustic pressures through
\begin{equation}
\label{eq:pressure}
p (\mathbf{x})  = p_{inc}(\mathbf{x}) +  p_{ref}(\mathbf{x}) \, ,
\end{equation}
and in the case of a plane wave of magnitude $P$ with wavenumber $k$ travelling in direction $\mathbf{d}$ with $|\mathbf{d}|=1$, the incident acoustic pressure is prescribed with $p_{inc}(\mathbf{x}) = Pe^{i\mathbf{k}\cdot\mathbf{x}}$ with the wavevector  $\mathbf{k} = k \mathbf{d}$. By setting the right-hand-side of \eqref{eqn:helmholtz} equal to the Dirac delta forcing function, the solution corresponds to the Helmholtz fundamental solution which is expressed as
\begin{equation}
\label{eqn:fundamental-solution}
G(\mathbf{x},\mathbf{y}) = \frac{e^{ikr}}{4\pi r} 
\end{equation}
with $r := |\mathbf{x}-\mathbf{y}|$ where $\mathbf{x}$ is the source point and $\mathbf{y}$ the field point.  The corresponding normal derivative of the kernel function of \eqref{eqn:fundamental-solution} is given by
\begin{equation}
\label{eqn:fundamental-solution-deriv}
\frac{\partial G(\mathbf{x},\mathbf{y})}{\partial n} = \frac{e^{ikr}}{4\pi r^2} (ikr - 1) \frac{\partial r}{\partial n}
\end{equation}
with $\partial (\cdot)/ \partial n \equiv \boldsymbol{\nabla} (\cdot) \cdot \mathbf{n}$.
By integrating (\ref{eqn:helmholtz}) over the domain $\Omega_f$ using (\ref{eqn:fundamental-solution}) as a weight function, applying the Green-Gauss theorem and then taking the limit as $\mathbf{x}$ approaches $\Gamma := \partial \Omega_f$, the acoustic boundary integral equation is expressed as
\begin{equation}
c(\mathbf{x})p(\mathbf{x})+\int_\Gamma \frac{\partial G(\mathbf{x},\mathbf{y})}{\partial n}p(\mathbf{y}) \,\mathrm{d}\Gamma(\mathbf{y}) = \int_\Gamma G(\mathbf{x},\mathbf{y}) \frac{\partial p(\mathbf{y})}{\partial n} \,\mathrm{d}\Gamma(\mathbf{y}) + p_{inc}(\mathbf{x}) \, ,
\label{eqn:CBIE}
\end{equation}
where $c(\mathbf{x})$ is a coefficient that depends on the geometry of the surface at the source point. For $\mathbf{x}$ located at a point with `smooth' geometry, $c(\mathbf{x}) = \frac{1}{2}$. 

\paragraph{Discretisation}
The boundary $\Gamma$ is defined in the usual piecewise manner through the union of elements, expressed as
\begin{equation}\label{eq:element-union}
\Gamma =\bigcup_{e=1}^{n_{el}} \Gamma_e \, ,
\end{equation}
where in the present work the set of elements $\{\Gamma_e\}_{e=1}^{n_{el}}$ is defined by the tessellation of the Loop subdivision surface.  The acoustic pressure and its normal derivative are discretised with subdivision basis functions $N_b(\boldsymbol \xi)$, cf.~\eqref{eq:loopSpline_physicalcoord},  
\begin{equation}
\label{eq:pressure-discretisation}
p^e(\boldsymbol{\xi}) = \sum_{b=1}^{n_v} N_{b}(\boldsymbol{\xi})p_{b}^e \quad \quad \frac{\partial p^e (\boldsymbol{\xi})}{\partial n} = \sum_{b=1}^{n_v}  N_{b}(\boldsymbol{\xi}) \frac{\partial p_{b}^e}{\partial n}
\end{equation}
where the coefficients $p_b^e$ and $\partial p_{b}^e / \partial n$ represent nodal coefficients of acoustic pressure and its normal derivatives respectively.
% Henceforth we will adopt the slightly abusive notation $N_{b}(\boldsymbol{v}) \equiv  N_{b}(\boldsymbol{u}(v,w))$. 
Using the tessellation \eqref{eq:element-union} and the approximations \eqref{eq:pressure-discretisation} the boundary integral equation given by \eqref{eqn:CBIE} is discretised as
\begin{equation}
\begin{aligned}
\label{eq:discretised-CBIE-1}
c(\mathbf{x})p(\mathbf{x})+&\sum_{e=1}^{n_{el}} \sum_{b=1}^{n_v} p^e_{b}  \int_{\Gamma_e} N_{b}(\boldsymbol{\xi}) \frac{\partial G(\mathbf{x},\mathbf{y(\boldsymbol{\xi})})}{\partial n} \,\mathrm{d}\Gamma_e (\boldsymbol{\xi}) &\\
= &\sum_{e=1}^{n_{el}} \sum_{b=1}^{n_v} \frac{\partial p^e_b}{\partial n} \int_{\Gamma_e} N_{b} (\boldsymbol{\xi}) G(\mathbf{x},\mathbf{y}(\boldsymbol{\xi})) \,\mathrm{d}\Gamma_e(\boldsymbol{\xi}) + p_{inc}(\mathbf{x}).&
\end{aligned}
\end{equation}

\paragraph{Collocation}
In order to generate a system of equations through \eqref{eq:discretised-CBIE-1} we adopt a collocation approach whereby the source point $\mathbf{x}$ is sampled at a discrete number of points on the boundary.  An alternative method is to apply a weighted residual (Galerkin) approach where \eqref{eq:discretised-CBIE-1} is multiplied by a suitable set of test functions and integrated over the boundary.   We adopt the former approach in the present work where significantly faster runtimes are achieved over an equivalent Galerkin implementation.

When spline based basis functions, like B-splines, NURBS, subdivision or T-spines, are employed different choices are possible for collocation points. For instance, Greville abscissae, Demko points and the maxima of B-splines have been considered in collocation methods that discretise directly the strong form of the governing equations~\cite{schillinger2013isogeometric}. In the present work we use the maxima of the subdivision basis functions, which are associated with the vertices of the tessellation, as the collocation points.  Note that due to the non-interpolatory nature  of the subdivision basis functions the coordinates of the collocation point located on the surface and the coordinates of the control vertices are not the same.

\begin{figure}[]
	\centering
	\includegraphics[width=0.6\textwidth]{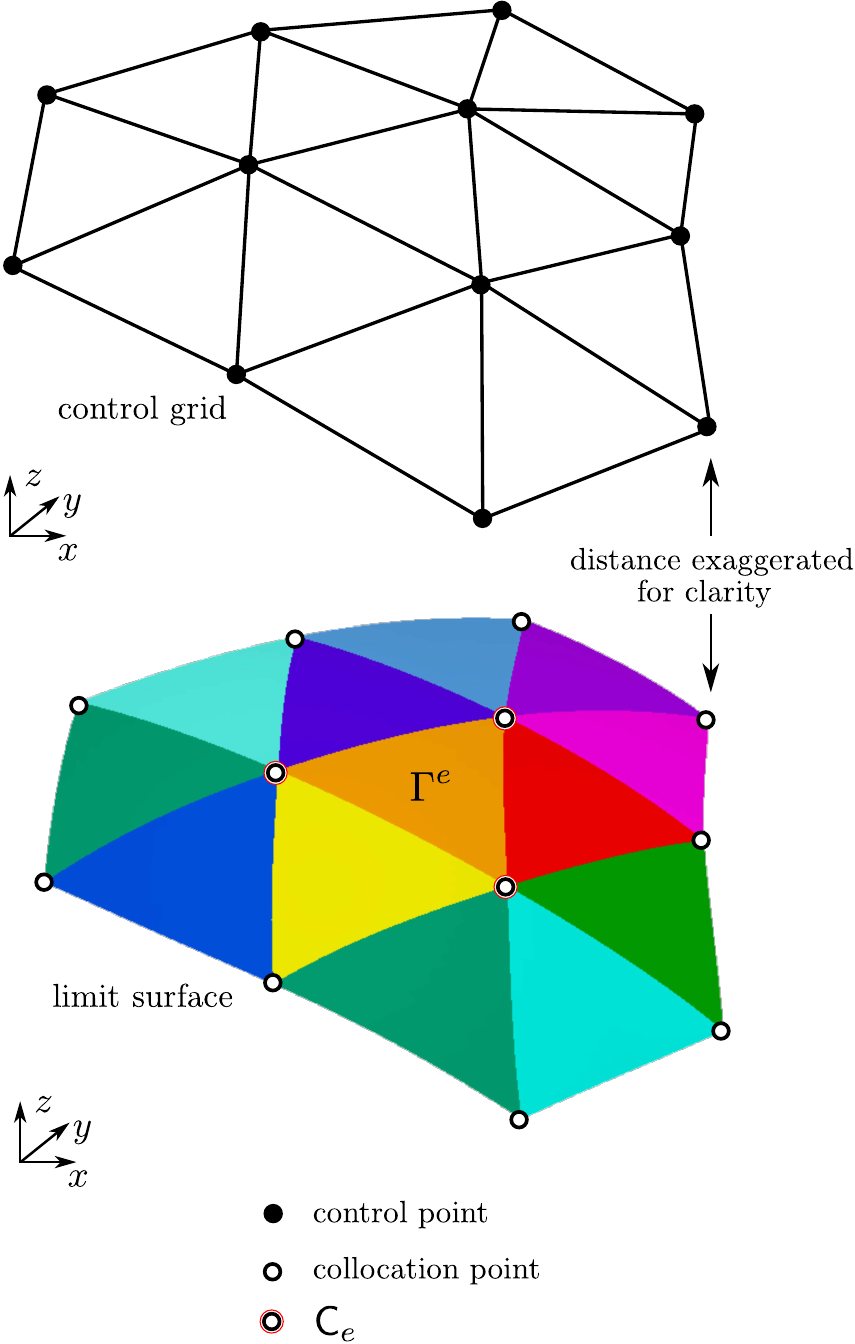}
	\caption{Illustration of a control mesh and associated limit surface for a regular patch using the Loop subdivision scheme. Control vertices are located in the control mesh at the intersection of edges and collocation points are defined as the corresponding points on the limit surface. The collocation points of the centre element specified by the set $\mathsf{C}_e$ are denoted by red circles.   It is important to note that control vertices and collocation points do not coincide.}
	\label{fig:collocation-point-limit-surface}
\end{figure}

For a given element $\Gamma_e$  we define the set of collocation points contained within the element through its three nodal points with the parametric coordinates
\begin{equation}
\mathsf{C}_e := \{\mathbf{x}^e(0,0), \mathbf{x}^e(1,0), \mathbf{x}^e(0,1) \}
\end{equation}
and construct the global set of collocation points as
\begin{equation}
\label{eq:collocation-global-set}
\mathsf{C} = \bigcup_{e=1}^{n_{el}} \mathsf{C}_e
\end{equation}
with duplicate entries discarded (i.e. each element of \eqref{eq:collocation-global-set} is unique).  For a regular patch of a Loop subdivision surface the scenario is depicted in Figure~\ref{fig:collocation-point-limit-surface}.   

\paragraph{System of equations}
By sampling over the set of collocation points $\mathsf{C}= \{\mathbf{x}_1, \mathbf{x}_2, \ldots \mathbf{x}_{n_{cp}}\}$ and employing an assembly operator $\assemble$ that maps a given element and local basis function index pairing $(e,b)$ to a global basis function index through $B = \assemble (e,b)$, the fully discrete boundary integral equation is written as
\begin{equation}
\label{eq:cbie-discrete-final}
\begin{aligned}
\frac{1}{2}p(\mathbf{x}_{A})+&\sum_{e=1}^{n_{el}} \sum_{b=1}^{n_v} p_{B}   \int_{\Gamma_e} N_{b}(\boldsymbol{\xi}) \frac{\partial G(\mathbf{x}_A,\mathbf{y}(\boldsymbol{\xi}))}{\partial n} \,\mathrm{d}\Gamma_e (\boldsymbol{\xi}) &\\
= &\sum_{e=1}^{n_{el}} \sum_{b=1}^{n_v} \frac{\partial p_B}{\partial n} \int_{\Gamma_e} N_{b} (\boldsymbol{\xi}) G(\mathbf{x}_A,\mathbf{y}(\boldsymbol{\xi})) \,\mathrm{d}\Gamma_e(\boldsymbol{\xi}) + p_{inc}(\mathbf{x}_A)&\\
&\qquad \qquad \qquad A = 1,2,\ldots, n_{cp}&\\
\end{aligned}
\end{equation}
where $p_B = p_{A(e,b)} = p_b^e$ and likewise for $\partial p_B / \partial n$. $n_{cp}$ is the number of collocation points. The final system of equations is then written as
\begin{equation}
\ary{H}\ary{p} = \ary{G}\ary{q} + \ary{p}_{inc}
\label{eqn:bem_planewave_hard}
\end{equation}
where  $\ary{p}$ and $\ary{q}$ represent vectors of global acoustic pressure and acoustic pressure normal derivative coefficients respectively, $\ary{H}$ and $\ary{G}$ are dense matrices with entries 
\begin{align}
\ary{H}_{AB} &= \frac{1}{2}N_B(\mathbf{x}_{A}) + \int_{\Gamma} N_{B}(\boldsymbol{y}) \frac{\partial G(\mathbf{x}_A,\mathbf{y})}{\partial n} \,\mathrm{d}\Gamma (\boldsymbol{y})\label{eq:H_AB-entry}\\
\ary{G}_{AB} &= \int_{\Gamma} N_{B} (\boldsymbol{y}) G(\mathbf{x}_A,\mathbf{y}) \,\mathrm{d}\Gamma(\boldsymbol{y}) \label{eq:G_AB-entry}
\end{align}
and $\ary{p}_{inc}$ is a vector with components $p_{inc}(\mathbf{x}_A)$. The integral in \eqref{eq:G_AB-entry} is found to be weakly singular and can therefore be integrated in a straightforward manner using e.g. polar integration (see \cite{scottsimpson2012} for details). In the present work we adopt a singularity subtraction technique \cite{liu1991some} whereby the integral in \eqref{eq:H_AB-entry} is regularised using the identity $\int_\Gamma \frac{\partial G^s(\mathbf{x}, \mathbf{y})}{\partial n} \, \mathrm{d}\Gamma(\mathbf{y})  = - \frac{1}{2}$ with $G^s(\mathbf{x}, \mathbf{y})$ denoting the `static' version (i.e. $k = 0$) of the Helmholtz kernel.

\subsection{Structural domain: finite element method for dynamic analysis of Kirchhoff-Love shells}
 
The mechanical response of thin-shells is approximated with a Kirchhoff-Love model.  In the following we provide a summary of the  corresponding weak form of  equilibrium equations and refer to~\cite{cirakortiz2000,Cirak:2001aa,Ciarlet:2005aa} for a more detailed presentation. The weak form, i.e. the virtual work expression, for the shell with a mid-surface~$\Gamma$ and displacements~$\vec u$ reads
\begin{equation}\label{eq:shellWeak}
	W_{mas} (\vec u, \vec v) + W_{int} (\vec u, \vec v) + W_{ext} (\vec v)= 0 \, .
\end{equation}
Here,  the three terms~$W_{mas} (\vec u, \vec v)$, $ W_{int} (\vec u, \vec v) $ and~$W_{ext} (\vec v)$ denote the virtual work contributions of the  inertia, internal and external forces respectively and $\vec v$ are the virtual displacements. For a thin-shell with aerial density $\rho_s$ the inertia contribution is given by
\begin{equation}
	W_{mas} (\vec u, \vec v) = \int_\Gamma \rho_s \frac{\partial^2 \vec u}{\partial t^2} \cdot \vec v  \mu \, \mathrm{d} \Gamma  \, 
\end{equation}
where $\mu$ is  a Jacobian taking care of integration across the thickness of the shell.  Only the acceleration of the mid-surface is considered. The contribution of the angular acceleration associated with the  shell mid-surface normal has been neglected. The internal work $W_{int} (\vec u, \vec v)$ consists of two parts, namely the membrane and bending parts:
 \begin{equation}
	W_{int}{ ( \vec u, \vec v) } =  \int_\Gamma  \vec \alpha (\vec u) : \vec E : \vec \alpha ( \vec v)  \mu \,\mathrm{d} \Gamma + h^2 \int_\Gamma  \vec \beta (\vec u) : \vec E : \vec \beta ( \vec v) \mu \,\mathrm{d} \Gamma \, .
\end{equation}
The membrane part, the first integral, depends on the constitutive tensor~$\vec E$ and change of the metric tensor~$\vec \alpha (\vec u)$ of the mid-surface between the reference and deformed configurations of the shell. The second integral representing the bending part is multiplied with the square of the shell thickness~$h$  and depends on change of curvature tensor~$\vec \beta(\vec u)$ between the reference and deformed configurations. Finally, the external virtual work for a shell loaded with pressure loading $p$, such as considered in this paper,  is given by 
 \begin{equation}
	W_{ext}{ ( \vec v) } = -  \int_\Gamma   p \vec n_S \cdot \vec v \,\mathrm{d} \Gamma \, 
\end{equation}
where $\vec n_S$ denotes the normal to the mid-surface.

The mid-surface displacements $\vec u$ are assumed to be time harmonic so that applying a Fourier transformation to the weak form~\eqref{eq:shellWeak} leads to its time-harmonic form. With a slight abuse of notation we denote in the following the Fourier transform of the displacements $\vec u$ with the same symbol; that is from now on $\vec u$ denotes the Fourier transform of displacements. The weak form for the shell in its time-harmonic form reads
\begin{equation} \label{eq:shellWeakHarmonic}
	-  \omega^2 \int_\Gamma \rho_s  \vec u  \cdot \vec v \mu \,\mathrm{d} \Gamma  + W_{int} (\vec u, \vec v) + W_{ext} (\vec v)= 0 
\end{equation}
with the angular velocity $\omega$. 

For discretising the displacements $\vec u$ and test functions $\vec v$ we use the same tessellation as for the fluid domain and the Loop subdivision basis functions $N_b (\boldsymbol \xi)$. 
%introduced in~\eqref{} and~\eqref{}
In each element $e$ the two fields $\vec u$  and $\vec v$ are approximated with 
\begin{equation}
	\vec u^e (\boldsymbol \xi) = \sum_{b=1}^{n_v} N_b(\boldsymbol \xi) \vec u^e_b \quad \text{and} \quad  \vec v^e (\boldsymbol \xi) = \sum_{b=1}^{n_v} N_b(\boldsymbol \xi) \vec v^e_b \, ,
\end{equation}   
 where the coefficients $\vec u^e_b$  and $\vec v^e_b$ can be interpreted as nodal quantities.  Introducing the discretisation into the time-harmonic weak form~\eqref{eq:shellWeakHarmonic} and evaluating the integrals with numerical quadrature yields, after linearisation, a system of equations that governs the time-harmonic behaviour of displacements 
\begin{equation}\label{eq:shell-matrix-formulation1}
( - \omega^2 \ary {M} + \ary K)\ary{u}  = \ary{f} \, ,
\end{equation}
where $\ary{K}$ is the stiffness matrix, $\ary{M}$ the mass matrix and $\ary{f}$ is the global force vector. For explicit expressions  for the stiffness matrix and implementation details see~\cite{cirakortiz2000}. In order to consider damping effects the discrete system of equations~\eqref{eq:shell-matrix-formulation1} can be augmented with a viscous Rayleigh damping term 
\begin{equation}\label{eq:shell-matrix-formulation2}
( - \omega^2 \ary {M} + i \omega (c_1 \ary K + c_2 \ary M )  +  \ary K)\ary{u}  = \ary{f}
\end{equation}
with $c_1$ and $c_2$ representing two experimentally determined constants. Finally, we simplify~\eqref{eq:shell-matrix-formulation2} to the form 
\begin{equation}\label{eq:shell-matrix-formulation3}
\ary{A} \ary{u} = \ary{f}.
\end{equation}

%
% Section - coupled formulation
%
\subsection{Coupled formulation}
\label{sec:bem-fem-coupling}
The pressure field in the acoustic fluid domain induces a force on the shell surface that is directed along the surface normal $\mathbf{n}_f$. The array of shell nodal forces $\ary{f}$ in~\eqref{eq:shell-matrix-formulation3}  due to a fluid pressure field interpolated according~\eqref{eq:pressure-discretisation}  is given by
\begin{equation}\label{eq:fluid-force}
\ary{f} =  \ary{\tilde n}_f  \int_\Gamma \ary{N}^T \ary{N}\, \mathrm{d}\Gamma \, \ary p  =  \ary{C}_{sf}\ary{p}
\end{equation}
with the matrix of vertex normals $\ary{\tilde n}_f$, matrix of subdivision basis functions~$\ary N$ and array of fluid nodal pressures $\ary p$ 
\begin{align} 
\ary{\tilde n}_f & = 
\begin{bmatrix}	
	\vec n_{1f} \cdot \vec e_1 & 0  & \ldots \\
	 \vec n_{1f} \cdot \vec e_2 & 0  & \ldots \\
	 \vec n_{1f} \cdot \vec e_3 & 0   & \ldots \\
	 0 & \vec n_{2f} \cdot \vec e_1  & \ldots \\
	 0 & \vec n_{2f} \cdot \vec e_2   & \ldots \\
	 0& \vec n_{2f} \cdot \vec e_3    & \ldots \\
	 \ldots & \ldots & \ldots \\
\end{bmatrix} 
\\
\ary{N} &= 
\begin{bmatrix}
N_1(\mathbf{x}) & N_2(\mathbf{x}) & \ldots & N_{n_{cp}}(\mathbf{x}) \label{eq:global-basis-vector}
\end{bmatrix}
\\
\ary{p} &= 
\begin{bmatrix}
p_1 &  p_2 & \ldots &  p_{n_{cp}}
\end{bmatrix}
\end{align}
where $(\vec e_1, \vec e_2, \vec e_3)$ are the three (orthogonal) base vectors such that each 
 column of $\ary {\tilde n}_f $ contains the normal components at each of the $n_{cp}$ control vertices in the mesh. The transfer matrix $\ary{C}_{sf}$ of dimension $3n_{cp} \times n_{cp}$ transfers nodal values from the fluid to the shell.

%Turning our attention now towards the relationship between acoustic pressure and shell mid-surface velocities, we define the acoustic admittance $\beta$ (which is in general a function of the position on the boundary and frequency).  
The normal components of the fluid and structural velocities denoted by $\ary{v}^{n}_f$ and $\ary{v}^{n}_s$ respectively can then be related to the acoustic pressure vector through 
%\begin{equation}\label{eq:admittance-relationship}
%\ary{v}^{n}_f - \ary{v}^{n}_s = \beta \ary{p} \, .
%\end{equation}

\begin{equation}\label{eq:admittance-relationship}
\ary{v}^{n}_f - \ary{v}^{n}_s = 0 .
\end{equation}

The acoustic pressure normal derivative $\ary{q}$ is known to be related to the fluid normal velocity $\ary{v}_f^n$ as
\begin{equation}\label{eq:acoustic-normal-deriv}
\ary{q} = - i\omega\rho_f\ary{v}^{n}_f \, , 
\end{equation}
where $\rho_f \equiv \rho$ is the density of fluid. We consider the problem as a fully coupled fluid structure interaction model and as such no velocity loss occurs between the fluid surface and structural surface but we remark that such models could be easily incorporated within our approach. The structural normal velocity $\ary{v}_s^n$ is related to the nodal displacements $\ary{u}$ through
\begin{equation}\label{eq:structural-velocity-relationship}
 \ary{v}^{n}_s =  i\omega\ary{C}_{fs}\ary{u} \, , 
\end{equation}
where $\ary{C}_{fs} = \ary{\tilde n}_f^{T}$.  
%transforms structural degrees of freedom to fluid degrees of freedom and is expressed as
%\begin{equation}\label{eq:coupled-matrix-Cfs}
%\ary{C}_{fs} = \ary{P}^{-1}\ary{C}_{sf}^{T}
%\end{equation}
%with the projection matrix $\ary{P}$ defined as
%\begin{equation}
%\ary{P} = \int_\Gamma \mathbf{N}^T\mathbf{N} \, \mathrm{d}\Gamma. \label{eq:projection-matrix}
%\end{equation}
%The use of $\ary{P}$ is required in \eqref{eq:coupled-matrix-Cfs} due the non-interpolatory nature of the Loop subdivision basis functions. 

Substituting \eqref{eq:admittance-relationship} and \eqref{eq:structural-velocity-relationship} into \eqref{eq:acoustic-normal-deriv} for $\ary{v}^{n}_f$ and $\ary{v}^{n}_s$ respectively,  the desired relationship between acoustic pressure normal derivatives and structural displacements is written as
%\begin{equation}
%\ary{q} =  - i\omega\rho\beta\ary{p} + \omega^2\rho\ary{C}_{fs}\ary{u} \;.
%\end{equation}

\begin{equation}
\ary{q} =  \omega^2\rho\ary{C}_{fs}\ary{u} \;.
\label{eq:acoustic-normal-deriv-final}
\end{equation}

%By introducing the diagonal matrix $\ary{Y} = - i\omega\rho\beta \ary{I}$ with $\ary{I}$ denoting the identity matrix, this is simplified to
%\begin{equation}\label{eq:acoustic-normal-deriv-final}
%\ary{q} =  \ary{Y}\ary{p} + \omega^2\rho\ary{C}_{fs}\ary{u}.
%\end{equation}
Using \eqref{eq:acoustic-normal-deriv-final} to substitute for $\ary{q}$ in the boundary element system of equations given by \eqref{eqn:bem_planewave_hard}, a coupled system for the acoustic problem is given by
%\begin{equation}
%[\ary{H} - \ary{G}\ary{Y}]\ary{p} = \ary{G} \omega^2\rho\ary{C}_{fs}\ary{u} + \ary{p}_{inc} \;.
%\label{eq:acoustic_coupled}
%\end{equation}

\begin{equation}
\ary{H} \ary{p} = \ary{G} \omega^2\rho\ary{C}_{fs}\ary{u} + \ary{p}_{inc} \;.
\label{eq:acoustic_coupled}
\end{equation}
Likewise, by substituting \eqref{eq:fluid-force} into \eqref{eq:shell-matrix-formulation3} a coupled system for the structural dynamics problem is written as
\begin{equation}\label{eq:structure-coupled}
\ary{A}\ary{u} = \ary{C}_{sf}\ary{p} + \ary{f}_s.
\end{equation}
where $\ary{f}_s$ contains nodal shell forces due to external loading other than fluid pressure. Finally, by combining \eqref{eq:acoustic_coupled} and \eqref{eq:structure-coupled}, the global coupled system of equations is expressed as

\begin{equation}\label{eq:coupled-system-1}
\left[
\begin{array}{cc}
\ary{A} & -\ary{C}_{sf}\\
-\omega^2\rho\ary{G}\ary{C}_{fs} & \ary{H}\\
\end{array}
\right]
\left[
\begin{array}{c}
\ary{u}\\
\ary{p}\\
\end{array}
\right] = 
\left[
\begin{array}{c}
\ary{f}_s\\
\ary{p}_{inc}\\
\end{array}
\right].
\end{equation}

%\begin{equation}\label{eq:coupled-system-1}
%\left[
%\begin{array}{cc}
%\ary{A} & -\ary{C}_{sf}\\
%-\omega^2\rho\ary{G}\ary{C}_{fs} & \ary{H} - \ary{G}\ary{Y}\\
%\end{array}
%\right]
%\left[
%\begin{array}{c}
%\ary{u}\\
%\ary{p}\\
%\end{array}
%\right] = 
%\left[
%\begin{array}{c}
%\ary{f}_s\\
%\ary{p}_{inc}\\
%\end{array}
%\right].
%\end{equation}

%
% Section : implementation 
%
\section{Implementation}
The system of equations given by \eqref{eq:coupled-system-1} is non-symmetric and contains dense partitions that arise from the dense matrices $\ary{H}$ and $\ary{G}$. Application of a direct solver would lead to inordinately long runtimes and excessive memory demands and we therefore outline a modified version of \eqref{eq:coupled-system-1} that is more amenable for computations.  We first rewrite \eqref{eq:coupled-system-1} using its Schur complement as

\begin{equation}\label{eq:system-equations-implementation-1}
\ary{H} \ary{p} =\omega^2\rho \ary{G}\ary{C}_{fs}\ary{A}^{-1}\left(\ary{f}_s +  \ary{C}_{sf}\ary{p} \right) + \ary{p}_{inc} \;,
\end{equation}
%\begin{equation}\label{eq:system-equations-implementation-1}
%[\ary{H}-\ary{G}\ary{Y}]\ary{p} =\omega^2\rho \ary{G}\ary{C}_{fs}\ary{A}^{-1}\left(\ary{f}_s +  \ary{C}_{sf}\ary{p} \right) + \ary{p}_{inc} \;,
%\end{equation}
where $\ary{u} = \ary{A}^{-1}( \ary{f}_s+\ary{C}_{sf}\ary{p})$ has been employed.  Defining the vector $\ary{q}_s$ which accounts for the contribution of acoustic velocities from the structural domain as
\begin{equation}\label{eq:qs-defn}
\ary{q}_s = \omega^2\rho\ary{C}_{fs}\ary{A}^{-1}\ary{f}_s
\end{equation}
and a global admittance matrix $\ary{Y}_C$ that represents the admittance effect caused by the structure as
\begin{equation}\label{eq:Yc-matrix-defn}
\ary{Y}_C = \omega^2\rho\ary{C}_{fs}\ary{A}^{-1}\ary{C}_{sf} \;.
\end{equation}
The system of equations of \eqref{eq:system-equations-implementation-1} is then written as

\begin{equation}\label{eq:system-equations-implementation-2}
[\ary{H}-\ary{G}\ary{Y}_C]\ary{p} = \ary{G}\ary{q}_s + \ary{p}_{inc} \;.
\end{equation}

%\begin{equation}\label{eq:system-equations-implementation-2}
%[\ary{H}-\ary{G}\ary{Y}-\ary{G}\ary{Y}_C]\ary{p} = \ary{G}\ary{q}_s + \ary{p}_{inc} \;.
%\end{equation}

When solving \eqref{eq:system-equations-implementation-2} two important considerations must be made: computation of $\ary{A}^{-1}$ % and  $\ary{P}^{-1}$
and efficient representation of the dense matrices $\ary{H}$ and $\ary{G}$. For the former, a sensible strategy is to approximate matrix inverses using a singular value decomposition (SVD) or a modal analysis approach in much the same manner as~\cite{fritze2005fem}. Several libraries exist which allow for efficient SVD computations  (e.g.~\cite{gough2009gnu, eigenweb, trilinos-overview}) but we note that for the examples in the present study no such approximations were required.

\begin{figure}[]
	\centering
	\includegraphics[width=0.8\textwidth]{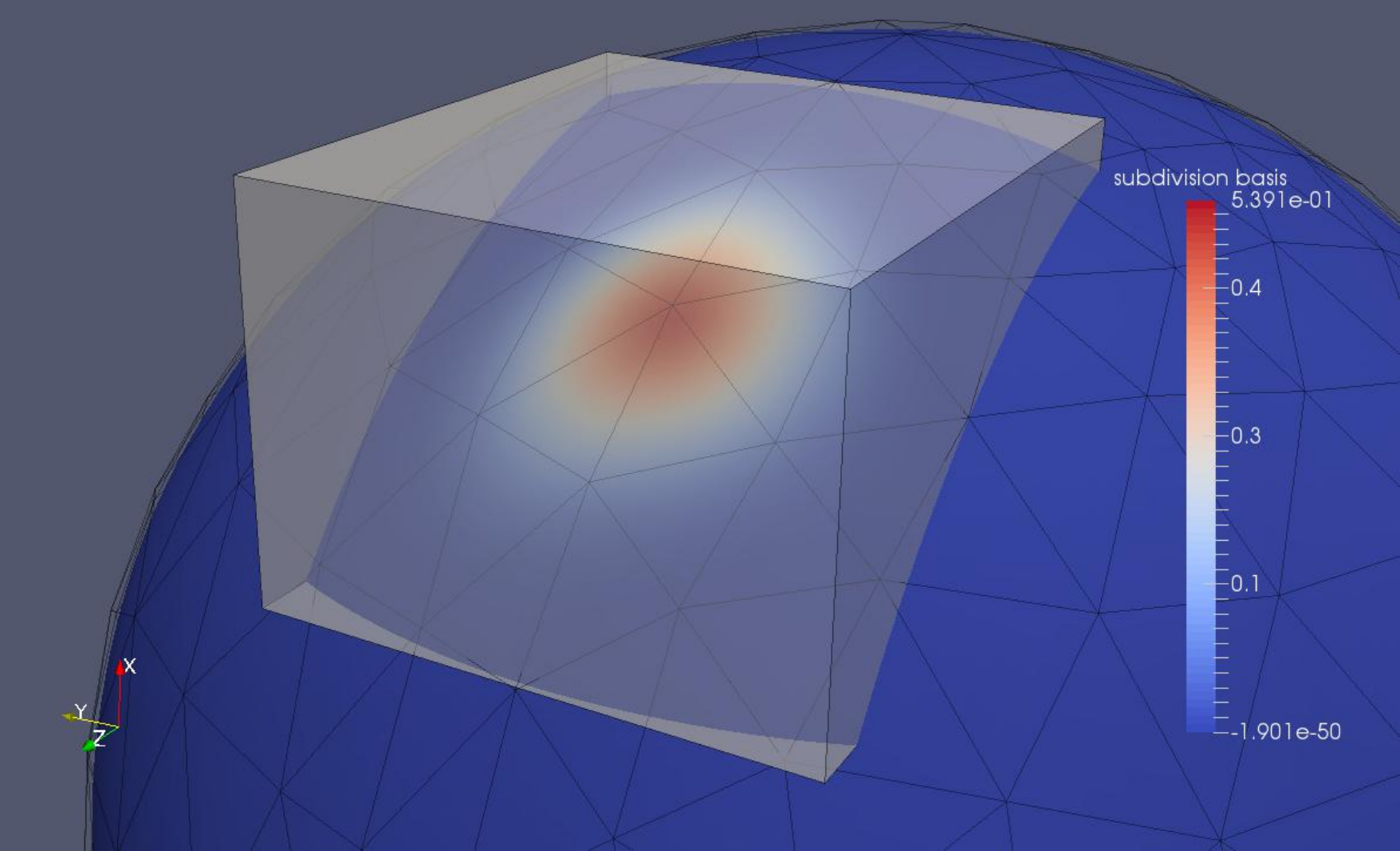}
	\caption{An example of a bounding box (translucent cuboid) defined by the support of a Loop subdivision basis function. The set of bounding boxes defined through the global set of subdivision basis functions over the surface is used to construct a cluster-tree that is used to build an $\mathscr{H}$-matrix approximation of the dense boundary element matrices.}
	\label{fig:bounding-box-closup}
\end{figure}

The representation of dense matrices requires a more involved implementation and in the present work we choose to adopt an $\mathscr{H}$-matrix approach which approximates dense matrices through low-rank approximations using the library HLibPro~\cite{kriemann2008hlibpro}. Without delving into the details of $\mathscr{H}$-matrix theory~\cite{hackbusch2004hierarchical, bebendorf2008hierarchical}, we simply state that the algorithm computes a low-rank approximation of $\ary{H}$ and $\ary{G}$ through a specified tolerance $\varepsilon$ by utilising a hierarchical `cluster tree' which separates terms into far-field (admissible) and near-field sets (non-admissible). The cluster tree is defined through a set of coordinates and bounding boxes related  to the underlying basis of the boundary element discretisation. In the present work we use the set of collocation points given by \eqref{eq:collocation-global-set} and the set of bounding boxes specified as $\mathsf{B} := \{Q_{min}(\mathbf{x}) : \mathbf{x} \in \textrm{supp}(N_A)\}_{A=1}^{n_cp}$ where $Q_{min}(\mathbf{x})$ represents the minimum bounding box containing $\mathbf{x}$. An example bounding box for a Loop subdivision basis function is illustrated in Figure~\ref{fig:bounding-box-closup}. Once low-rank approximations are computed for $\ary{H}$ and $\ary{G}$, we solve \eqref{eq:system-equations-implementation-2} using a GMRES iterative solver in combination with an approximation of the inverse operator computed through a triangular factorization.

We note that due to the relatively large support of the Loop subdivision basis there is a reduction in the number of admissible interactions over traditional discretisations, but this price is outweighed by the superior accuracy of the high-order basis.

\section{Numerical tests}

\subsection{Plane wave scattering problem: elastic spherical shell}

\begin{figure}[]
	\centering
	\includegraphics[width=0.8\textwidth]{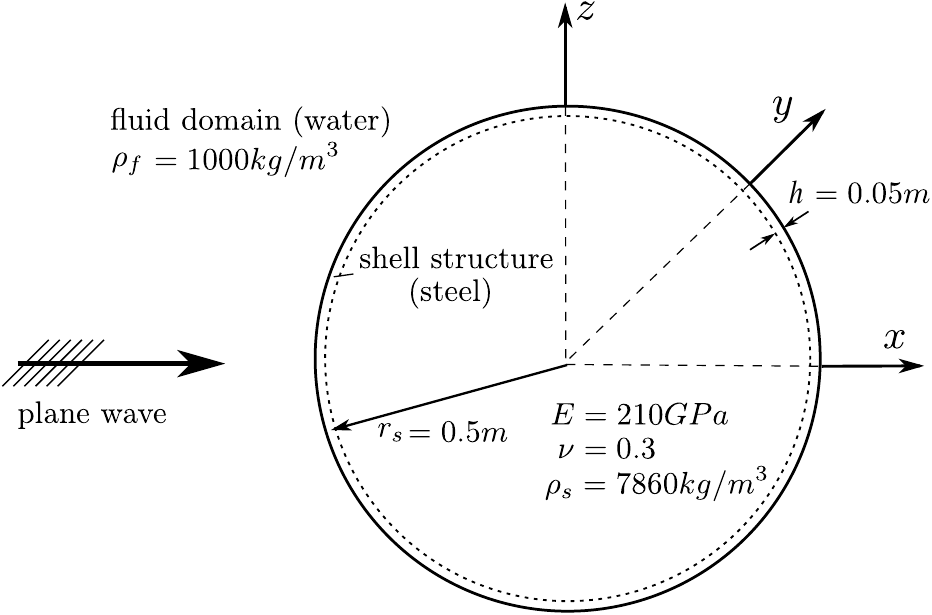}
	\caption{Coupled structural-acoustic problem of a plane wave impinged on a spherical shell immersed in an infinite fluid domain.}
	\label{fig:sphere-shell-coupled-problem}
\end{figure}

The problem of a plane wave impinged on an elastic spherical shell immersed in an infinite fluid  domain is illustrated in Figure~\ref{fig:sphere-shell-coupled-problem} with Table~\ref{tab:spherical-scattering-properties} specifying all geometry and material properties adopted in the present study.  The same material properties can be assumed in all included numerical examples unless specified otherwise.  We choose a plane wave of unit magnitude travelling in the positive $x$ direction given by $p_{inc} = e^{ikx}$ (i.e. $P=1$, $\mathbf{d}= (1,0,0)^T$). 
%\begin{figure}[]
%	\centering
%	\begin{subfigure}{0.8\textwidth}
%		\centering
%		\includegraphics[width=0.7\textwidth]{coupledshell.eps}
%		\caption{Geometry and material properties.}
%		\label{fig:sphere-shell-coupled-problem}
%	\end{subfigure}\\
%	\vspace{2ex}
%	\begin{subfigure}{0.8\textwidth}
%		\centering
%		\includegraphics[width=0.9\textwidth]{sample_points_shell.eps}
%		\caption{Location of sample points located at a radius $r=5m$ on the $x\mbox{-}y$ plane (drawing not to scale).  The number of sample points shown here is for illustration purposes only with the number of points used dependent on the wavenumber chosen.}
%		\label{fig:sphere-shell-sample-pts}
%	\end{subfigure}
%	\caption{Coupled structural-acoustic problem of a plane wave impinged on a spherical shell immersed in an infinite fluid domain.}
%\end{figure}

\begin{figure}[]
	\centering
	\includegraphics[width=0.6\textwidth]{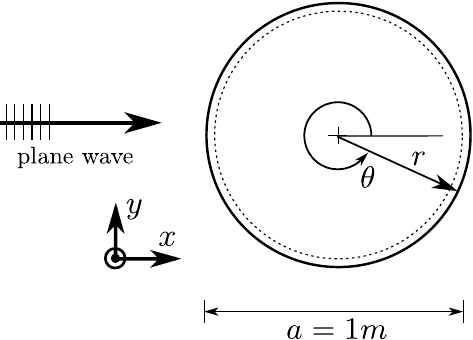}
	\caption{ Spherical shell scattering study: incident wave direction and polar coordinate system defined in the $x\mbox{-}y$ plane}
	\label{fig:polar-coordinate-system}
\end{figure}

\begin{table}[]
	\centering
	\caption{Material and geometry parameters adopted for the elastic sphere problem.}
	\label{tab:spherical-scattering-properties}
	\begin{tabular}{@{}cccc@{}}
		\toprule
		Parameter name      & Symbol  & Value & Unit                   \\\midrule
		density (water)    & $\rho_f$  & 1000  & $kg/m^3$ \\
		speed of sound(water)     & $c$       & 1482  & $m/s$                    \\
		density (steel)    & $\rho_s$    & 7860  & $kg/m^3$ \\
		Young's modulus    & $E$       & 210   & $GPa$                    \\
		Poisson's ratio    & $\nu$     & 0.3   & -                      \\
		sphere radius   & $r_s$       & 0.5   & $m$                      \\
		shell thickness & $h$       & 0.05  & $m$                    
		\\
		plane wave magnitude & $P$  & 1     & $Pa$
		\\ \bottomrule
	\end{tabular}
\end{table}

%The problem we consider to validate the coupling formulation of the Boundary Element Method and the Finite Element Method is an elastic sphere impinged by a plane wave under water. The problem is described as Fig. \ref{mypaper_fig:coup} shown. The spherical shell has a radius of $0.5 m$ and a uniform thickness of $0.05 m$. The material of the shell is steel, whose Young's modulus is $210 GPa$ and Poisson's ratio is $0.3$. The density of the steel is $7860 kg/m^3$. The spherical shell is enclosed by water, which has the density of $1000 kg/m^3$. The plane wave has a magnitude $A = 1$ with a direction $\mathbf{d} = \{1,0,0\}$.

The solution to this problem can be determined analytically~\cite{junger1986sound} which we reproduce here for completeness.  The total acoustic pressure $p_t \equiv p$ can be decomposed into scattered and elastic components as
\begin{equation}\label{eq:total-pressure-scat-ela}
p_{t} = p_{scat} + p_{ela}
\end{equation}
where $p_{scat}$ is the acoustic pressure that would result from scattering over a rigid sphere and $p_{ela}$ is the radiated acoustic pressure resulting from elastic shell vibrations.   Defining a polar coordinate system $(r,\theta)$ that lies in the $x\mbox{-}y$ plane as shown in Figure~\ref{fig:polar-coordinate-system}, the scattered and elastic pressure components can be expressed as
\begin{equation}\label{eq:scattered-pressure}
p_{scat}(r,\theta) = p_0\sum_{n=1}^{\infty}{-\frac{i^{n}(2n+1)j_n'(kr_s)}{h_n'(kr_s)}}P_n(\cos{\theta})h_n(kr)
\end{equation}
and
\begin{equation}\label{eq:elastic-pressure}
p_{ela}(r,\theta) = p_0\sum_{n=1}^{\infty}{\frac{i^{n}(2n+1)\rho_f c}{(Z_n+z_n)[kr_sh_n'(kr_s)]^2}}P_n(\cos{\theta})h_n(kr)
\end{equation}

where $p_0 \equiv P$ is the incident wave magnitude, $P_n$ is the $n$th Legendre function, $h_n$ and $h_n'$ are the $n$th Hankel function and its derivative respectively, and $j_n'$ is the derivative of the $n$th spherical Bessel function.  $Z_n$ denotes the \textit{invacuo} modal impedance of the spherical shell calculated as

\begin{equation}
Z_n = - \frac{i \rho_s c_p}{\Omega}\frac{h}{r_s}\frac{[\Omega^2 - (\Omega_n^{(1)})^2][\Omega^2 - (\Omega_n^{(2)})^2]}{[\Omega^2-(1+\beta^2)(\nu + \lambda_n -1)]}
\end{equation}

where $\lambda_n = n(n+1)$, $\Omega = \omega\frac{r_s}{c_p}$ is a dimensionless driving frequency, $\beta^2 = \frac{h^2}{12r_s^2}$, $c_p$ is the velocity of compressional waves in the structure given by 
\begin{equation}
c_p = \sqrt{\frac{E}{(1 - \nu^2)\rho_s}}
\label{eqn:c_p}
\end{equation}
and $\Omega_n^{(1)}$ and $\Omega_n^{(2)}$ are dimensionless natural frequencies of the spherical shell determined from the two positive roots of the polynomial equation

\begin{equation}
\begin{split}
\Omega_n^4 - [1+3\nu+\lambda_n - \beta^2(1 - \nu - \lambda_n^2 - \nu\lambda_n)]\Omega_n^2  + (\lambda_n - 2)(1 - \nu^2)
+\\ \beta^2[\lambda_n^3 - 4\lambda_n^2 + \lambda_n(5 - \nu^2) - 2(1 - \nu^2)] = 0.
\end{split}
\label{eqn:nondimensial_frequencies}
\end{equation}

Finally, $z_n$ is the modal specific acoustic impedance expressed as

\begin{equation}
z_n = i\rho_f c\frac{h_n(kr_s)}{h'_n(kr_s)}.
\end{equation}

\subsubsection{Geometrical error study}
\label{sec:sphere-geometry-study}

\begin{figure}[h]
\centering
\begin{subfigure}{0.49\textwidth}
\centering
  \includegraphics[height=5.5cm]{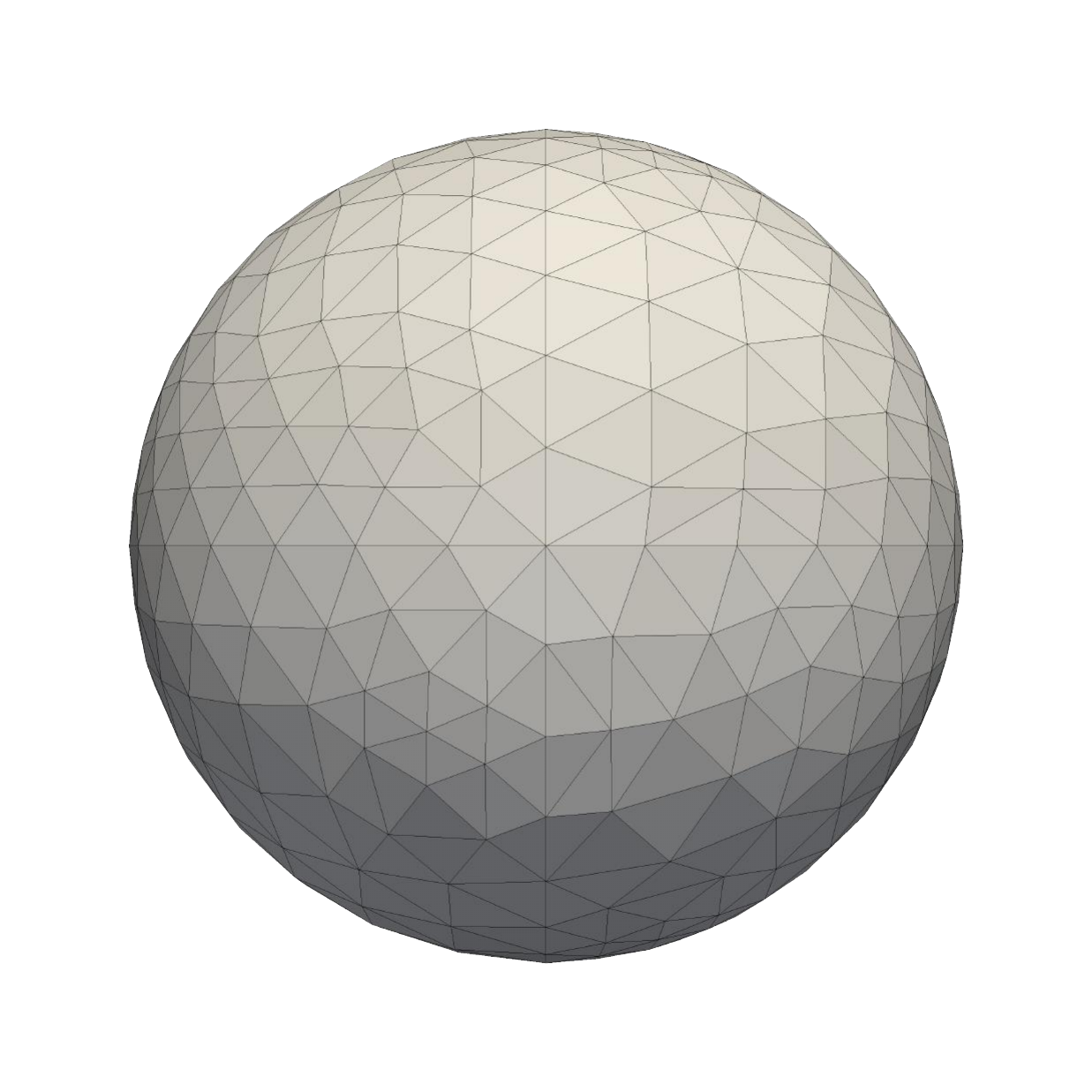}
  \caption{Control mesh.}
  \label{fig:initial-sphere-mesh-grid}
\end{subfigure}
\begin{subfigure}{0.49\textwidth}
\centering
  \includegraphics[height=5.5cm]{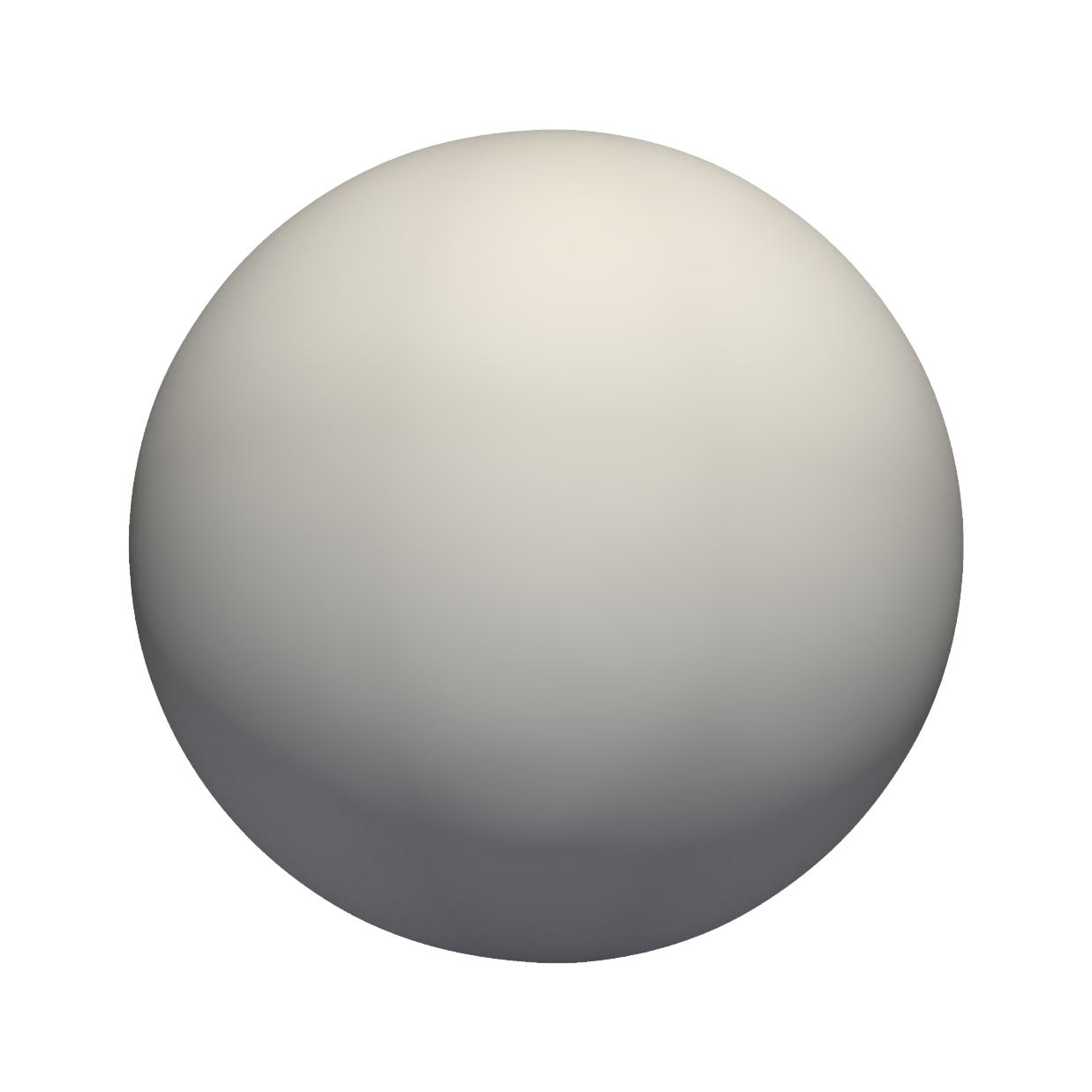}
  \caption{Limit surface.}
  \label{fig:initial-sphere-mesh-surface}
\end{subfigure}
\caption{The initial coarse Loop subdivision discretisation with 438 vertices used to generate control meshes (a) through to (d) in Table~\ref{tab:sphere-meshes-subd} and Table~\ref{tab:sphere-meshes-proj}.The control vertices are placed such that they lie on a sphere with diameter $a=1$.}
\label{fig:initial-sphere-mesh}
\end{figure}

The first numerical study we conduct is to verify convergence characteristics of the present method to the analytical solution given by Equations~\eqref{eq:total-pressure-scat-ela} to \eqref{eq:elastic-pressure}.  We construct the system of equations given by~\eqref{eq:system-equations-implementation-2} using our Loop subdivision discretisation procedure and choose a relatively low normalised wavenumber of $ka=10$, where $a$ is the diameter of the spherical shell.  Four Loop subdivision control girds are generated from an initial control mesh shown in Figure~\ref{fig:initial-sphere-mesh-grid} using two strategies:
\begin{enumerate}
	\item \textbf{Subdivision}: the Loop subdivision refinement algorithm of~\cite{loop1987smooth} is applied to the initial coarse control mesh to generate two successively refined control meshes (a) and (b).  The limit surface of the coarse control mesh and meshes given in Table~\ref{tab:sphere-meshes-subd}~(a) and~(b) is  identical (see Figure~\ref{fig:initial-sphere-mesh-surface}) but exhibits a non-negligible geometry errors due to its deviation from a sphere. For analysis evidently the approximation basis provided by control mesh (b) provides a richer approximation space over (a).
	\item \textbf{Least squares fitting}: two successively refined control meshes (c) and (d) are generated by performing a least square fitting, or $L_2$ projection, of the subdivision surface to a sphere with diameter $a = 1$. With this refinement strategy the geometry error is successively reduced but the same approximation basis as for the equivalent control meshes (a) and (b) is generated.
\end{enumerate}

\begin{figure}[h]
\centering
\begin{subfigure}{0.49\textwidth}
\centering
  \includegraphics[height=5.5cm]{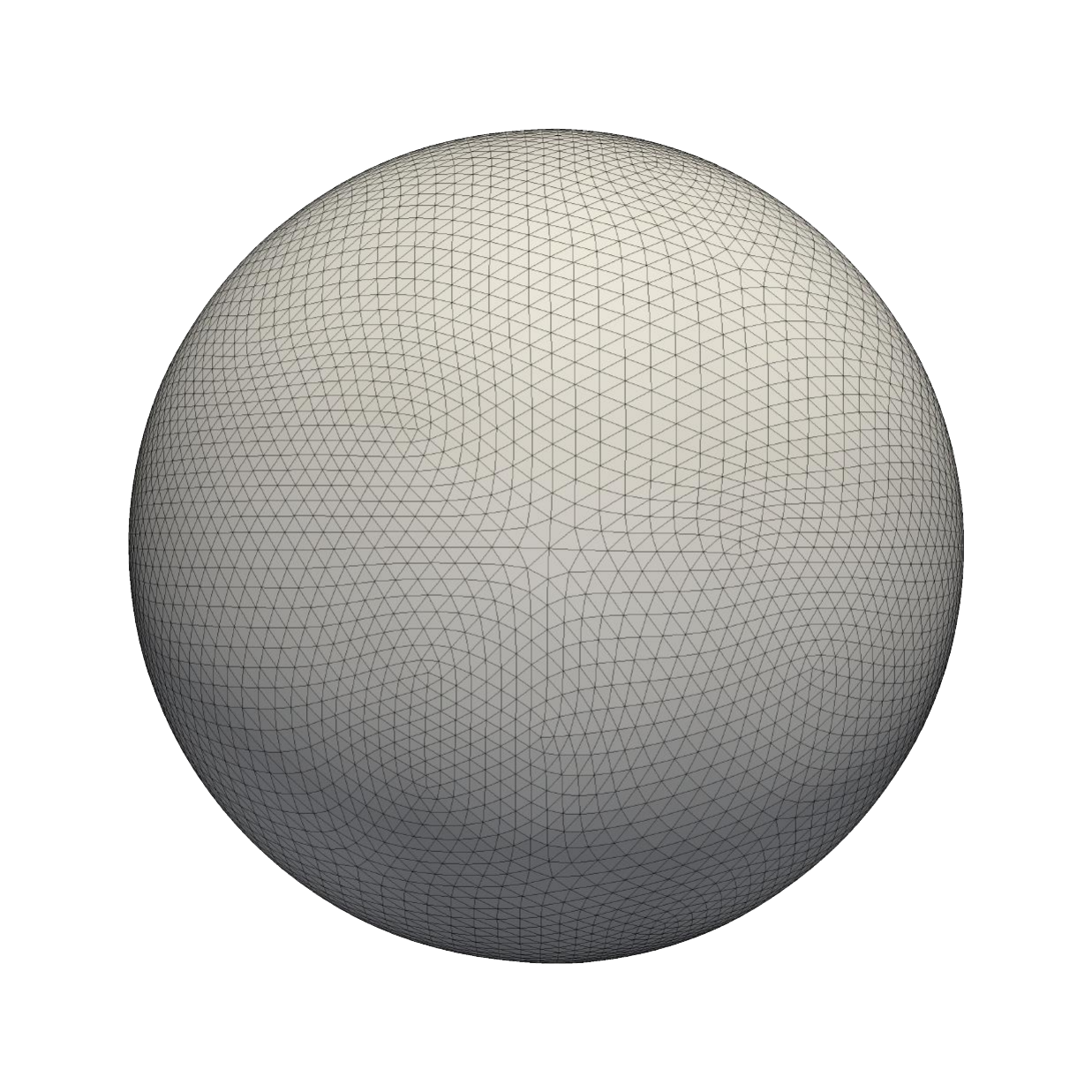}
  \caption{Control mesh.}
  \label{fig:sphere-mesh-d-grid}
\end{subfigure}
\begin{subfigure}{0.49\textwidth}
\centering
  \includegraphics[height=5.5cm]{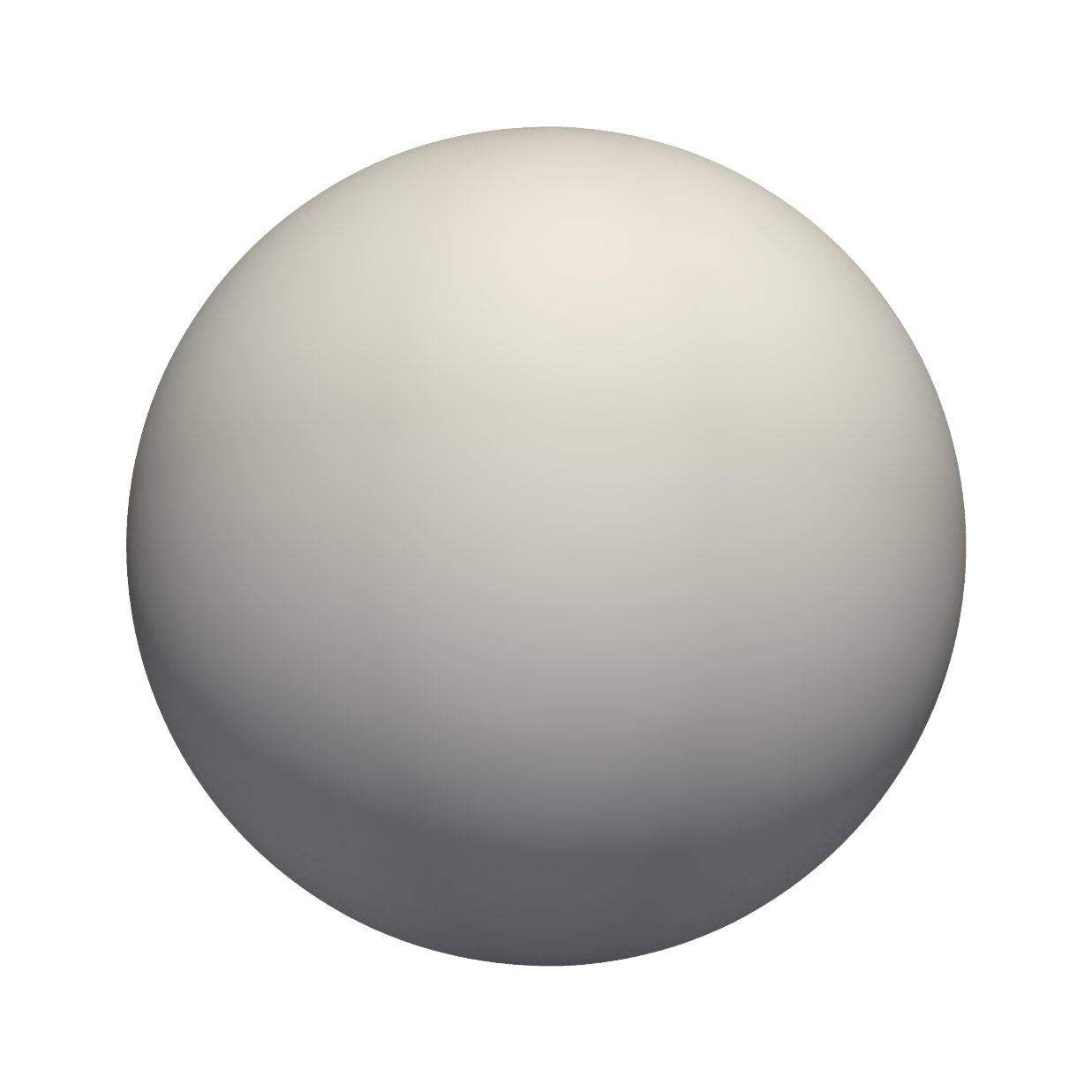}
  \caption{Limit surface.}
  \label{fig:sphere-mesh-d-surface}
\end{subfigure}
\caption{The Loop subdivision discretisation of a sphere geometry with 6978 vertices corresponding to control mesh (d) in Table~\ref{tab:sphere-meshes-proj}. The control vertex positions are determined such that the limit surface approximates a sphere with diameter $a=1$.}
\label{fig:sphere-mesh-d}
\end{figure}

For each control mesh we calculate the relative geometry error $\varepsilon_{g}$ of the limit surface as
\begin{equation}
\varepsilon_{g} = \frac{|| \mathbf{x}^h - \mathbf{x}||_0}{||\mathbf{x}||_0}
\end{equation}
where $\mathbf{x}^h$ and $\mathbf{x}$ are physical coordinates on the Loop subdivision surface and analytical surface respectively with
\begin{equation}
||\mathbf{\cdot}||_0 := \left(\int_{\Gamma} (\mathbf{\cdot})^2\, \mathrm{d}\Gamma \right)^{1/2}.
\end{equation}
Figure~\ref{fig:initial-sphere-mesh} illustrates the initial control mesh with its associated limit surface. Similar illustrations are shown in Figure~\ref{fig:sphere-mesh-d} for control mesh (d).  Table~\ref{tab:sphere-meshes-subd} and \ref{tab:sphere-meshes-proj} details each of the control meshes for both refinement strategies showing that geometry error remains constant during subdivision refinement (the limit surface is independent of subdivision refinement) and converges to zero when control vertices are $L_2$ projected onto the analytical sphere surface.
%\begin{table}[h]
%	\centering
%	\caption{Mesh properties}
%	\label{tab:sphere-meshes}
%	\begin{tabular}{@{}cccccc@{}}
%		\toprule
%		control grid               & initial   & (a)                & (b)               & (c)        & (d)         \\ \midrule
%		refinement method      & -         & \textbf{subdivision} & \textbf{subdivision} & \textbf{projection}       &   \textbf{projection}       \\
%		degrees of freedom & 438      & 1746               & 6978               & 1746       & 6978        \\
%		$n_e$ & 872       & 3488    & 13952    & 
%		3488               &13952\\
%		$\varepsilon_{g}$     & 0.93\% & 0.93\%          & 0.93\%   & 0.24\% & 0.06\%         \\ \bottomrule
%	\end{tabular}
%\end{table}

\begin{table}[h]
	\centering
	\caption{Mesh properties for subdivision refinement}
	\label{tab:sphere-meshes-subd}
	\begin{tabular}{@{}cccccc@{}}
		\toprule
		control mesh               & initial   & (a)                & (b)  \\  \hline
		number of vertices & 438      & 1746               & 6978                      \\
		$n_e$ & 872       & 3488    & 13952    \\
		$\varepsilon_{g}$     & 0.93\% & 0.93\%          & 0.93\%   \\ \bottomrule
	\end{tabular}
\end{table}

\begin{table}[h]
	\centering
	\caption{Mesh properties for subdivision refinement and subsequent $L_2$ projection}
	\label{tab:sphere-meshes-proj}
	\begin{tabular}{@{}cccccc@{}}
		\toprule
		control mesh               & initial   & (c)        & (d)         \\  \hline
		number of vertices & 438       & 1746       & 6978        \\
		$n_e$ & 872       & 3488    & 13952    \\
		$\varepsilon_{g}$     & 0.93\%   & 0.24\% & 0.06\%         \\ \bottomrule
	\end{tabular}
\end{table}

\begin{figure}[h]
	\centering
	\includegraphics[width=14cm]{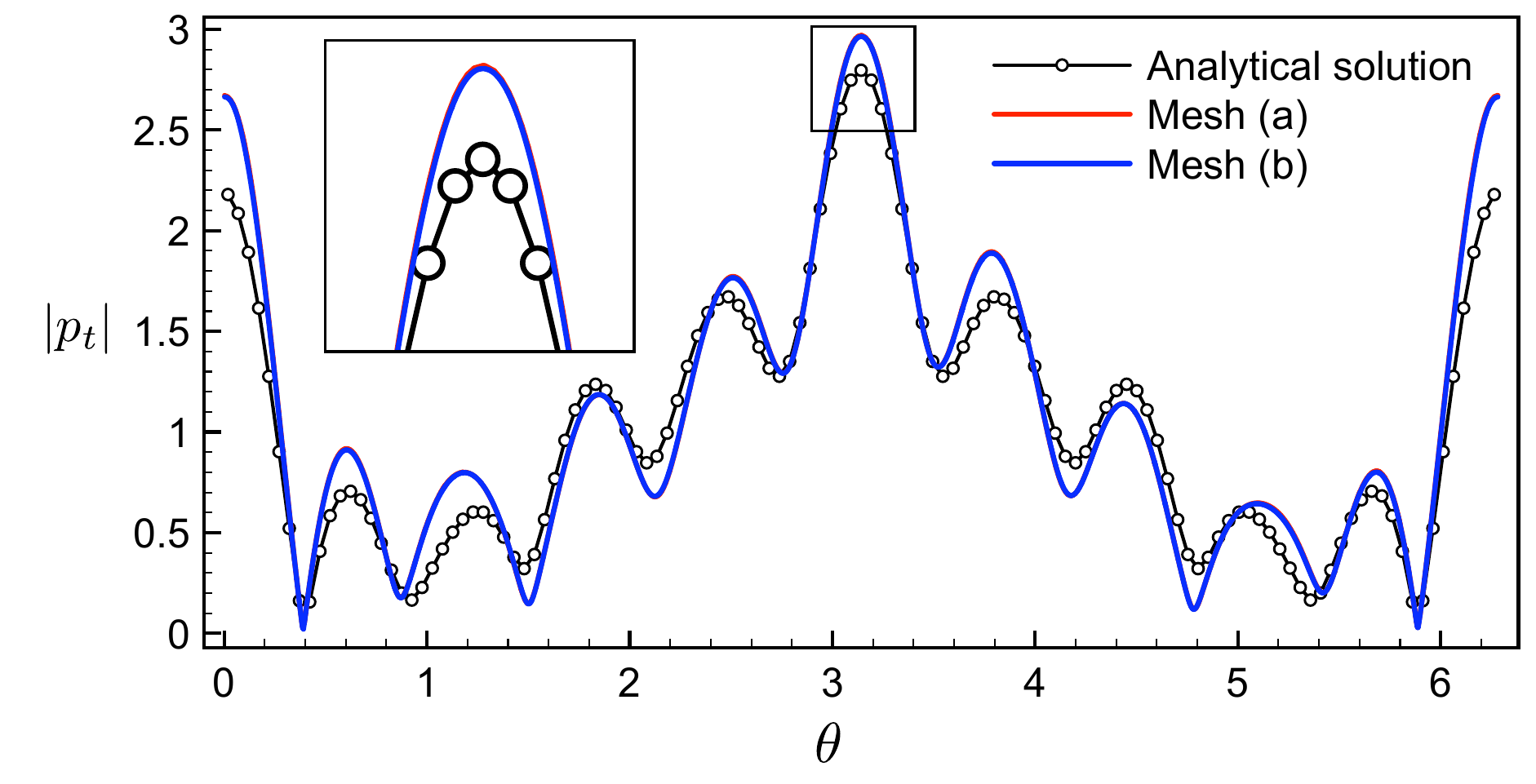}
	\caption{Coupled sphere shell problem,  $ka = 10$: surface acoustic potential magnitude along $x\mbox{-}y$ plane using control meshes (a) and (b). The inset image illustrates convergence to a solution which does not correspond to the analytical solution due to the non-negligible geometrical error of the limit surface.}
	\label{fig:sphere-meshab-results}
\end{figure}

\begin{figure}[h]
	\centering
	\includegraphics[width=14cm]{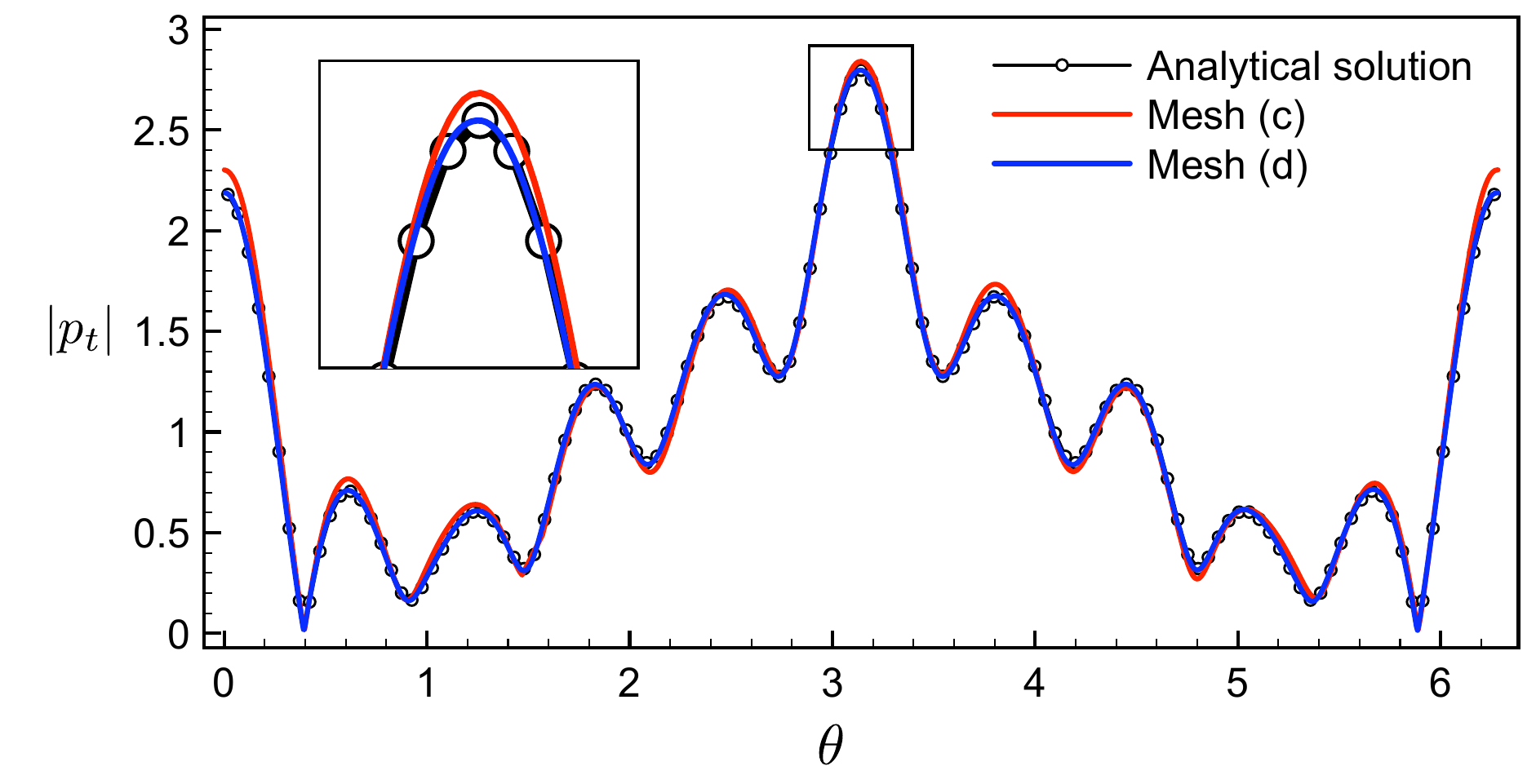}
	\caption{Coupled sphere shell problem,  $ka = 10$: surface acoustic potential magnitude along $x\mbox{-}y$ plane using control meshes (c) and (d). $L_2$ projection of control vertices onto the analytical sphere surface leads to a reduced geometrical error in the limit surface and convergence to the analytical solution.}
	\label{fig:sphere-meshcd-results}
\end{figure}

% *******************************
% ****** Sphere results *********
% *******************************
We compute the magnitude of the complex-valued total acoustic pressure $|p_t| = \left( \textrm{Re}(p_t)^2 + \textrm{Im}(p_t)^2 \right)^{1/2}$ at sample points located  on the $x\mbox{-}y$ plane of the sphere surface. Results for control meshes (a) and (b) are shown in Figure~\ref{fig:sphere-meshab-results} which illustrate convergence to a solution that is distinct from the analytical solution.  In contrast, the results for control meshes (c) and (d) (Figure~\ref{fig:sphere-meshcd-results}) illustrate convergence to the analytical solution and thus demonstrate the importance of controlling geometry error in the context of coupled structural-acoustic problems.  The relatively small geometry error of 0.93\% leads to sample point solution errors in the order of 20\% and therefore care must be taken to ensure that the limit surface provides an accurate representation of required model geometry\textemdash further application of the Loop subdivision refinement algorithm will not overcome the inherent error induced by the incorrect geometry representation. In the case that the limit surface is an accurate representation of the sphere our method converges to the analytical coupled solution and thus verifies our implementation. We remark that for scenarios in which an exact sphere representation is required there exist special subdivision schemes but we do not consider these in the present study.

% Results for meshes (c) and (d)

\clearpage

\subsubsection{Medium frequency problems}
We now consider the ability of the present method to handle medium frequency problems and determine the upper frequency limits of our Loop subdivision discretisation.  It is well-known that when traditional low-order discretisations are used for wave problems involving medium to high frequencies high resolution meshes must be used in order to control approximation error.  A common rule of thumb is to use ten elements per wavelength that can lead to extremely large systems of equations in the case of complex geometries and high frequencies and it is therefore desirable to use a higher order basis that reduces the number of elements per wavelength and associated system of equations.

We use the Loop subdivision discretisations generated through control meshes (c) and (d) shown in Section~\ref{sec:sphere-geometry-study} (see Table~\ref{tab:sphere-meshes-proj} for details) and determine the upper frequency limit of each discretisation by applying a set of increasing normalised wavenumbers $ka=10,30,40,50,60,80$.  Defining a set of sample points on the sphere surface aligned with the $x\mbox{-}y$ plane as $\mathsf{S} = \{s_1, s_2, \ldots, s_{n_{sample}}\}$ we introduce the maximum pointwise error in the sample set as

\begin{equation}\label{eq:max-pointwise-error}
\max_{\mathbf{s} \in \mathsf{S}} \frac{\left| |p_t^h(\mathbf{s})| - |p_t(\mathbf{s})| \right|}{{\left| | p_t(\mathbf{s})| \right|}_{\infty}}
\end{equation}

where $p_t^h$ and $p_t$ represent numerical and analytical total acoustic pressures respectively.

\begin{figure}[h]
	\centering
	\includegraphics[width=14cm]{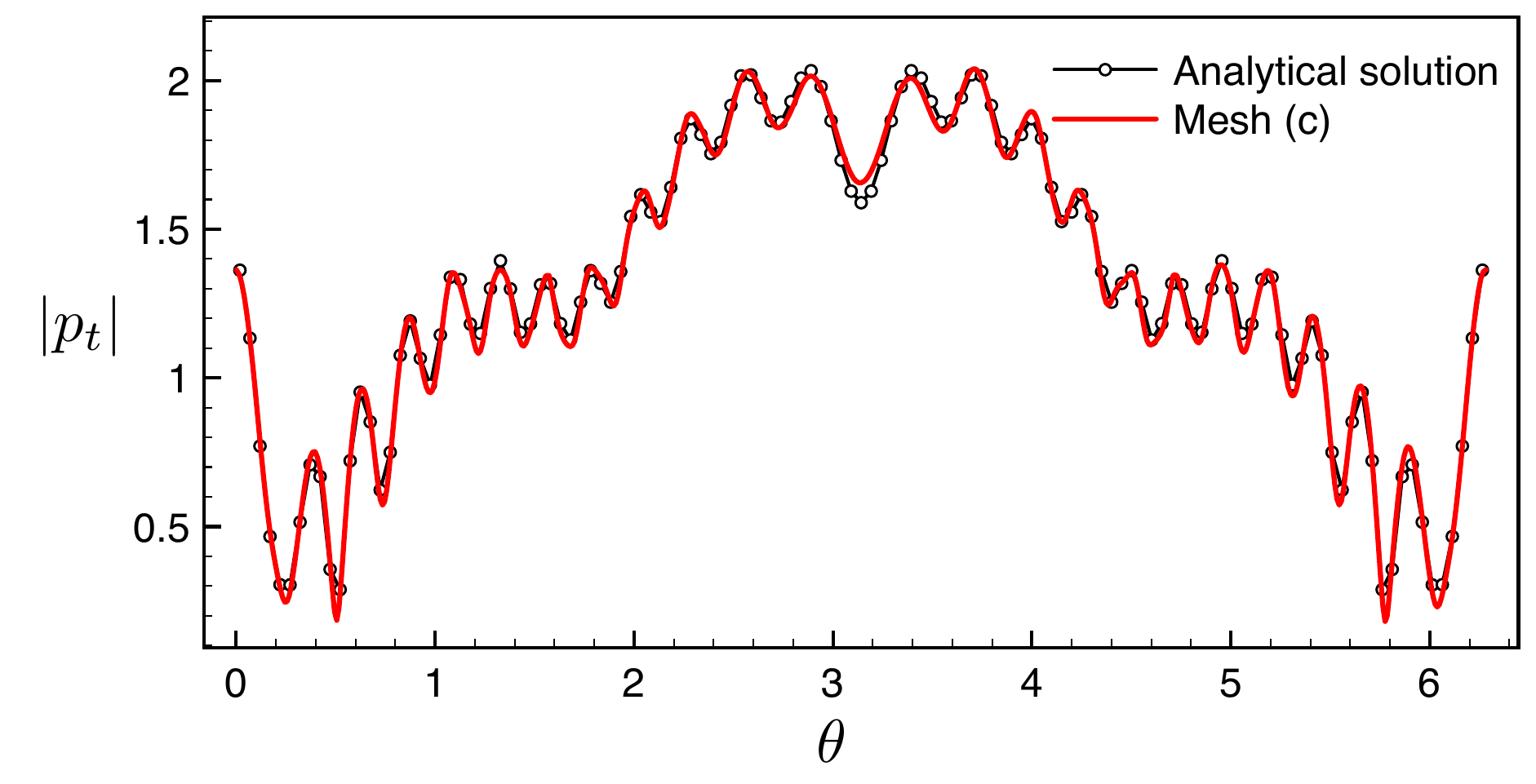}
	\caption{Coupled sphere study: surface acoustic potential magnitude profile along $x\mbox{-}y$ plane for $ka=30$, control mesh (c).}
	\label{fig:sphere-mesh-c-k30}
\end{figure}

\begin{figure}[h]
	\centering
	\includegraphics[width=14cm]{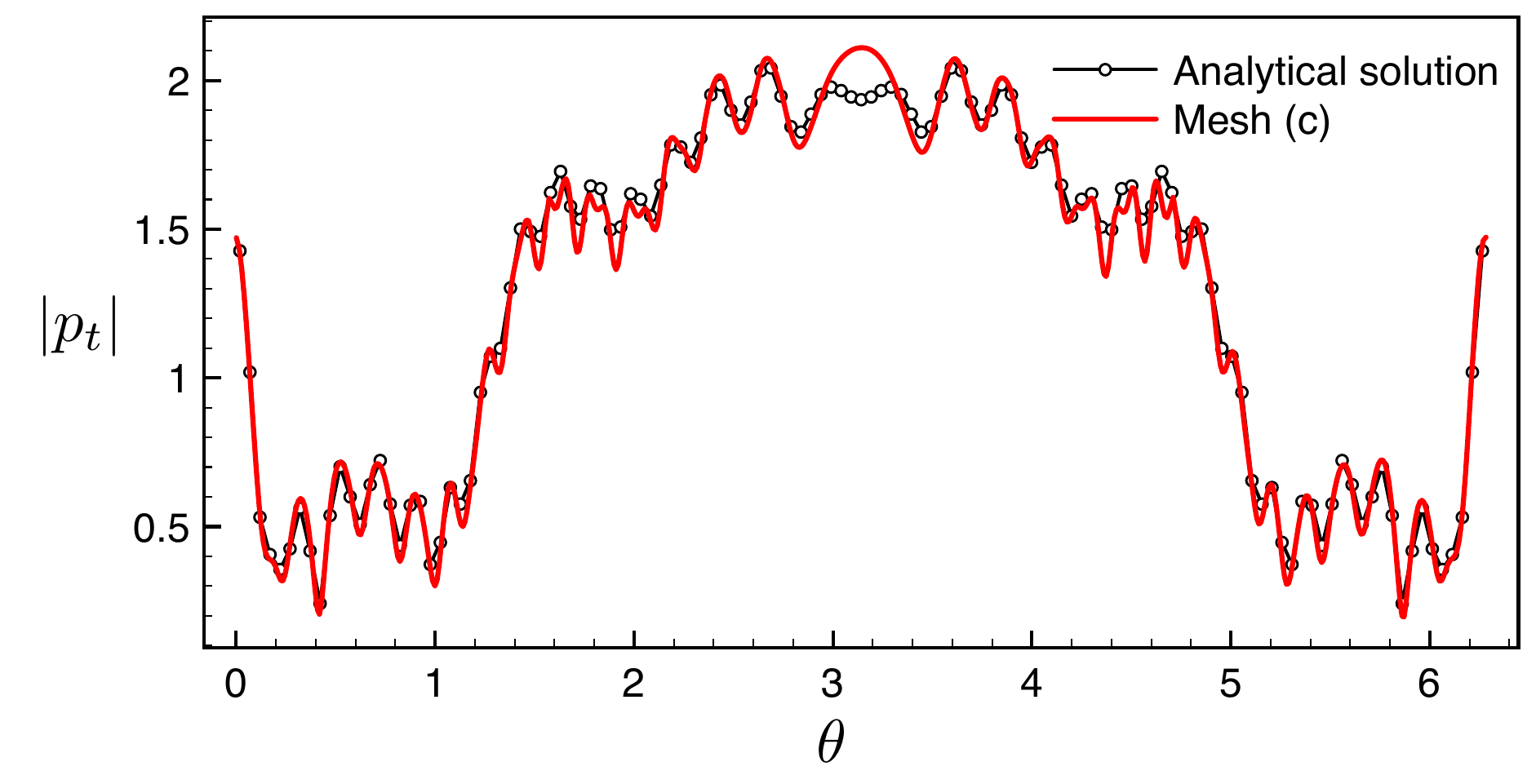}
	\caption{Coupled sphere study: surface acoustic potential magnitude along $x\mbox{-}y$ plane for $ka=40$, control mesh (c).}
	\label{fig:sphere-mesh-c-k40}
\end{figure}

\begin{figure}[h]
	\centering
	\includegraphics[width=14cm]{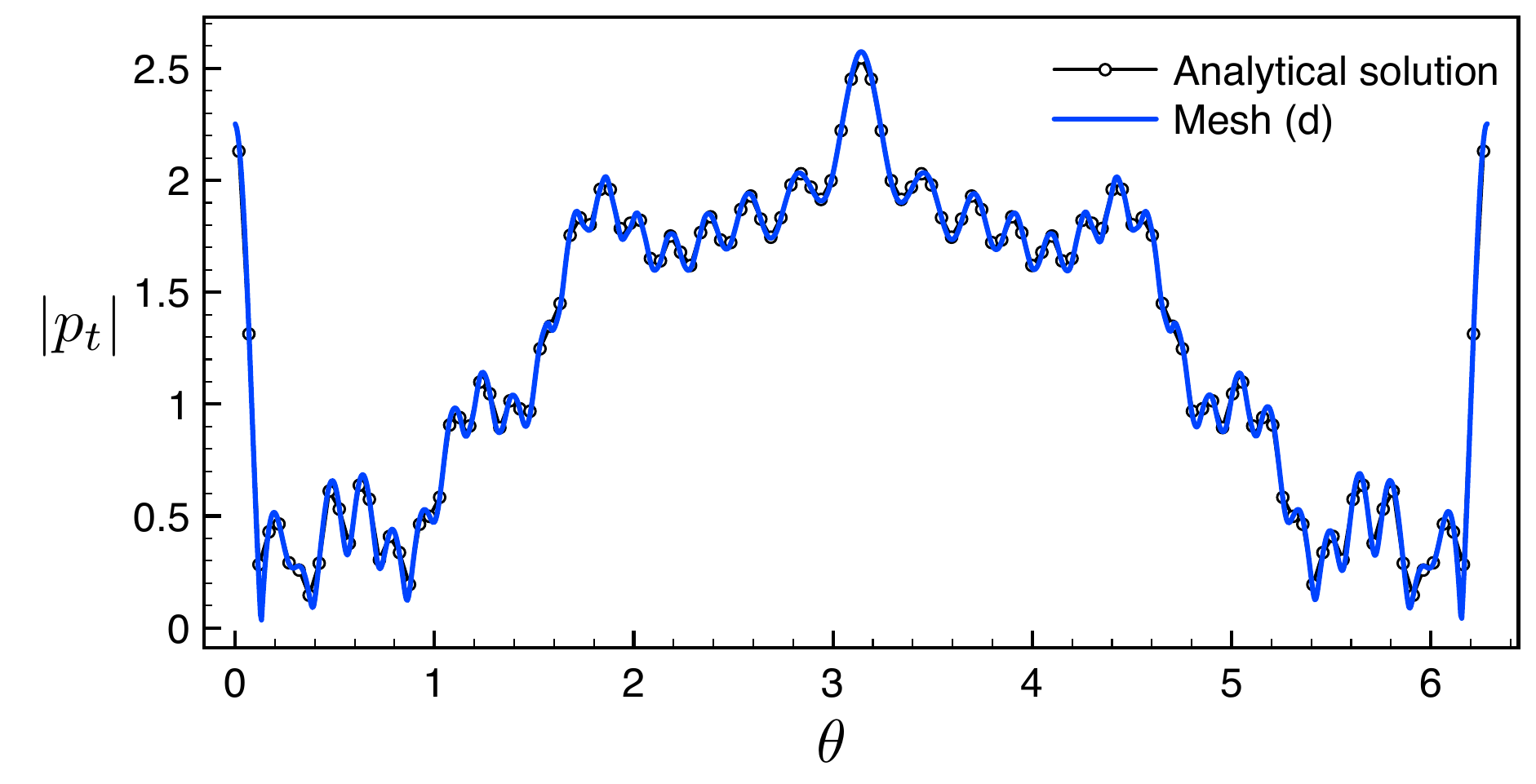}
	\caption{Coupled sphere study: surface acoustic potential magnitude along $x\mbox{-}y$ plane for $ka=50$, control mesh (d).}
	\label{fig:sphere-mesh-d-k50}
\end{figure}

\begin{figure}[h]
	\includegraphics[width=14cm]{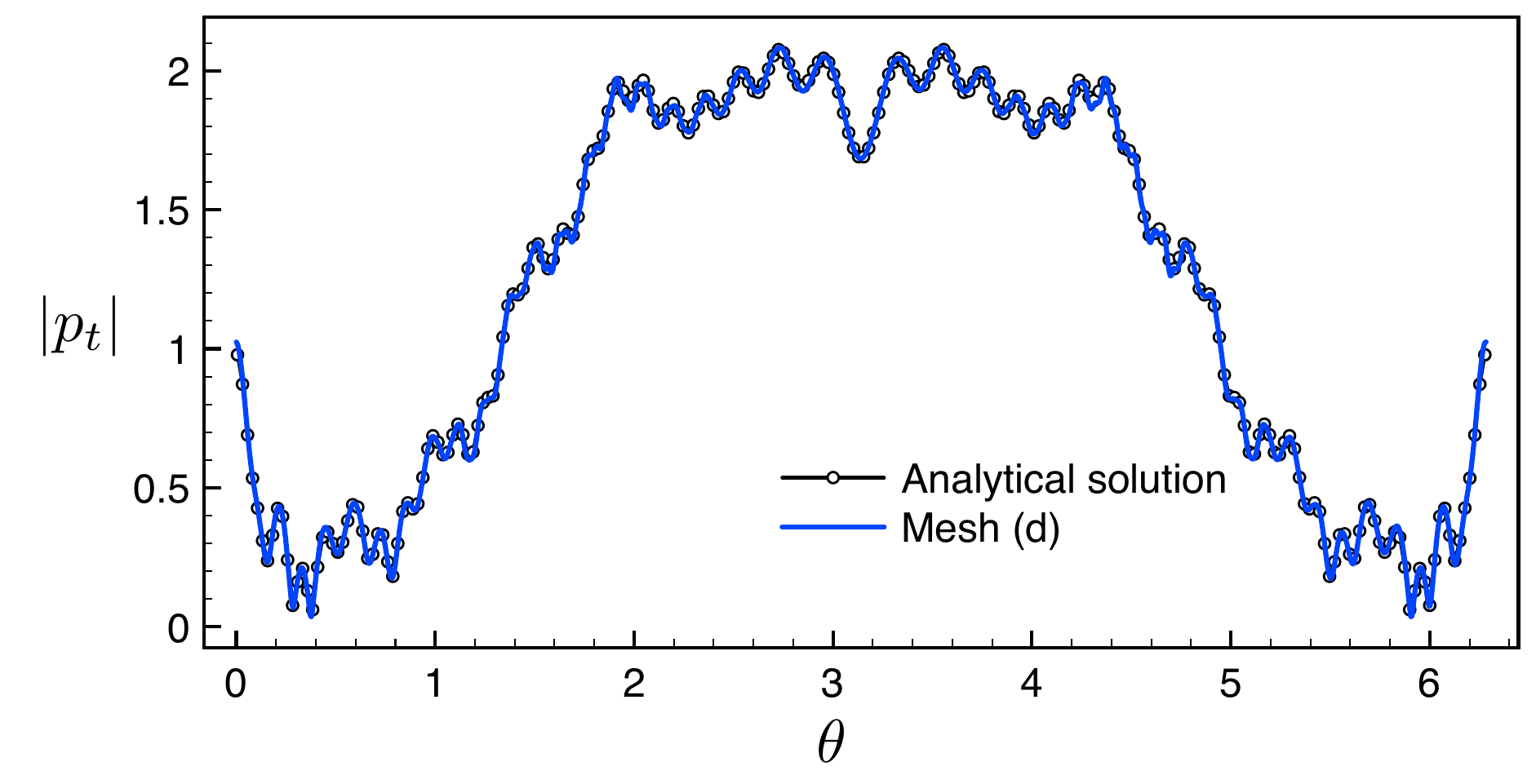}
	\caption{Coupled sphere study: surface acoustic potential magnitude along $x\mbox{-}y$ plane for $ka=60$, control mesh (d).}
	\label{fig:sphere-mesh-d-k60}
\end{figure}

\begin{figure}[h]
	\centering
	\includegraphics[width=14cm]{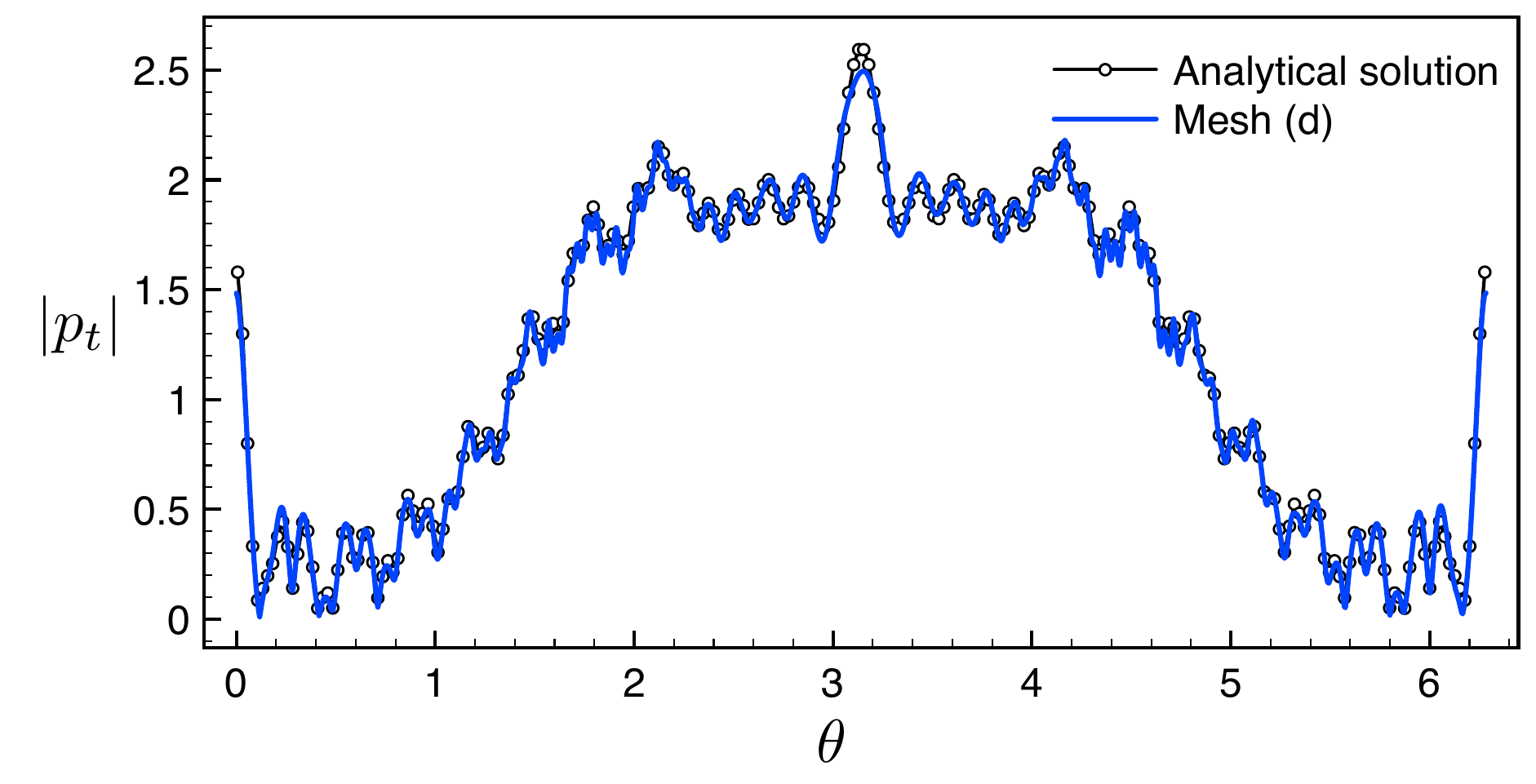}
	\caption{Coupled sphere study: surface acoustic potential magnitude along $x\mbox{-}y$ plane for $ka=80$, control mesh (d).}
	\label{fig:sphere-mesh-d-k80}
\end{figure}

\begin{figure}[h]
	\centering
	\includegraphics[width=0.5\textwidth]{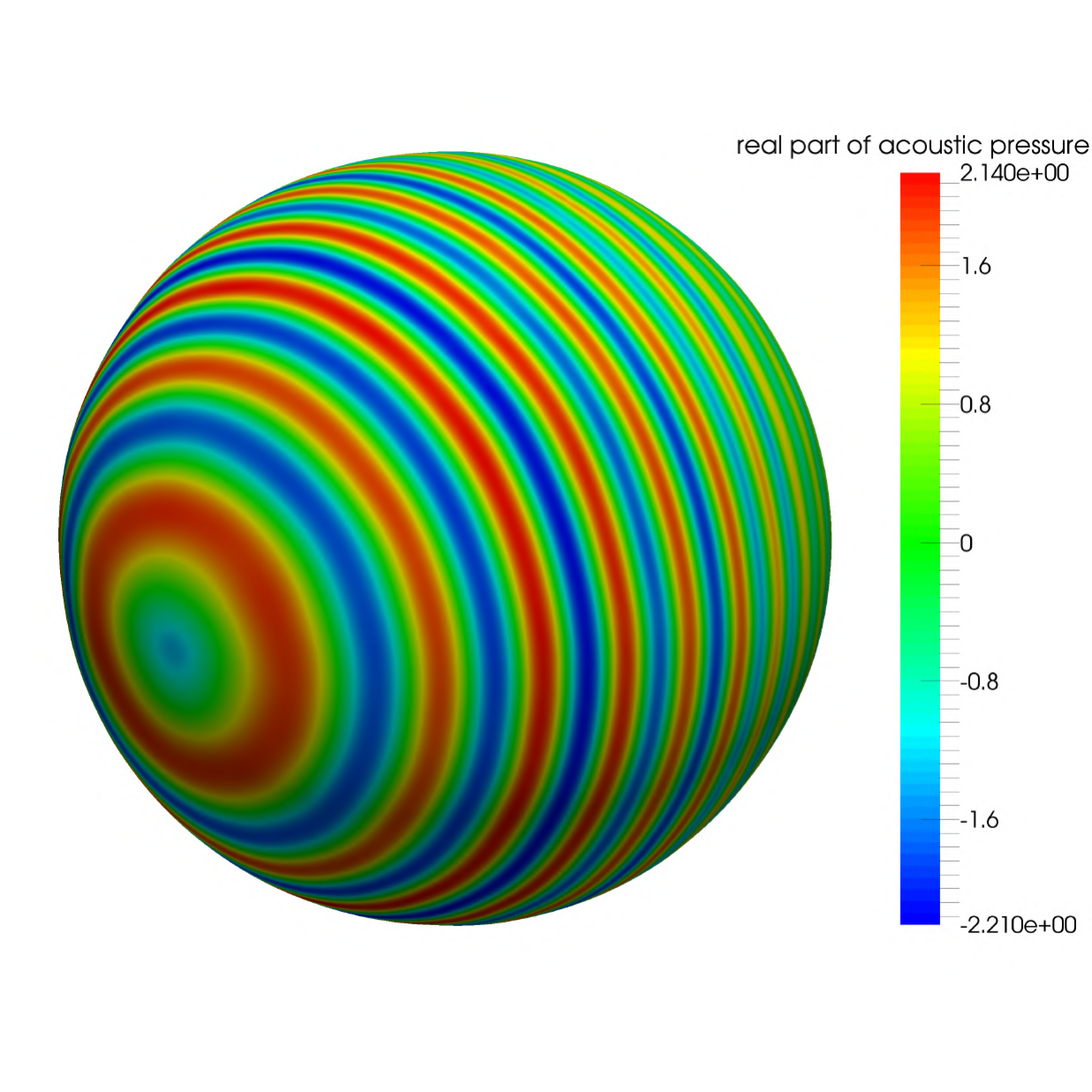}
	\caption{Coupled sphere study: $\textrm{Re}(p)$ (real component of acoustic pressure) for $ka=80$, control mesh (d).}
	\label{fig:sphere-real-pressure-surfaceplot-k80}
\end{figure}

Pointwise errors calculated through \eqref{eq:max-pointwise-error} for discretisations generated through control meshes (c) and (d) are tabulated in Table~\ref{tab:sphere-pointwise-errors} for each wavenumber.  The approximate number of elements per wavelength is detailed for each case.  Plots of $|p_t|$ against analytical solutions for $ka=30,40$ with control mesh (c) are shown in Figures~\ref{fig:sphere-mesh-c-k30} and \ref{fig:sphere-mesh-c-k40} respectively (Figure~\ref{fig:sphere-meshcd-results} illustrates results for $ka=10$) and similar plots are also shown for $ka=50,60,80$ for control mesh (d) in Figures~\ref{fig:sphere-mesh-d-k50}, \ref{fig:sphere-mesh-d-k60} and \ref{fig:sphere-mesh-d-k80} respectively.   Figure~\ref{fig:sphere-real-pressure-surfaceplot-k80} illustrates the real part of the total acoustic pressure on the surface on the sphere for $ka=80$. 

% ******************************************************
% ****** medium frequency pointwise errors *************
% ******************************************************
\begin{table}[h]
	\centering
	\caption{Maximum pointwise errors calculated through Equation~\eqref{eq:max-pointwise-error} for a set of increasing normalised wavenumbers applied to the coupled sphere problem. The numbers in parentheses indicates the approximate number of elements per wavelength. A dash indicates a result with a large maximum pointwise error that indicates insufficient resolution of mesh for the given wavenumber.}
	\label{tab:sphere-pointwise-errors}
	\begin{tabular}{ccccccc}
		\hline
		\multirow{2}{*}{\textbf{\begin{tabular}[c]{@{}c@{}}control \\ mesh\end{tabular}}} & \multicolumn{6}{c}{$ka$}      \\ \cline{2-7} 
		& 10 & 30 & 40 & 50 & 60 & 80 \\ \hline 
		(c) & 0.0427 (16)   & 0.0354 (6)   & 0.0853 (4)  & -   & -   &  -  \\
		(d) & 0.0090 (32)    & 0.0182 (11)   & 0.0194 (8)    & 0.0126 (7)  & 0.0208 (6)  & 0.0553  (4) \\ \hline
	\end{tabular}

\end{table}

% *******************************************
% ****** medium frequency plots *************
% *******************************************

An inspection of total acoustic pressures profiles and maximum pointwise  errors for control mesh (c) reveals the expected general trend of increasing errors for increasing wavenumbers.  Figure~\ref{fig:sphere-mesh-c-k30} indicates a discrepancy between the analytical and numerical potential magnitudes at $\theta=\pi$ which is in fact caused by geometry error (further subdivision refinement of control mesh (c) indicates that the solution has converged).  For $ka=40$ large errors are encountered ( $\sim8.5\%$ maximum pointwise error) and the discretisation is deemed insufficiently fine.  The more refined basis generated through control mesh (d) allows for much higher wavenumbers and it is found that for wavenumbers up to $ka=60$ low errors are obtained with maximum pointwise errors less than $2.1\%$. For $ka=80$ we find that the maximum wavenumber for control mesh (d) is reached but remark that even for this relatively high wavenumber a maximum pointwise error of $\sim 5.5\%$ is seen. Based on the present results we advise  a conservative estimate of six elements per wavelength using a Loop subdivision discretisation to obtain maximum pointwise errors less than $5\%$.  The guidance of six elements per wavelength is in agreement with the study of \cite{marburg2002six}.

\clearpage
% **************************************
% ****** Coupling effect study *********
% **************************************
\subsubsection{Coupling effect study}

\begin{figure}[h]
	\centering
	\begin{subfigure}[b]{0.45\textwidth}
		\includegraphics[width=\textwidth]{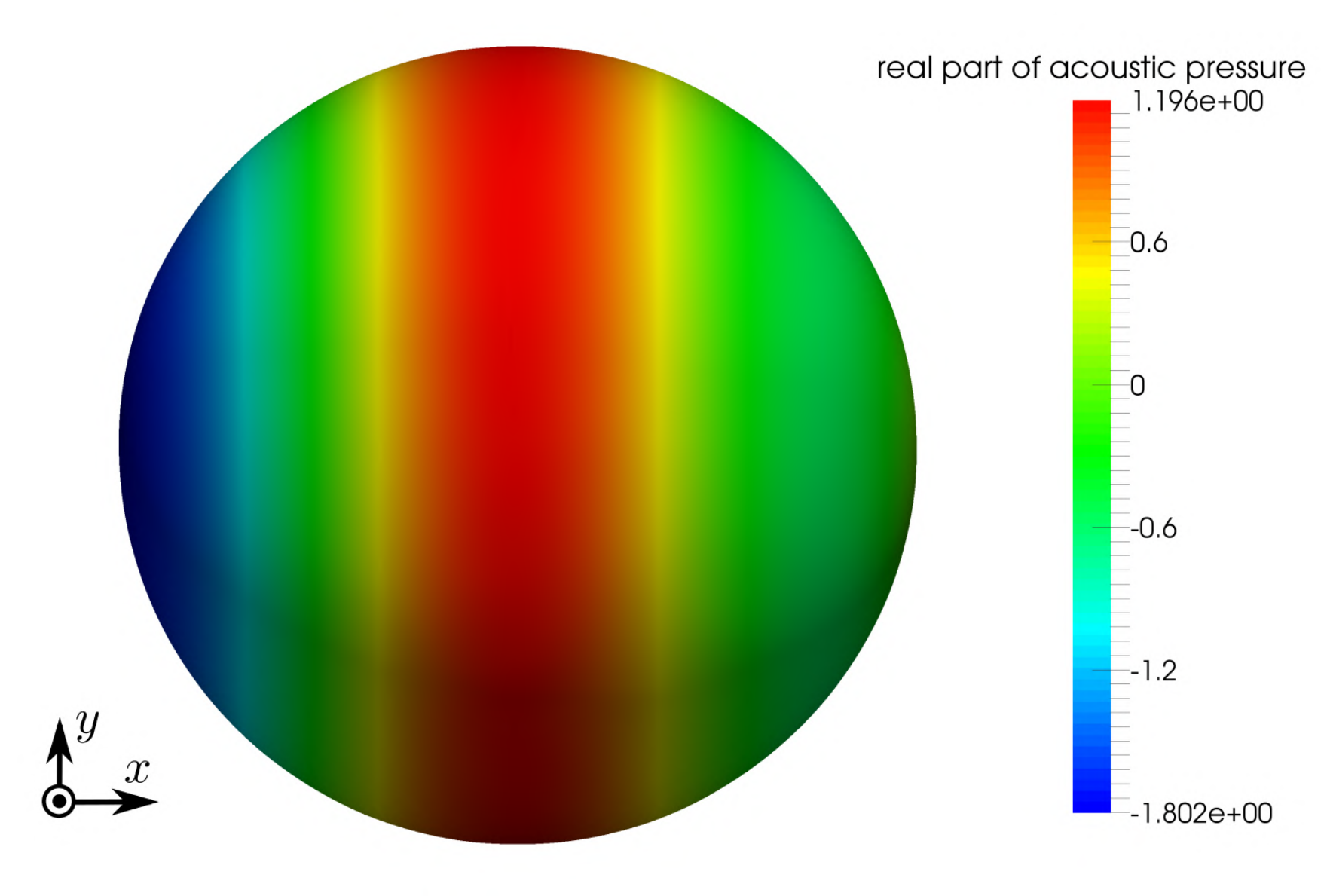}
		\caption{$\textrm{Re}(p)$, $x\mbox{-}y$ plane.}
		\label{hard-sphere-real-acoustic-potential-xy}
	\end{subfigure}
	\begin{subfigure}[b]{0.45\textwidth}
		\includegraphics[width=\textwidth]{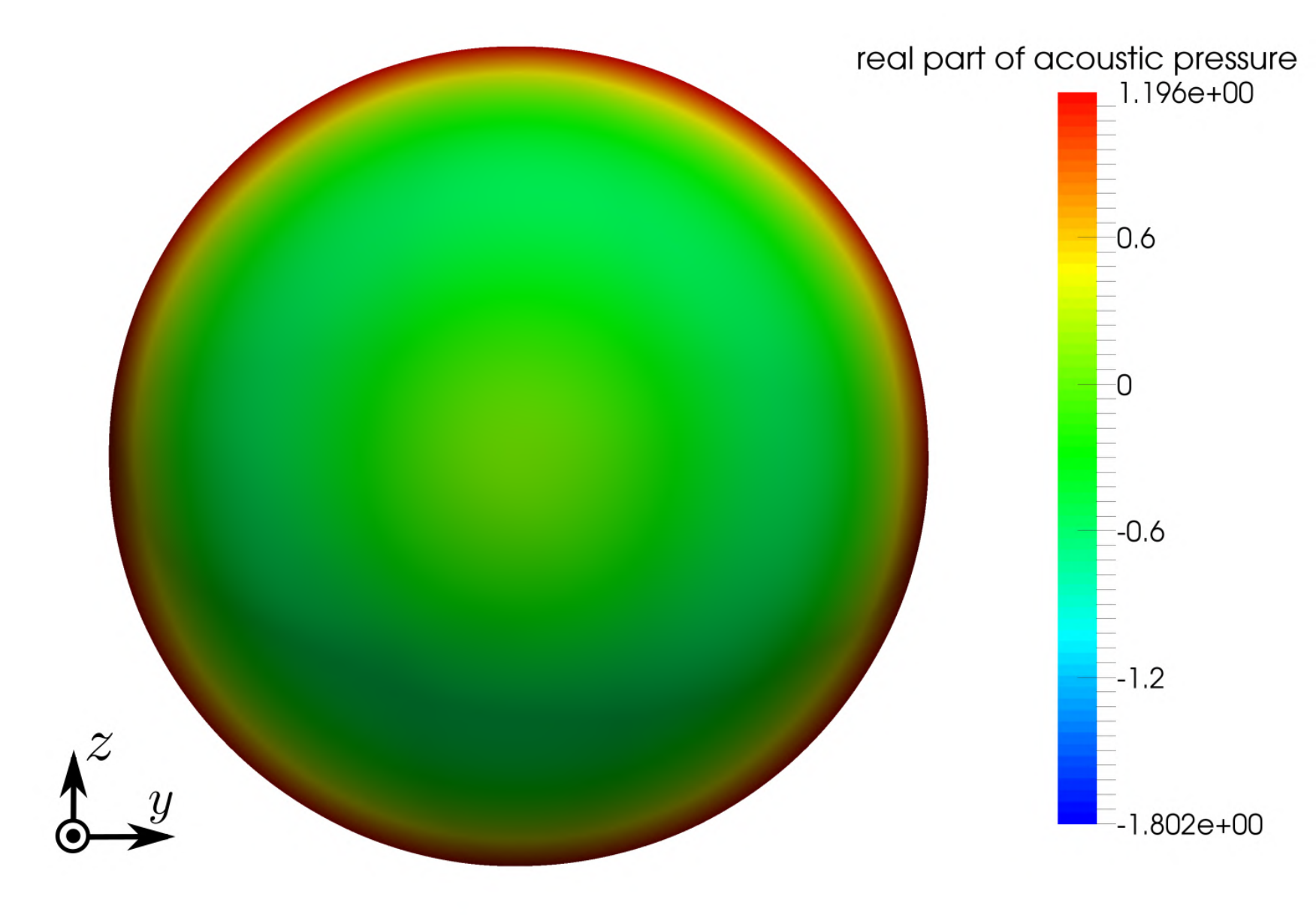}
		\caption{$\textrm{Re}(p)$, $y\mbox{-}z$ plane.}
		\label{hard-sphere-real-acoustic-potential-yz}
	\end{subfigure}
	
	\vspace{2ex}

	\begin{subfigure}[b]{0.6\textwidth}
	\includegraphics[width=\textwidth]{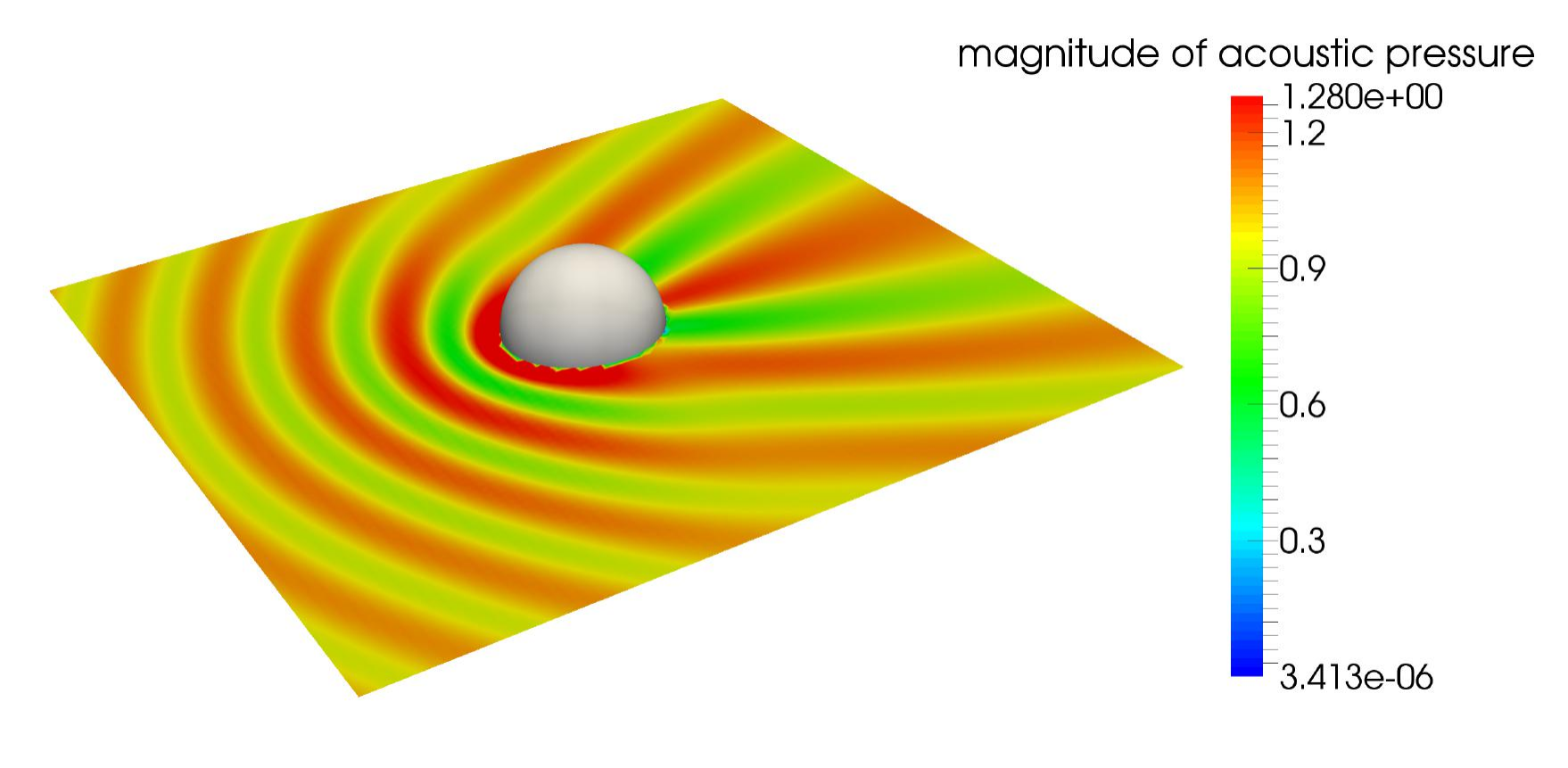}
	\caption{$|p|$ on the $x\mbox{-}y$ plane.}
	\label{fig:hard-sphere-xy-scattering-slice}
	\end{subfigure}
	\caption{Acoustic pressure plots for the coupled sphere scattering problem with acoustically hard surface, $ka=6$. Shell displacements are equal to zero in this case.}
	\label{fig:hard-sphere-results}
\end{figure}

\begin{figure}[h]
	\centering
	\begin{subfigure}[b]{0.45\textwidth}
		\includegraphics[width=\textwidth]{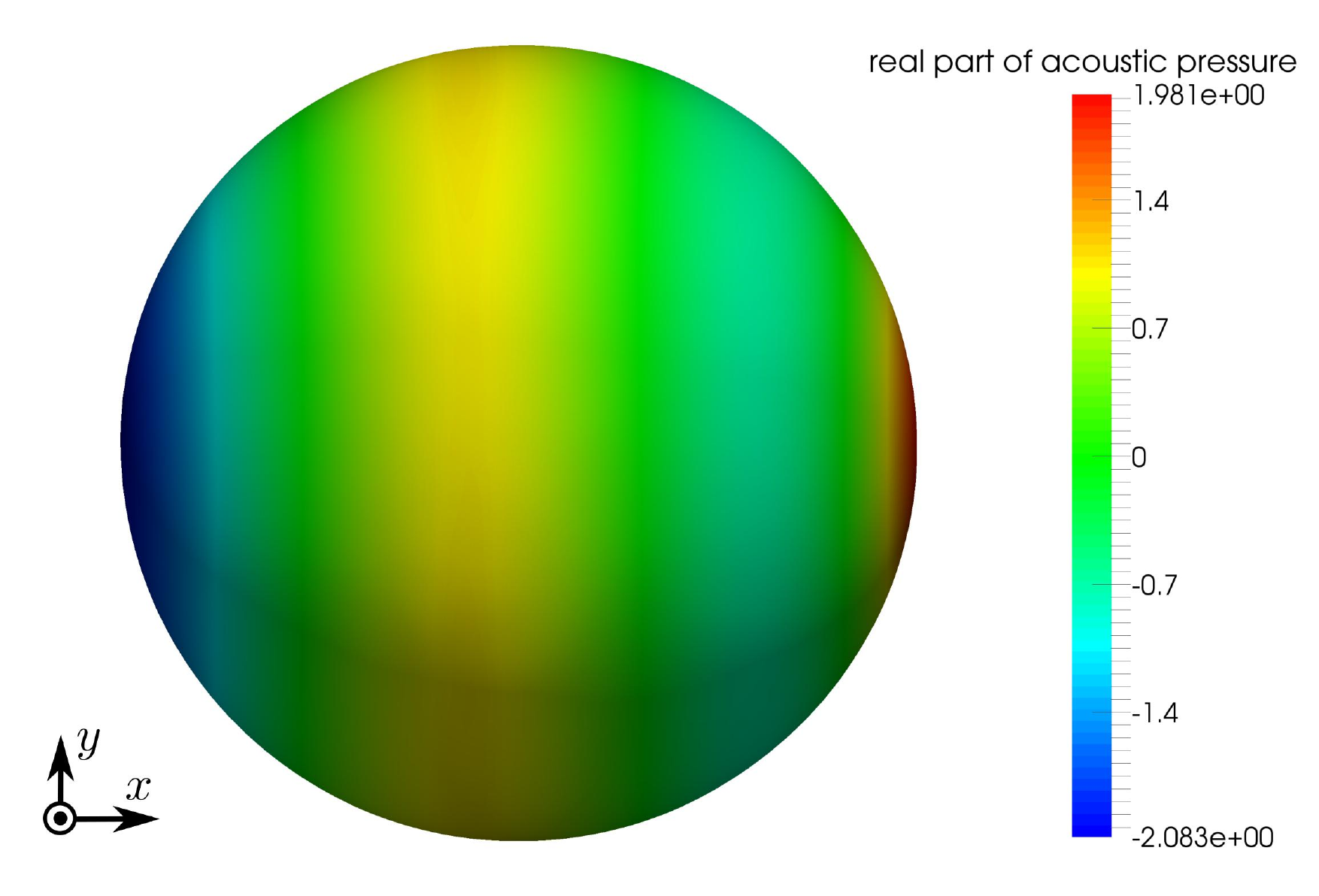}
		\caption{$\textrm{Re}(p)$, $x\mbox{-}y$ plane.}
		\label{fig:elastic-sphere-t0-1-real-acoustic-potential-xy}
	\end{subfigure}
	\begin{subfigure}[b]{0.45\textwidth}
		\includegraphics[width=\textwidth]{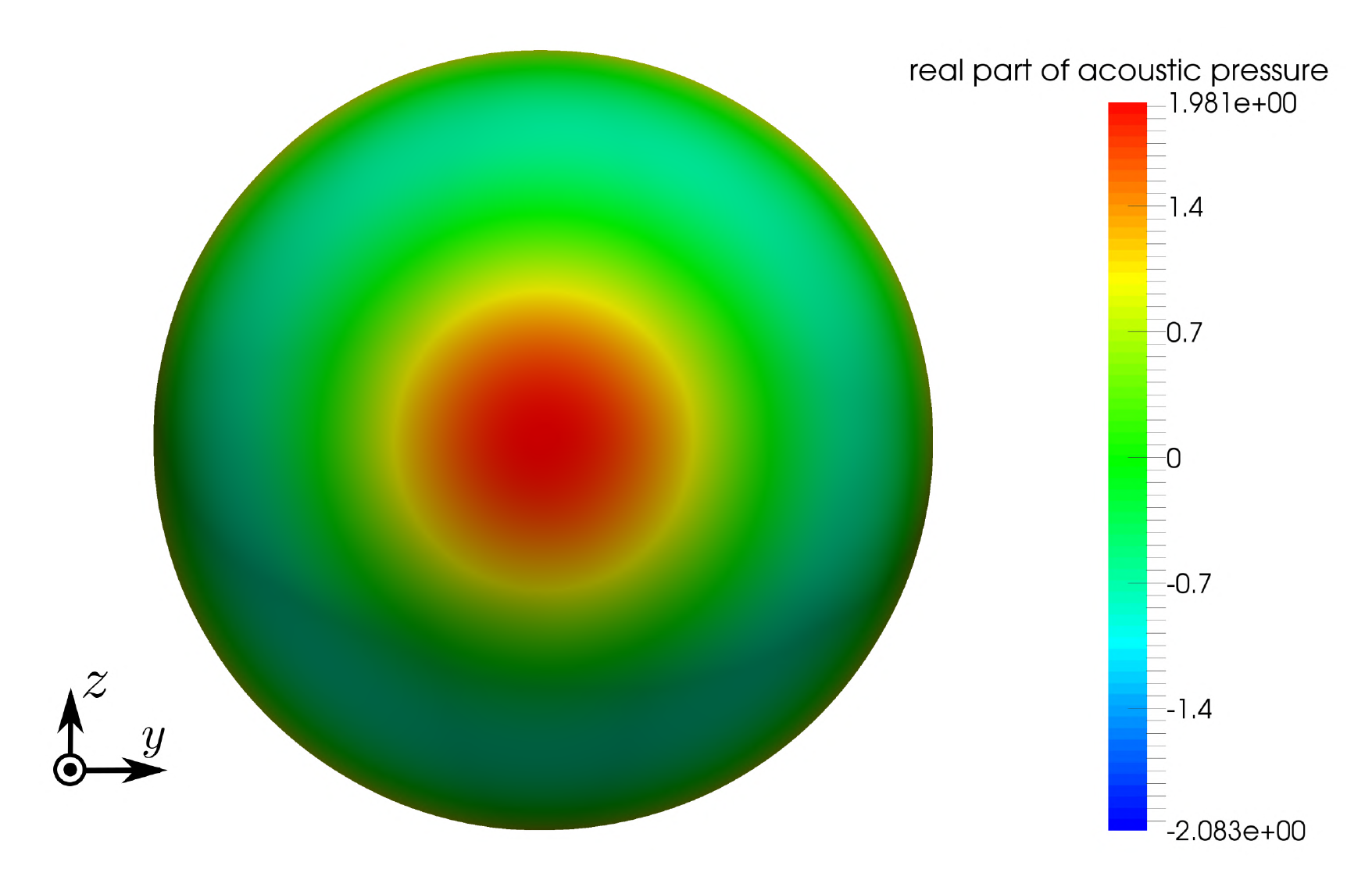}
		\caption{$\textrm{Re}(p)$, $y\mbox{-}z$ plane.}
		\label{fig:elastic-sphere-t0-1-real-acoustic-potential-yz}
	\end{subfigure}

	\vspace{2ex}
	
	\begin{subfigure}[b]{0.45\textwidth}
		\includegraphics[width=\textwidth]{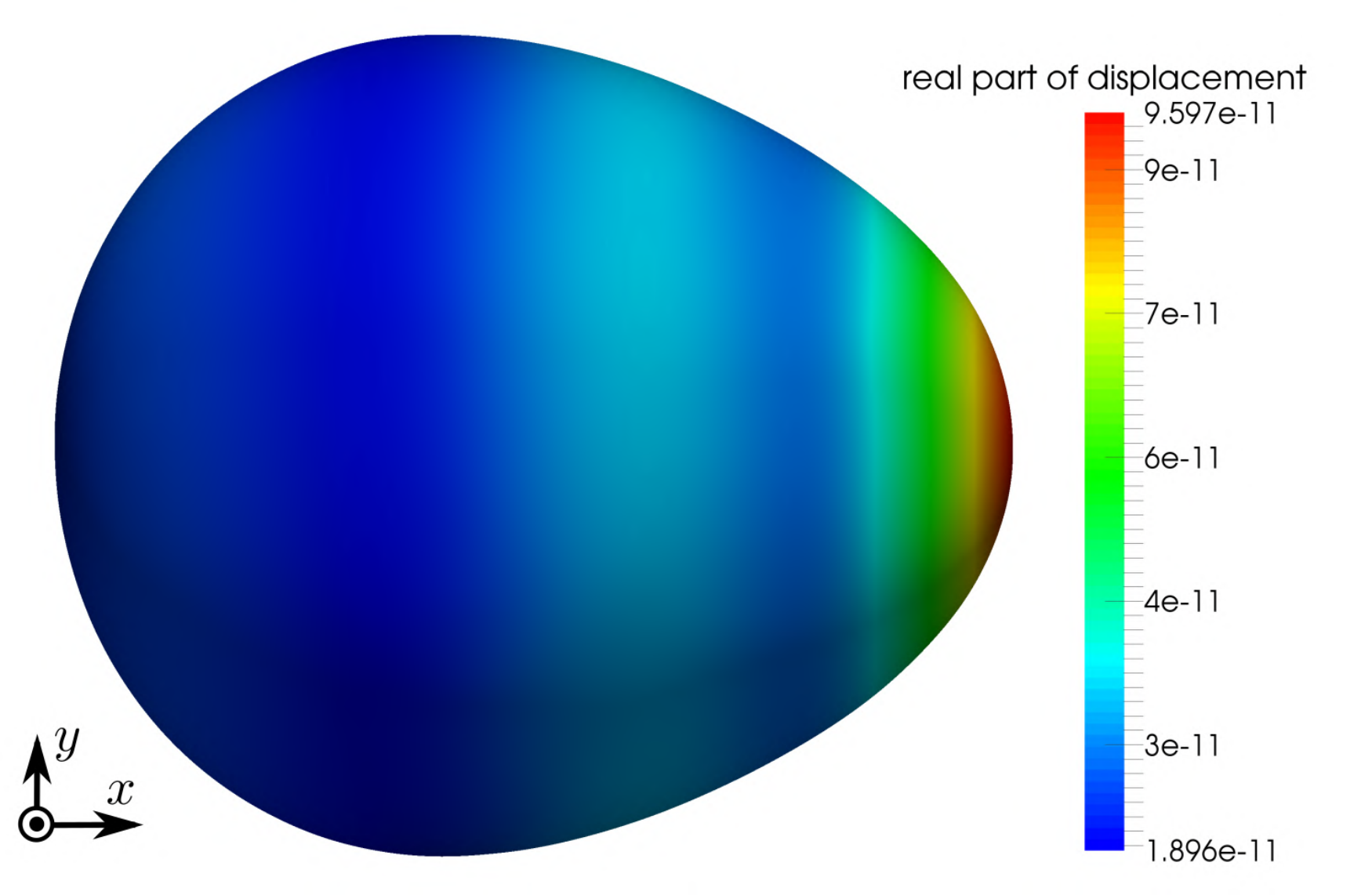}
		\caption{$\textrm{Re}(\mathbf{u})$, $x\mbox{-}y$ plane.}
		\label{fig:elastic-sphere-t0-1-real-displacement-xy}
	\end{subfigure}
	\begin{subfigure}[b]{0.45\textwidth}
		\includegraphics[width=\textwidth]{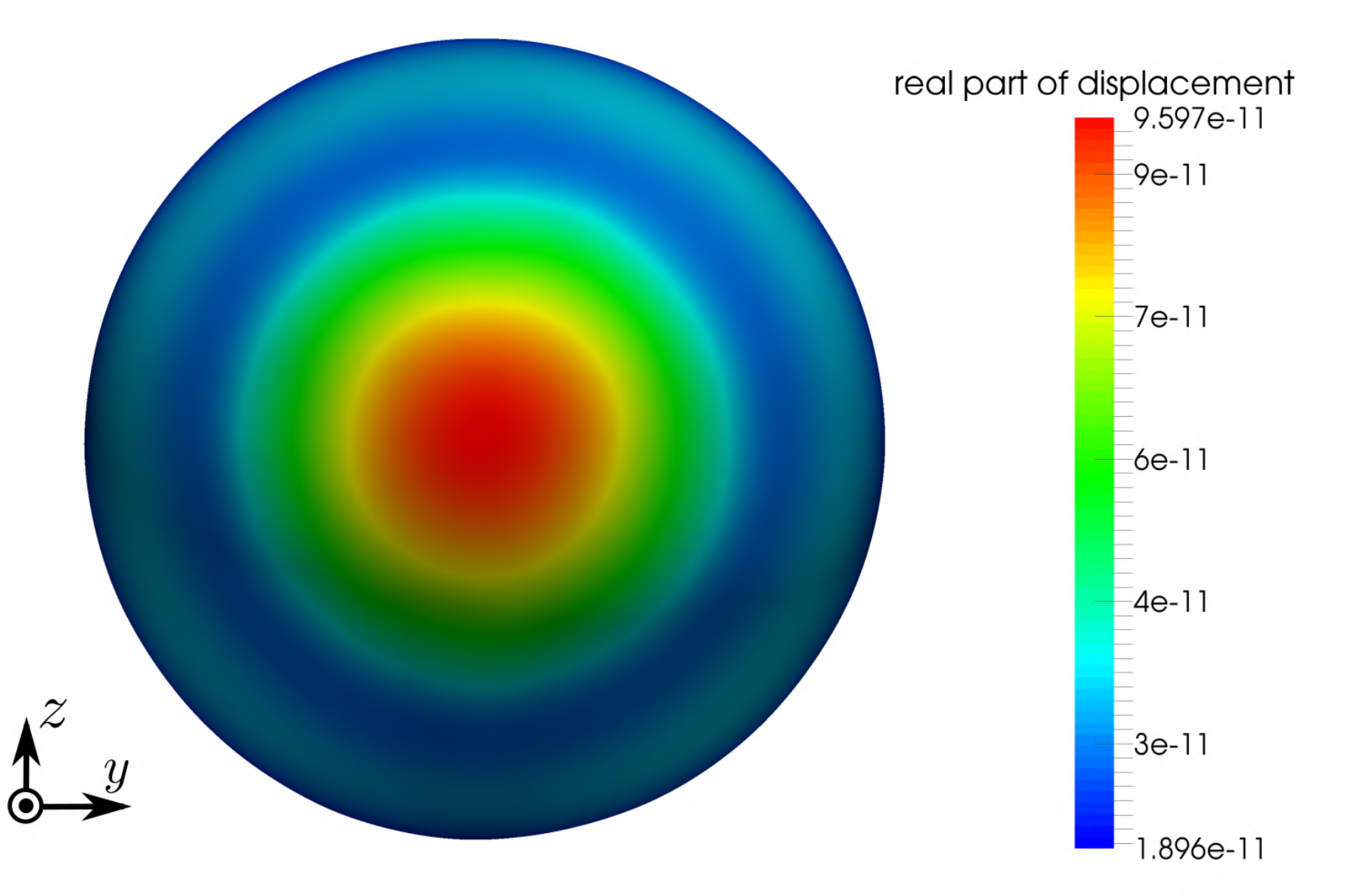}
		\caption{$\textrm{Re}(\mathbf{u})$, $y\mbox{-}z$ plane.}
		\label{fig:elastic-sphere-t0-1-real-displacement-yz}
	\end{subfigure}

	\vspace{2ex}

	\begin{subfigure}[b]{0.8\textwidth}
		\includegraphics[width=\textwidth]{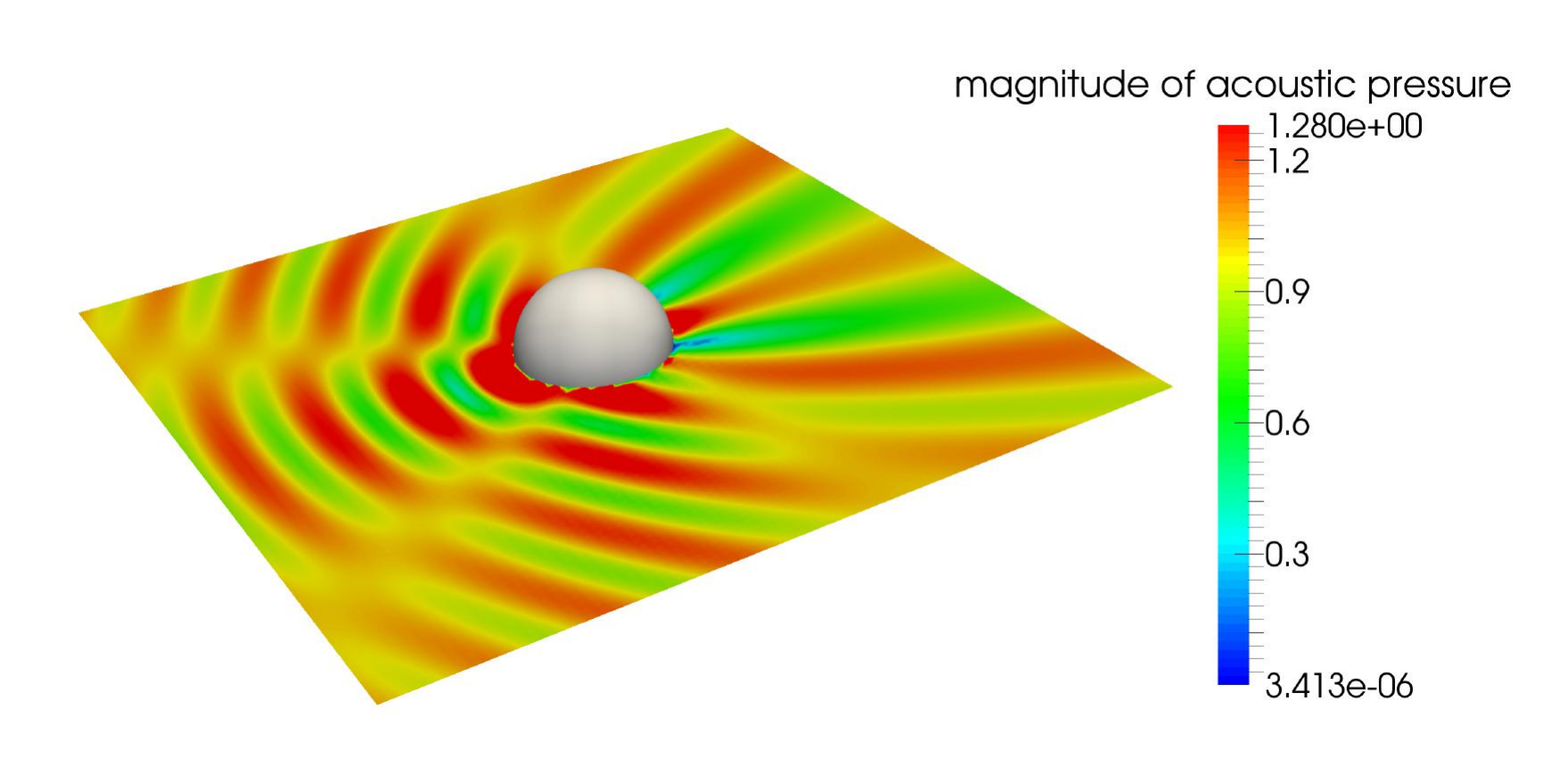}
		\caption{$|p|$ sampled on the $x\mbox{-}y$ plane.}
		\label{fig:elastic-sphere-t0-1-xy-scattering-slice}
	\end{subfigure}

	\caption{Acoustic pressure and displacement plots for the coupled sphere scattering problem with $h=0.1m$, $ka=6$.}
	\label{fig:elastic-sphere-t0-1}
\end{figure}

\begin{figure}[h]
	\centering
	\begin{subfigure}[b]{0.45\textwidth}
		\includegraphics[width=\textwidth]{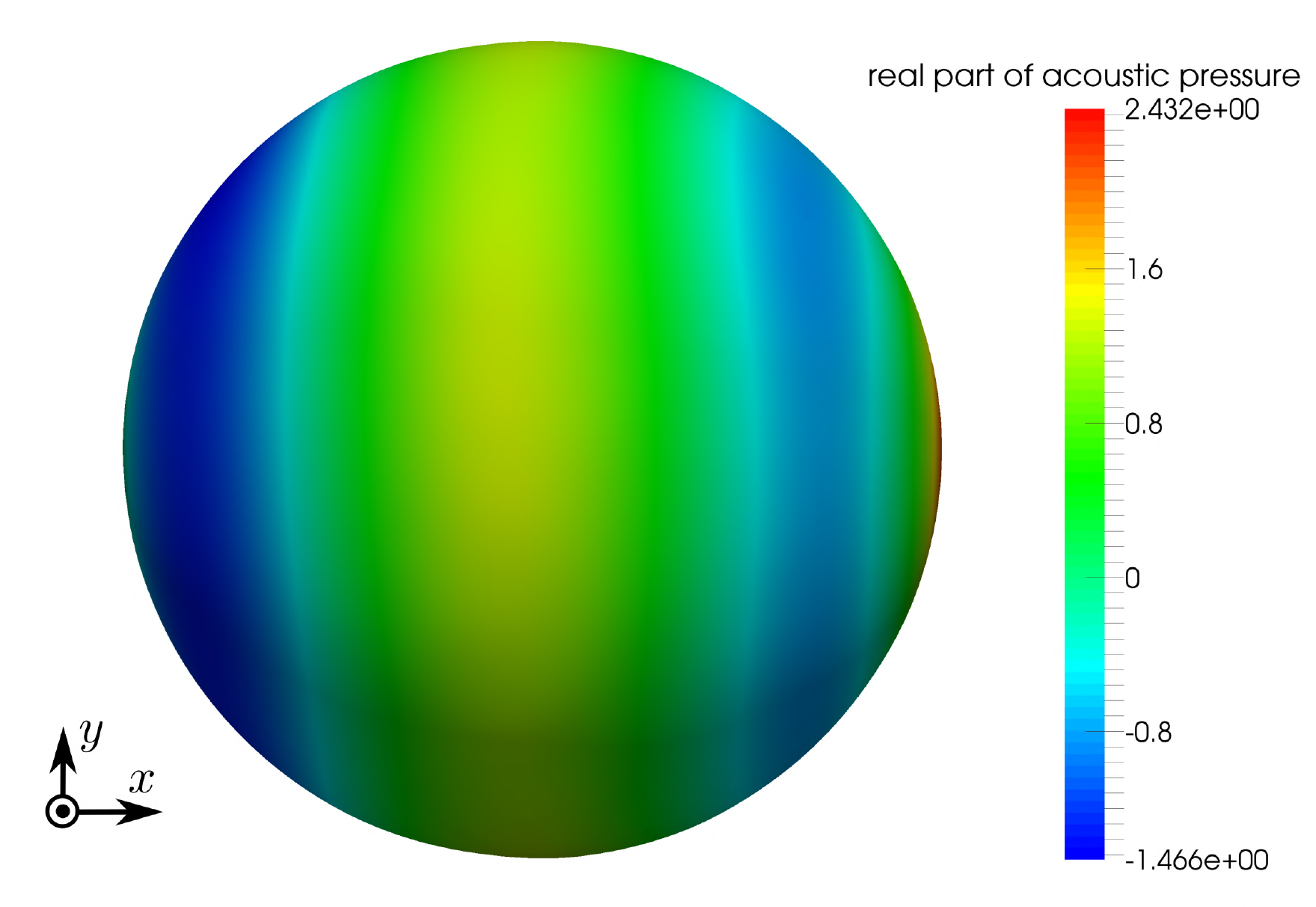}
		\caption{$\textrm{Re}(p)$, $x\mbox{-}y$ plane.}
		\label{fig:elastic-sphere-t0-05-real-acoustic-potential-xy}
	\end{subfigure}
	\begin{subfigure}[b]{0.45\textwidth}
		\includegraphics[width=\textwidth]{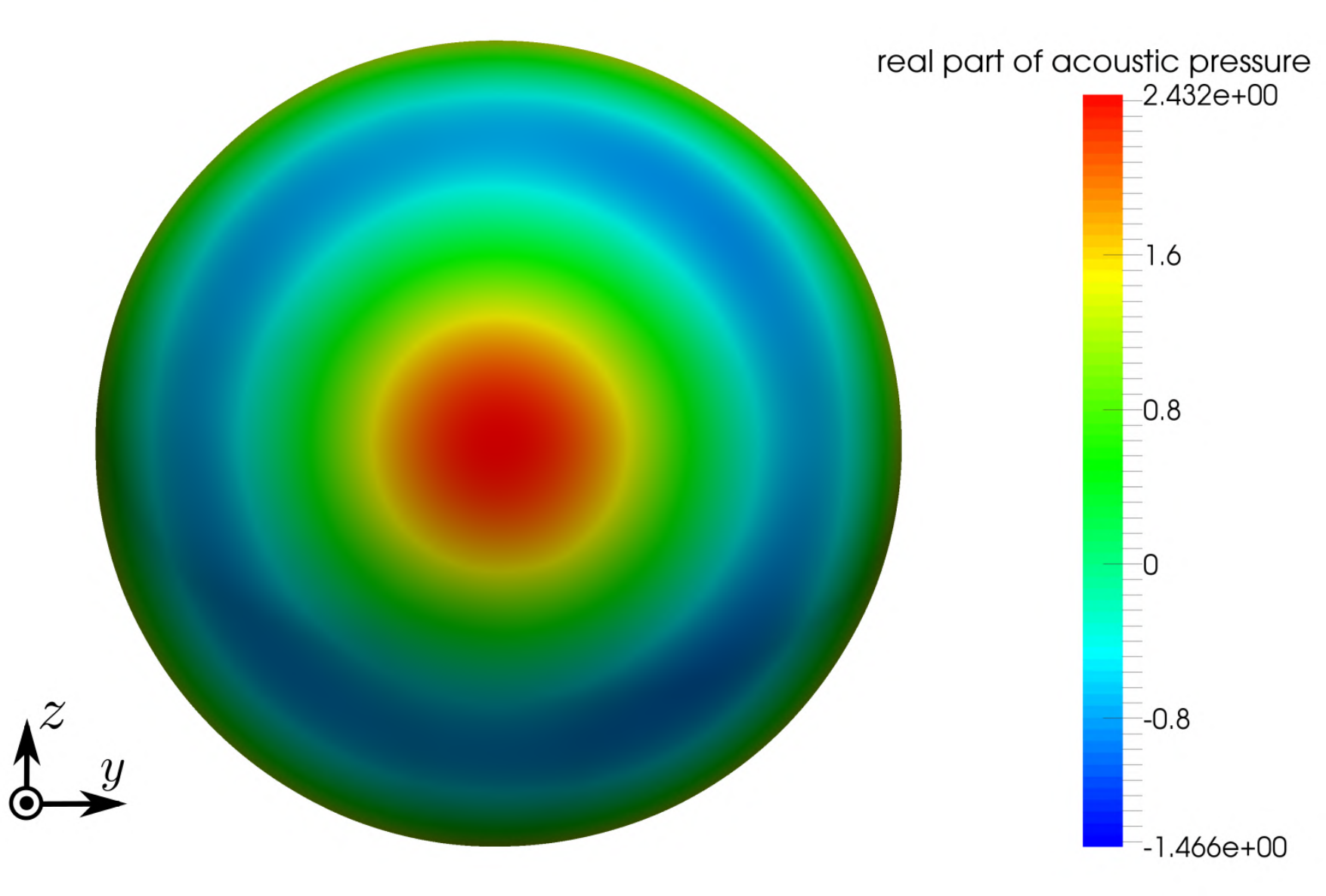}
		\caption{$\textrm{Re}(p)$, $y\mbox{-}z$ plane.}
		\label{fig:elastic-sphere-t0-05-real-acoustic-potential-yz}
	\end{subfigure}
	
	\vspace{2ex}
	
	\begin{subfigure}[b]{0.45\textwidth}
		\includegraphics[width=\textwidth]{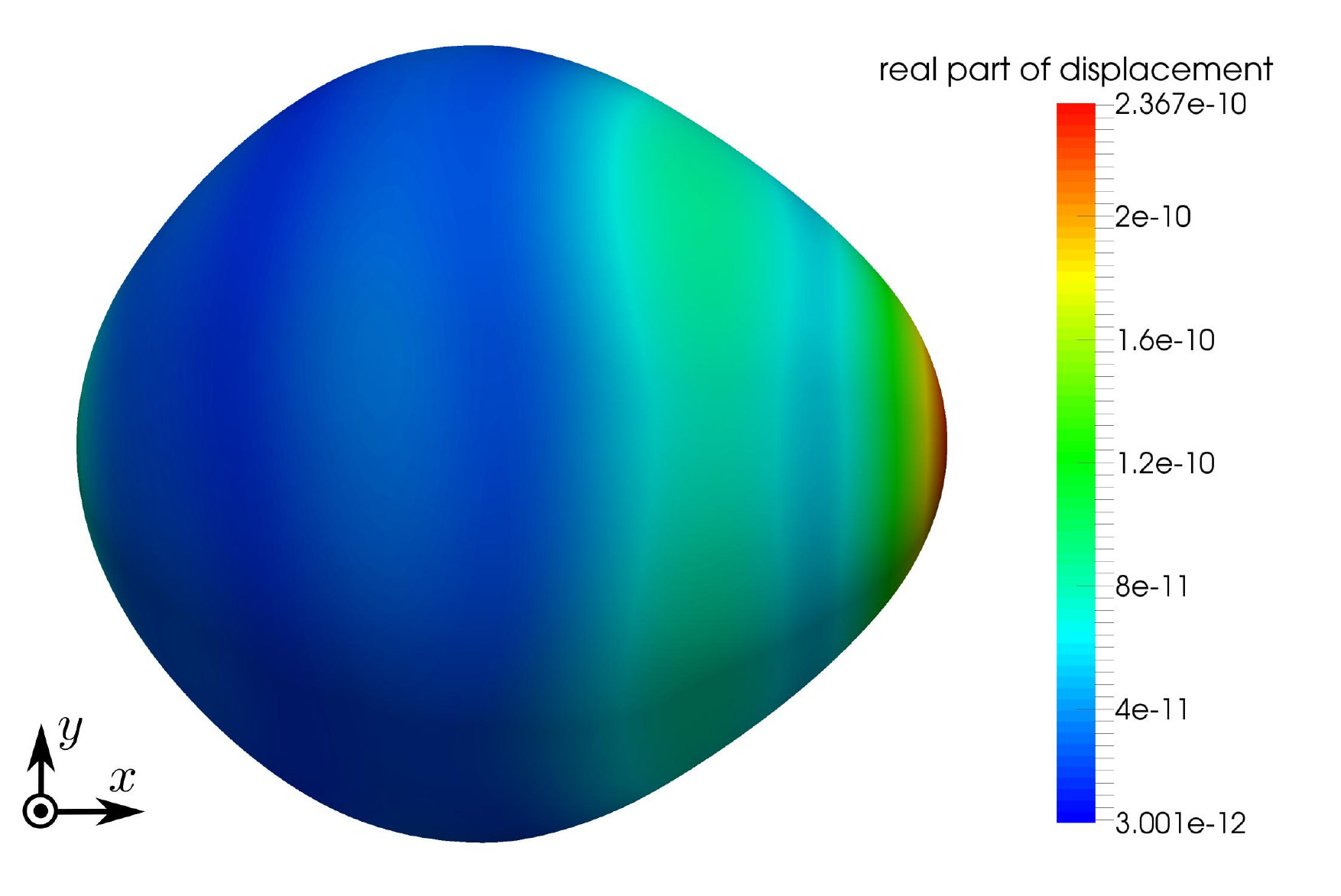}
		\caption{$\textrm{Re}(\mathbf{u})$, $x\mbox{-}y$ plane.}
		\label{fig:elastic-sphere-t0-05-real-displacement-xy}
	\end{subfigure}
	\begin{subfigure}[b]{0.45\textwidth}
		\includegraphics[width=\textwidth]{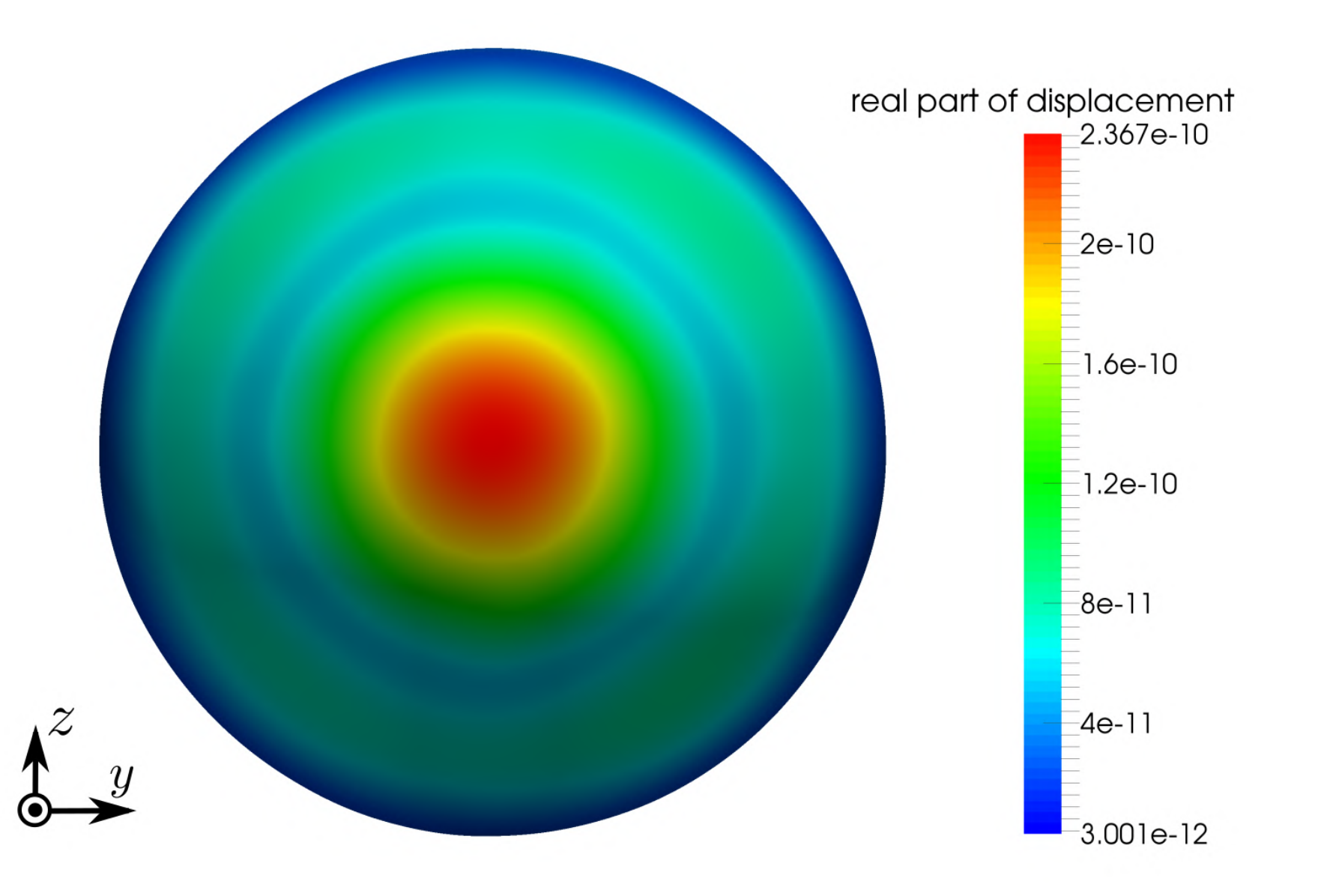}
		\caption{$\textrm{Re}(\mathbf{u})$, $y\mbox{-}z$ plane.}
		\label{fig:elastic-sphere-t0-05-real-displacement-yz}
	\end{subfigure}
	
	\vspace{2ex}
	
	\begin{subfigure}[b]{0.8\textwidth}
		\includegraphics[width=\textwidth]{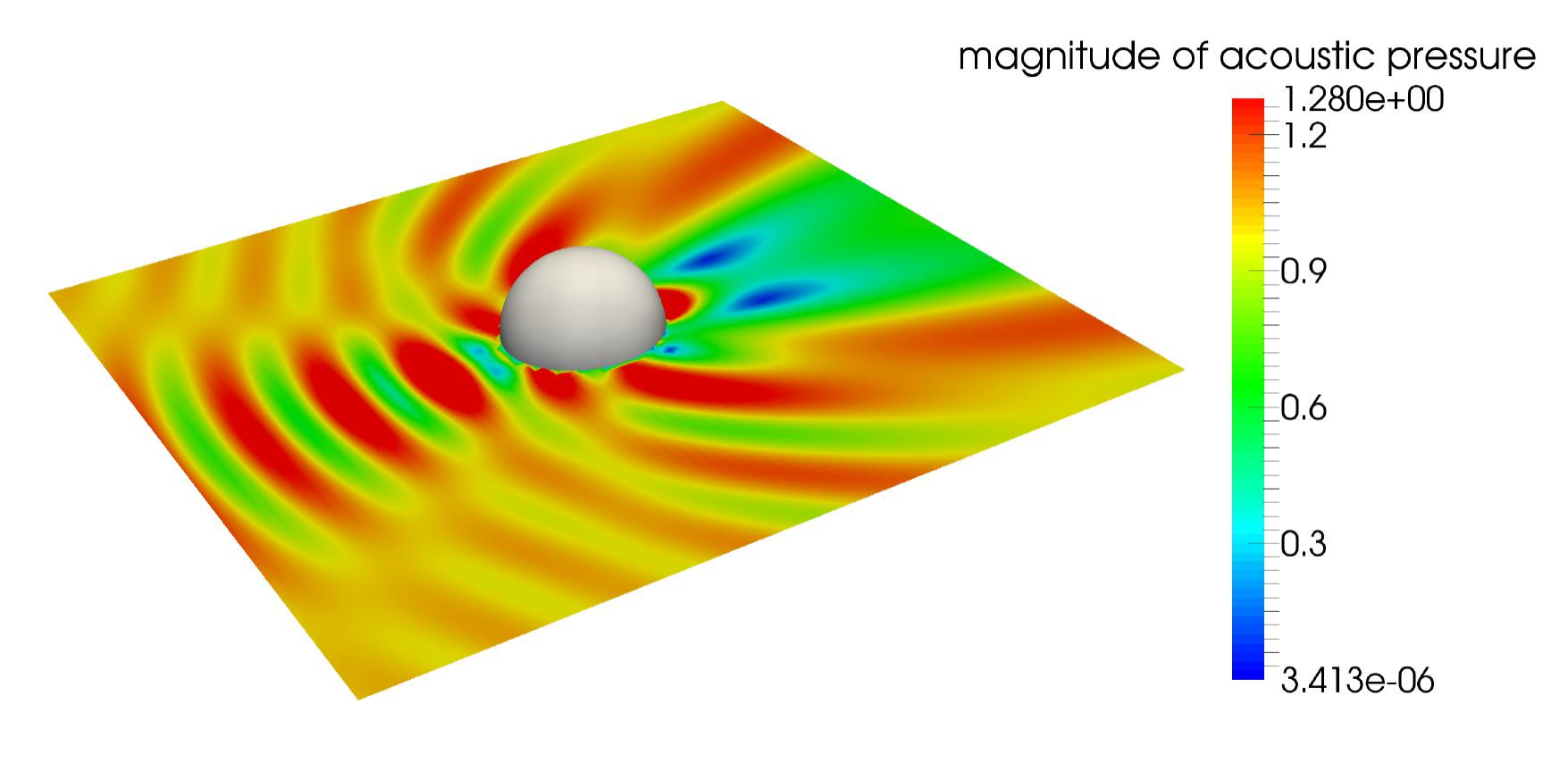}
		\caption{$|p|$ sampled on the $x\mbox{-}y$ plane.}
		\label{fig:elastic-sphere-t0-05-xy-scattering-slice}
	\end{subfigure}
	
	\caption{Acoustic pressure and displacement plots for the coupled sphere scattering problem with $h=0.05m$, $ka=6$.}
	\label{fig:elastic-sphere-t0-05}
\end{figure}

To demonstrate the coupling effect created by elastic displacements of the shell structure we conduct three studies using the sphere geometry with our Loop subdivision discretisation approach:
\begin{enumerate}
	\item\label{item:coupling-effect-hard} \textbf{Acoustically hard surface}: we specify the surface as acoustically hard which eliminates the coupling effect. The normal derivative of the acoustic potential is set to zero on the sphere surface (i.e. $\partial p/ \partial n = 0$ for all $\mathbf{x} \in \Gamma$) and no displacement of the shell occurs.
	\item\label{item:coupling-effect-t-0-1} \textbf{Elastic shell, $\boldsymbol{h}$=$\boldsymbol{0.1m}$}: a fully coupled system is formed where shell vibrations result in a radiated acoustic pressure.
	\item\label{item:coupling-effect-t-0-05} \textbf{Elastic shell, $\boldsymbol{h}$=$\boldsymbol{0.05m}$}: the same study as for Item \ref{item:coupling-effect-t-0-1} but with a reduced shell thickness. 
\end{enumerate}
We specify a low wavenumber of $ka=6$ in keeping with the acoustic  scattering study of~\cite{simpsonscott2014} and use the initial control mesh to generate a Loop subdivision discretisation with approximately thirteen elements per wavelength.  Results for scenarios \ref{item:coupling-effect-hard}, \ref{item:coupling-effect-t-0-1} and \ref{item:coupling-effect-t-0-05} are illustrated in Figures~\ref{fig:hard-sphere-results}, \ref{fig:elastic-sphere-t0-1} and \ref{fig:elastic-sphere-t0-05} respectively.  A comparison of the far-field acoustic potential magnitude for the hard sphere case against the elastic cases (Figures~\ref{fig:hard-sphere-xy-scattering-slice}, \ref{fig:elastic-sphere-t0-1-xy-scattering-slice} and \ref{fig:elastic-sphere-t0-05-xy-scattering-slice}) reveals the noticeable effect of the radiated acoustic pressure caused by shell vibrations.  It is also clear that a decrease in shell thickness has a substantial effect on the scattered profile where a region of low pressure is created in the shadow region (see Figure~\ref{fig:elastic-sphere-t0-05-xy-scattering-slice}) which is not observed in the other cases.  In addition, an inspection of Figures~\ref{fig:elastic-sphere-t0-1-real-acoustic-potential-xy}, \ref{fig:elastic-sphere-t0-1-real-acoustic-potential-yz}, \ref{fig:elastic-sphere-t0-05-real-acoustic-potential-xy} and \ref{fig:elastic-sphere-t0-05-real-acoustic-potential-yz} reveals that the maximum value of $Re(p)$ shifts from the illuminating region to the shadow region when the shell thickness is reduced.  We finally remark that as expected, shell displacements increase when shell thickness is reduced (Figures~\ref{fig:elastic-sphere-t0-1-real-displacement-xy},\ref{fig:elastic-sphere-t0-1-real-displacement-yz}, \ref{fig:elastic-sphere-t0-05-real-displacement-xy}  and \ref{fig:elastic-sphere-t0-05-real-displacement-yz}) with maximum values consistently obtained at a position opposite to the incident point (first point of contact on the surface by the incident wave).

% hard sphere results

% t=0.1 sphere plots

% t=0.05 sphere plots

\clearpage
% *****************************************************
% ********** Submarine model ***********
% *****************************************************
\subsection{Submarine model}
\begin{figure}[h]
\centering
\includegraphics[width=0.7\textwidth]
{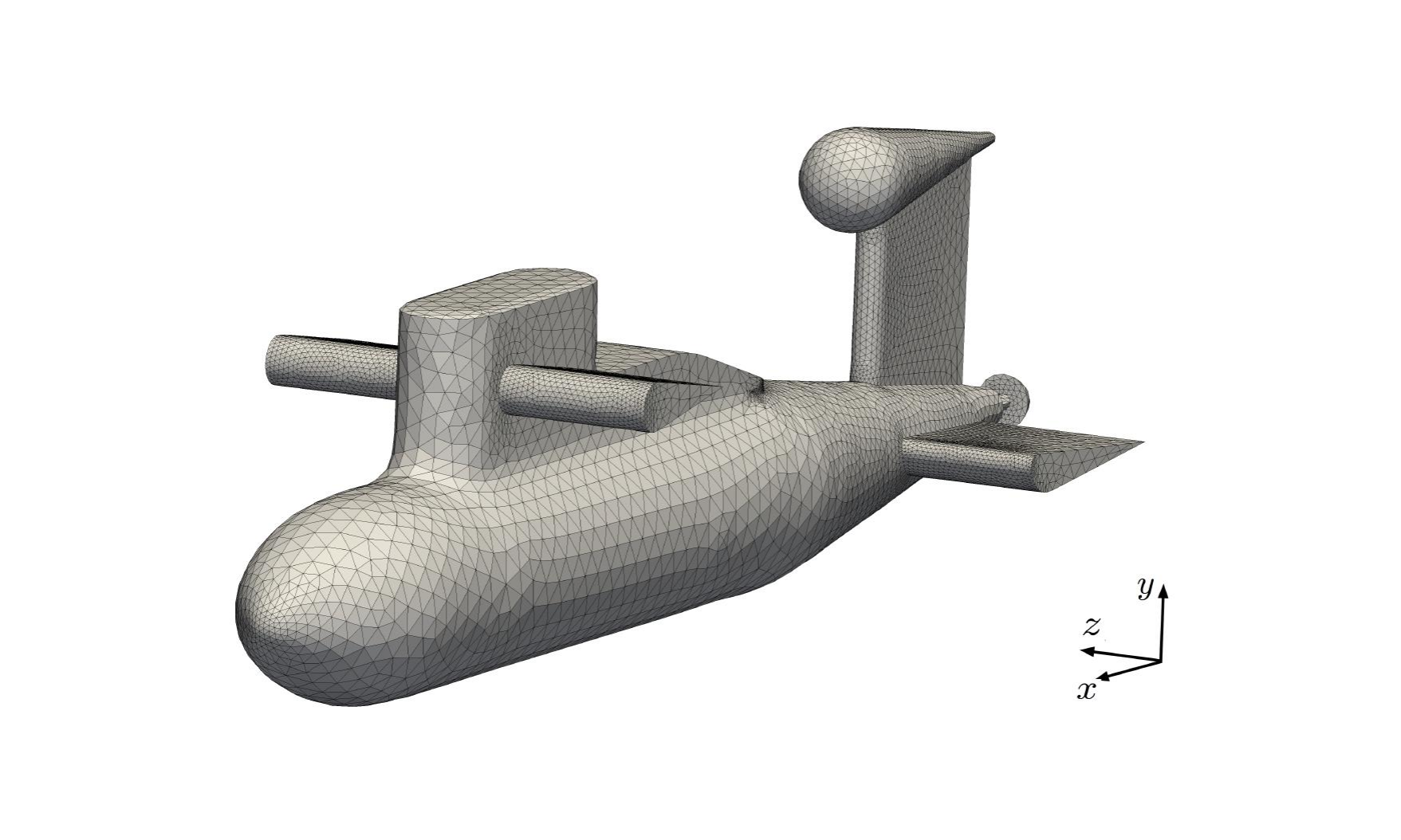}
\caption{Submarine problem: control mesh with $9,510$ vertices.}
\label{fig:submarine-mesh}
\end{figure}

\begin{figure}[h]
\centering
\includegraphics[width=0.7\textwidth]
{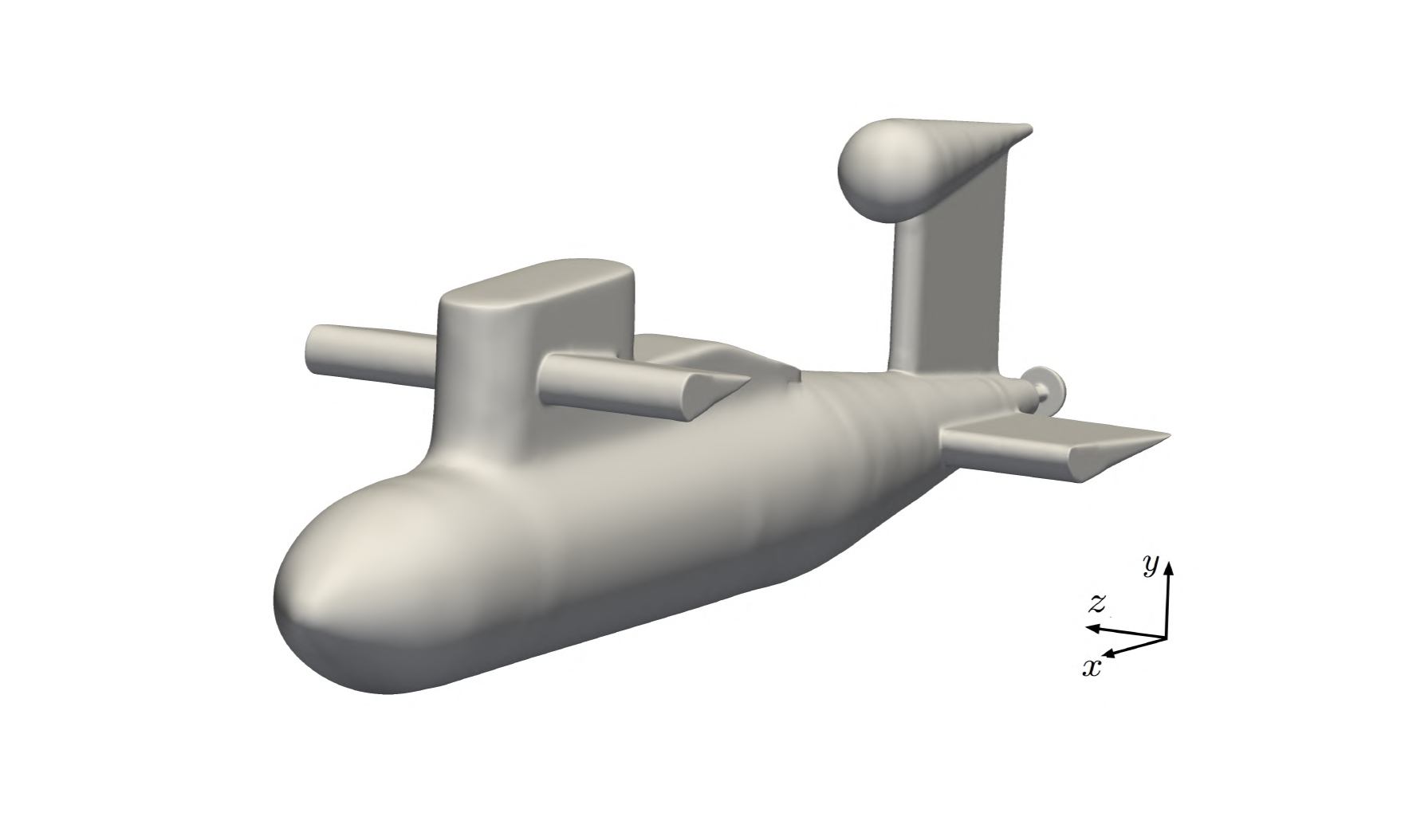}
\caption{The smooth limit surface of the submarine model generated through Loop subdivision.}
\label{fig:submarine-surface}
\end{figure}

We now consider the ability of the present method to provide high-order discretisations of arbitrarily complex geometries. The first model we consider is a submarine model with a control mesh illustrated in Figure~\ref{fig:submarine-mesh} and limit surface as shown in Figure~\ref{fig:submarine-surface}. The minimum bounding box for this model is defined by $[x_i^{min},x_i^{max}]^3 =  [-51.3,41.0]\narrowtimes[-58.4,17.8]\narrowtimes[-11.8,11.8]$. The model is symmetric about the $x\mbox{-}y$ plane. The material properties as specified in Table~\ref{tab:spherical-scattering-properties} are applied with a shell thickness of $h=0.5m$. We construct a forcing function through an incident plane wave with unit magnitude travelling in the negative $x$-direction and a normalised wavenumber of $ka = 46.15$. 

\begin{figure}[h]
	\centering
	\begin{subfigure}[b]{0.49\textwidth}
		\includegraphics[width=\textwidth]{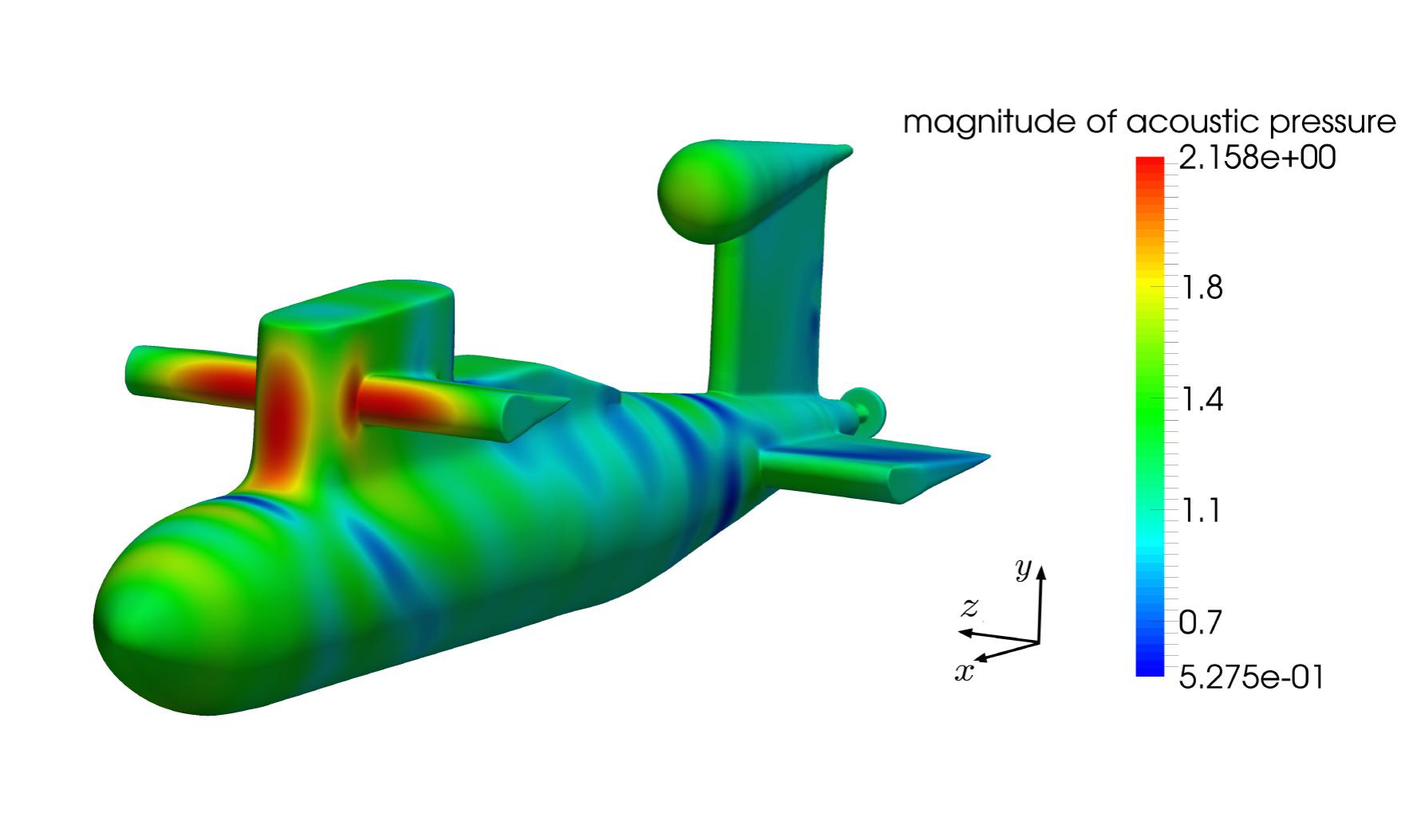}
		\caption{Acoustically hard surface, perspective view.}
		\label{fig:scattering_3D_submarine}
	\end{subfigure}
	\begin{subfigure}[b]{0.49\textwidth}
		\includegraphics[width=\textwidth]{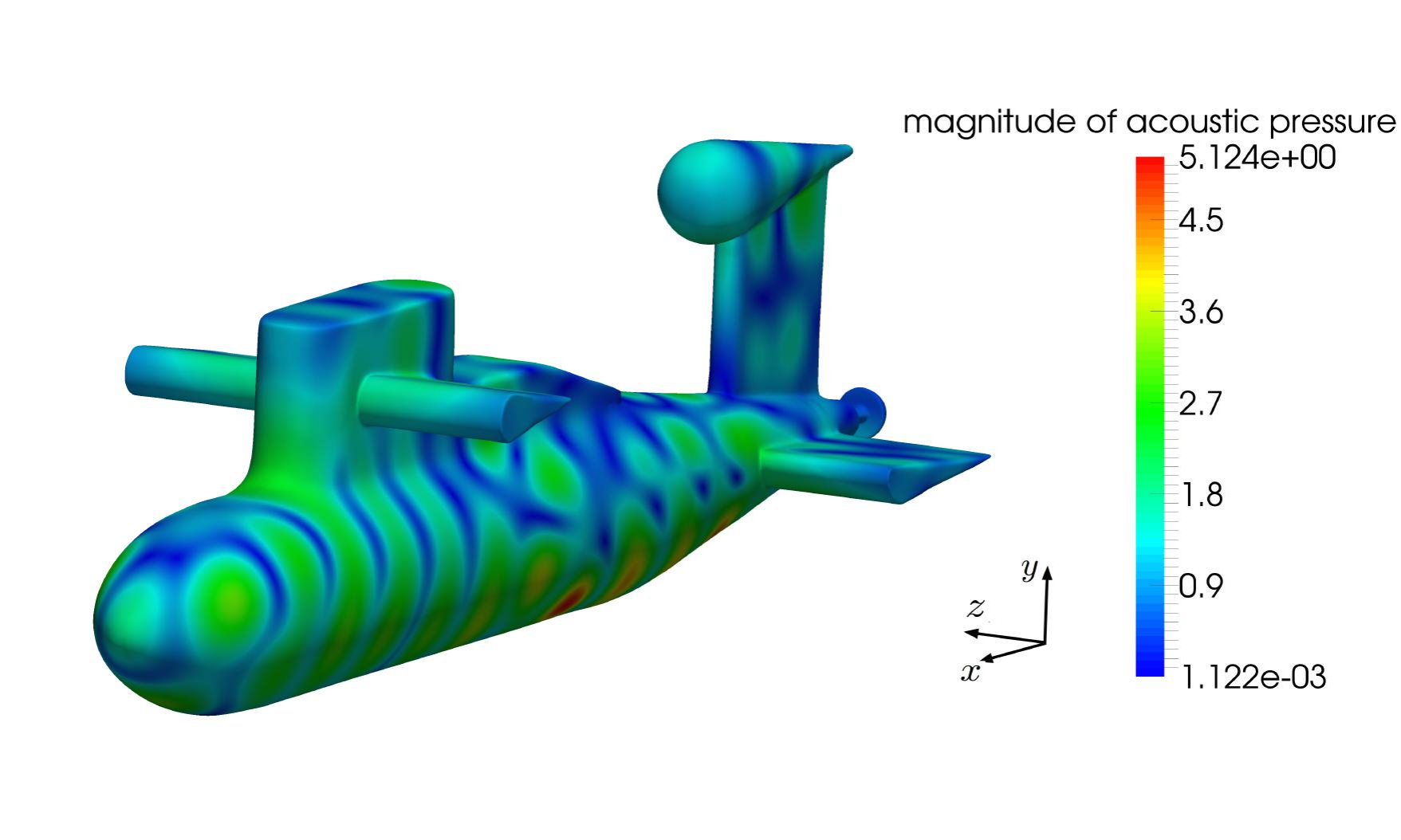}
		\caption{Coupled elastic shell, perspective view.}
		\label{fig:coupled_3D_submarine}
	\end{subfigure}
	
	\vspace{2ex}
	
	\begin{subfigure}[b]{0.49\textwidth}
		\includegraphics[width=\textwidth]{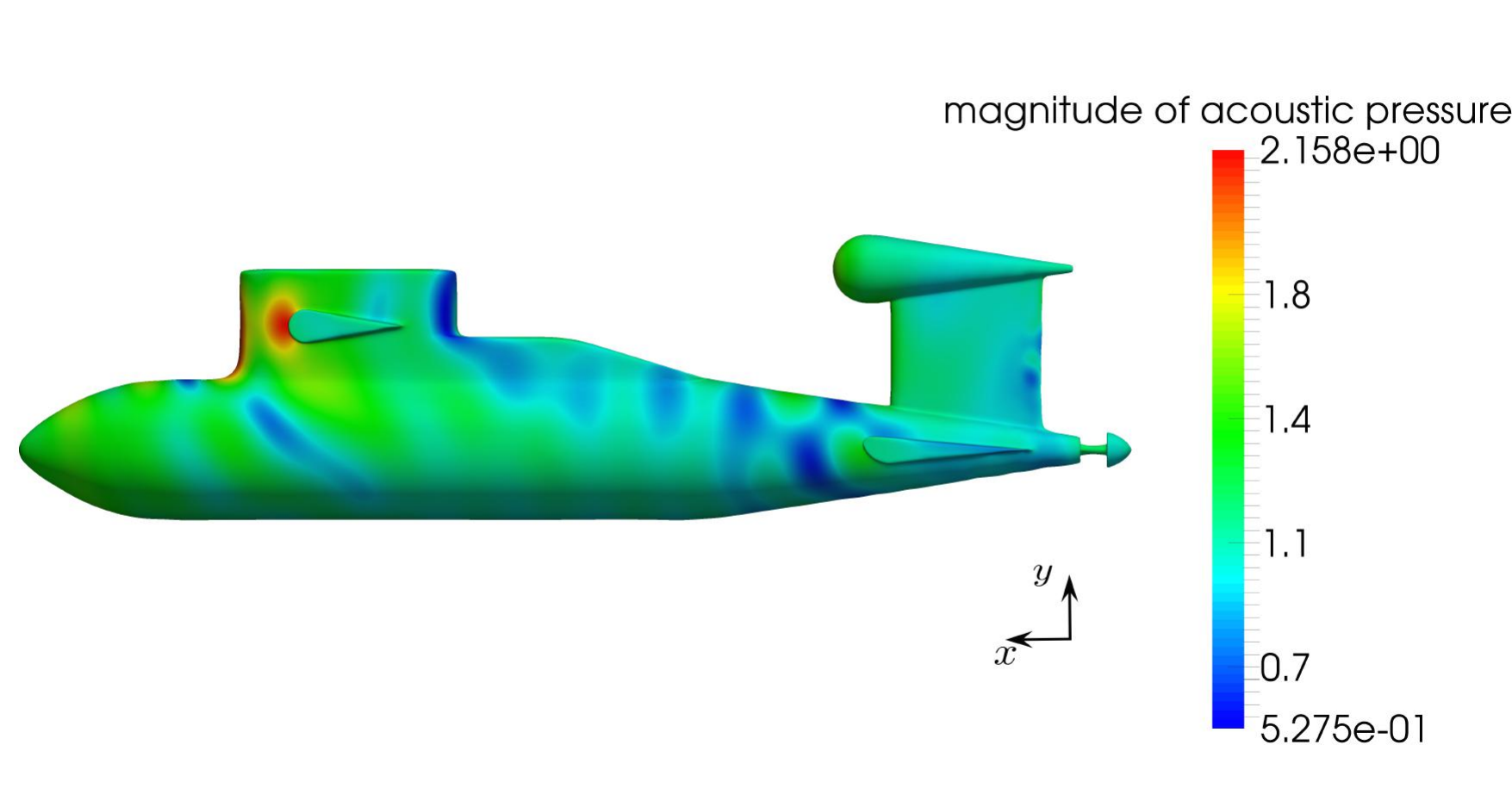}
		\caption{Acoustically hard surface, $x\mbox{-}y$ plane.}
		\label{fig:scattering_xy_submarine}
	\end{subfigure}
	\begin{subfigure}[b]{0.49\textwidth}
		\includegraphics[width=\textwidth]{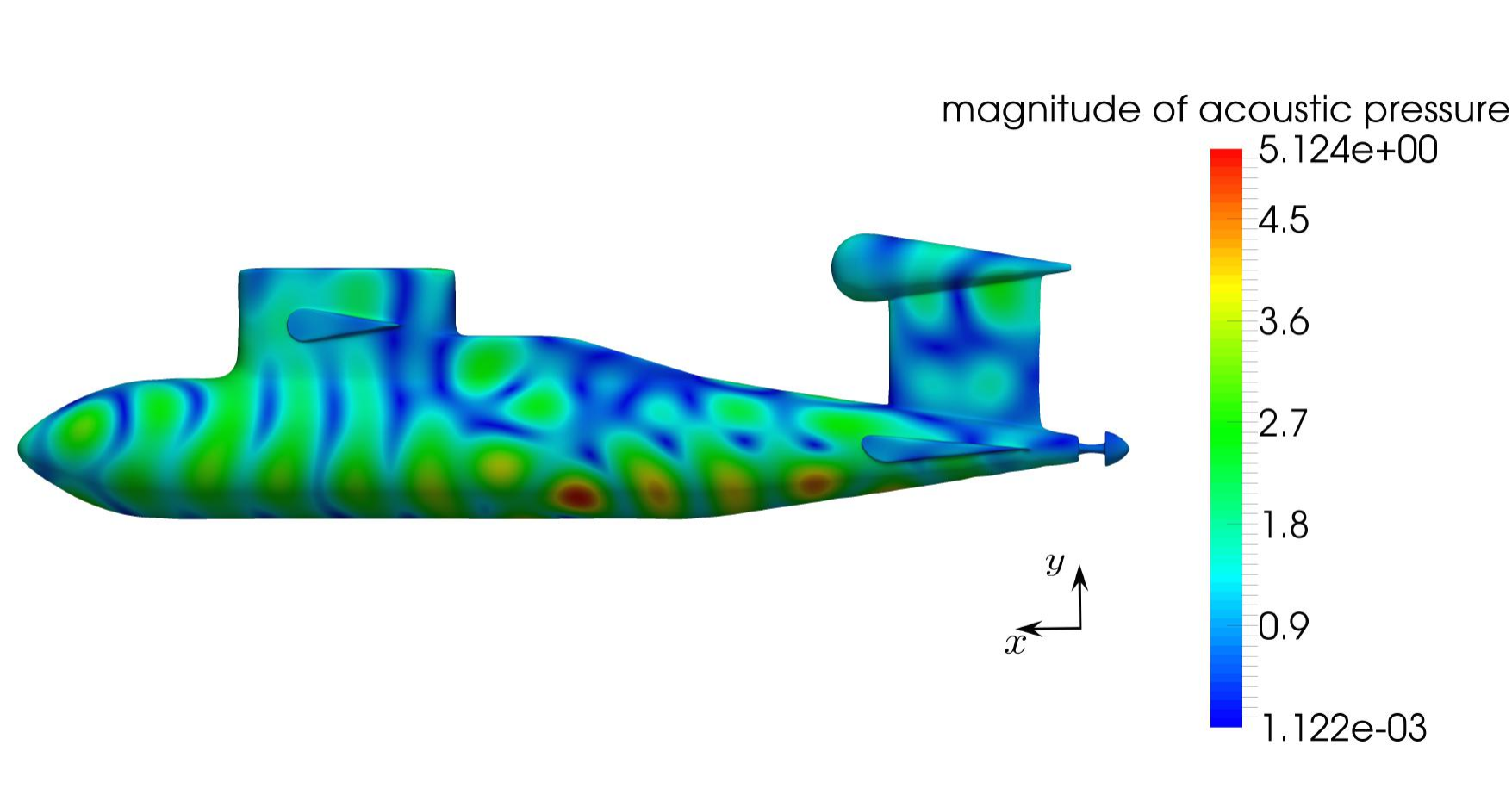}
		\caption{Coupled elastic shell, $x\mbox{-}y$ plane.}
		\label{fig:coupled_xy_submarine}
	\end{subfigure}

	\vspace{2ex}
		
	\begin{subfigure}[b]{0.49\textwidth}
		\includegraphics[width=\textwidth]{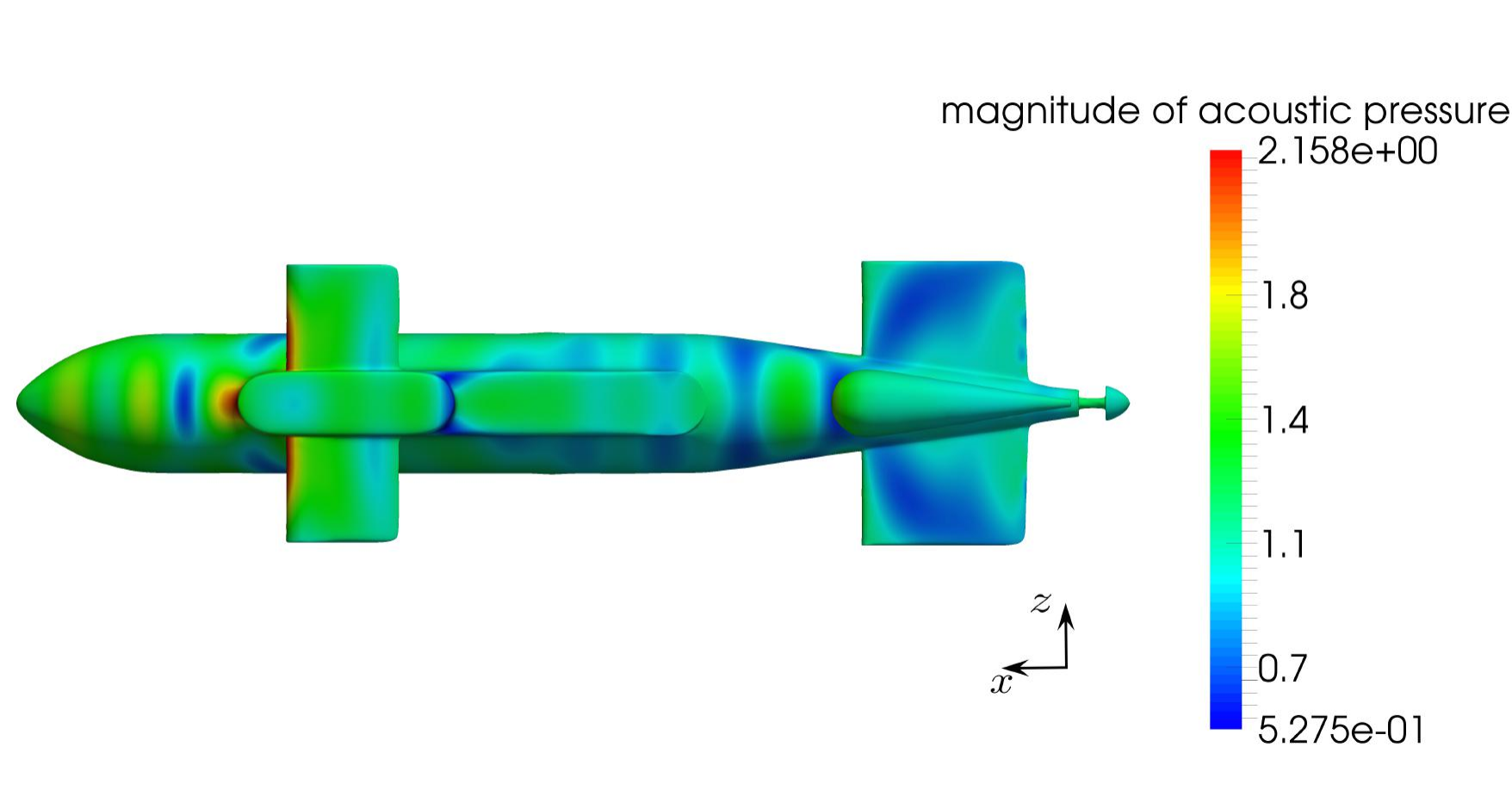}
		\caption{Acoustically hard surface, $x\mbox{-}z$ plane.}
		\label{fig:scattering_xz_submarine}
	\end{subfigure}
	\begin{subfigure}[b]{0.49\textwidth}
		\includegraphics[width=\textwidth]{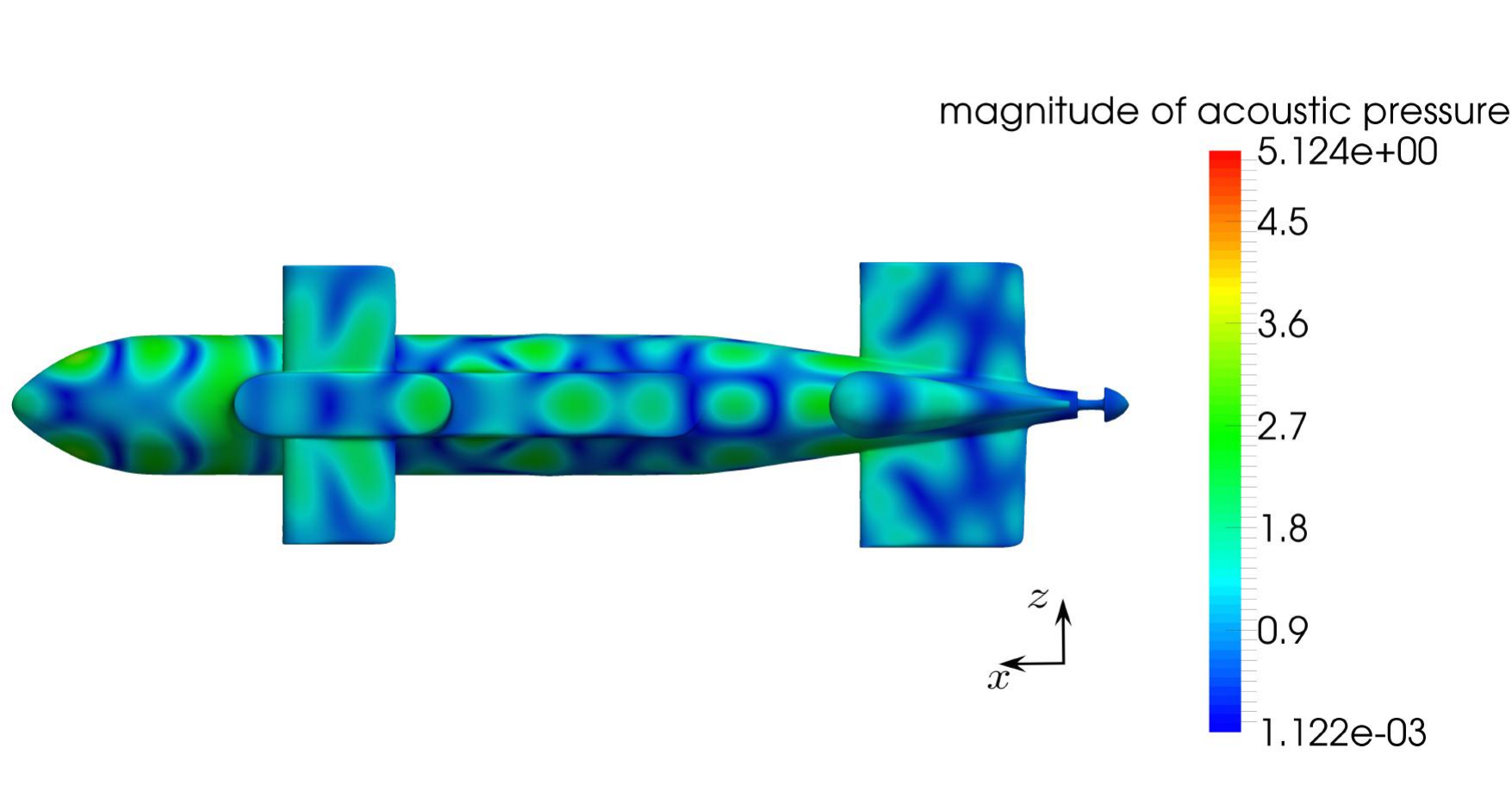}
		\caption{Coupled elastic shell, $x\mbox{-} z$ plane.}
		\label{fig:coupled_xz_submarine}	
	\end{subfigure}
	
	\caption{ Submarine model: acoustic pressure magnitude profiles for both an acoustically hard surface (i.e. $\partial p/ \partial n = 0$) and coupled shell formulation with a normalised wavenumber of $ka = 46.15$.}
	\label{fig:submarine-comparision}
\end{figure}

A comparison between profiles of the acoustic pressure magnitude are illustrated in Figure~\ref{fig:submarine-comparision} for both an acoustically hard surface and elastic shell formulation with the effect on the radiated acoustic pressure caused by shell vibrations apparent.  We primarily use this example to illustrate the ability of our approach to handle large models of industrial relevance making  use of smooth geometry representations that are generated by CAD software.

\clearpage

% *****************************************************
% ********** Complex topology example ***********
% *****************************************************

\subsection{Complex topology example}
%\subsubsection{Complex topology example}

\begin{figure}[h]
    \centering
    \begin{subfigure}{0.5\textwidth}
        \includegraphics[width=\textwidth]{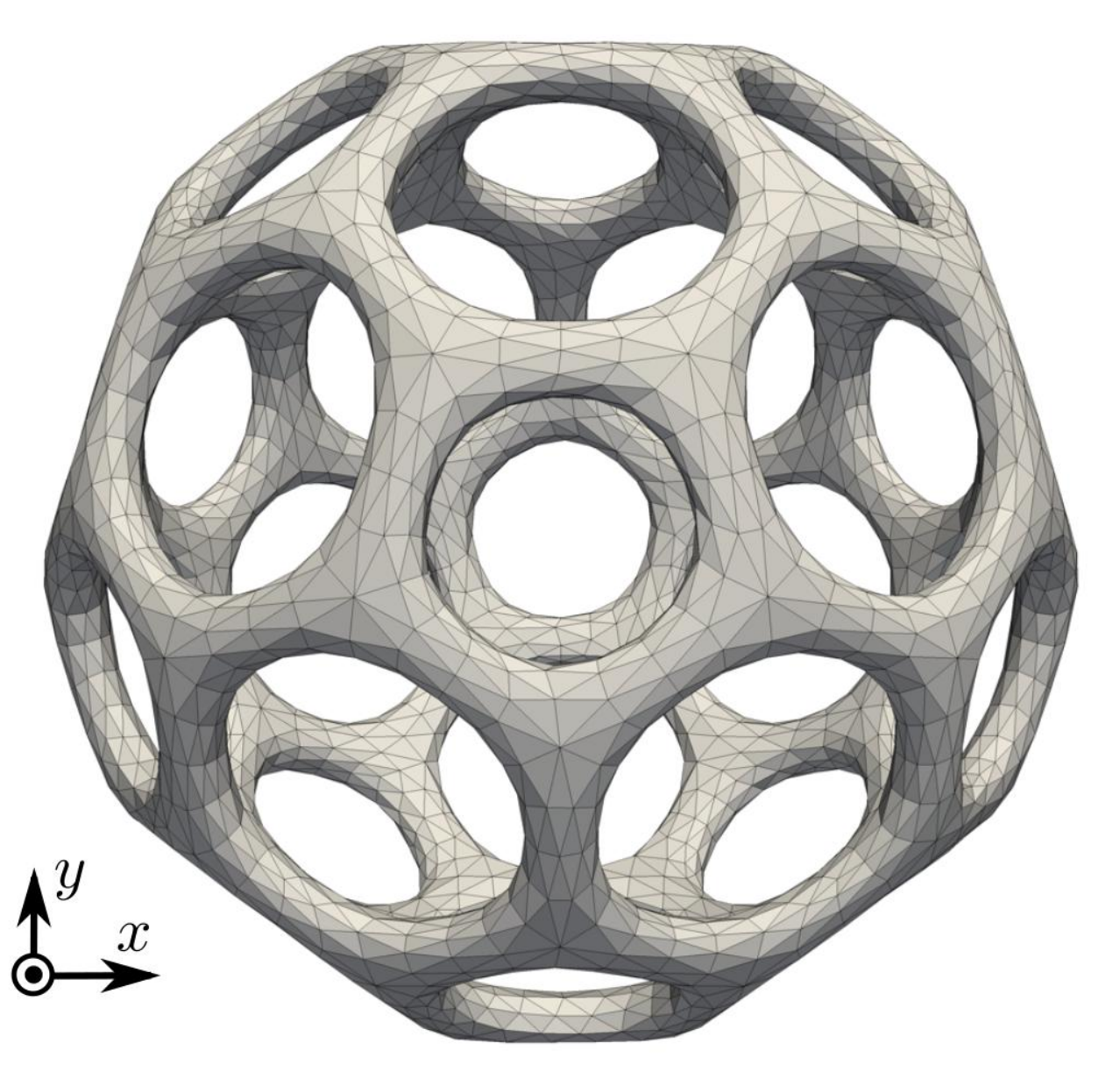}
        \caption{Control mesh.}
        \label{fig:molecule-control-grid}
    \end{subfigure}
    
    \begin{subfigure}{0.5\textwidth}
        \includegraphics[width=\textwidth]{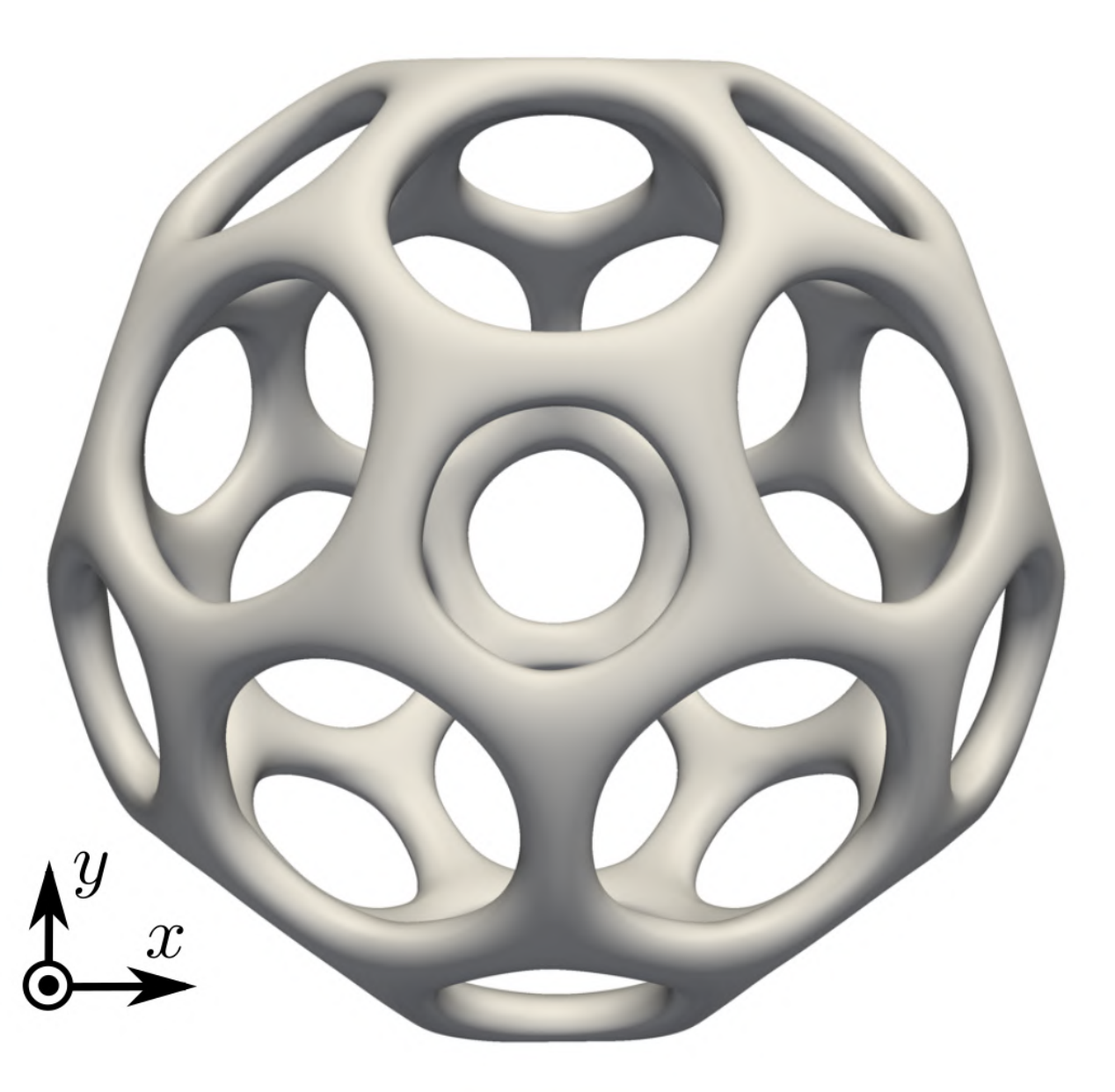}
        \caption{Limit surface.}
        \label{fig:molecule-limit-surface}
    \end{subfigure}
\caption{Complex topology example: Loop subdivision discretisation with 3,940 vertices. }
    \label{fig:molecule-discretisation}
\end{figure}

\begin{figure}[h]
    \centering
    \begin{subfigure}{0.7\textwidth}
        \includegraphics[width=\textwidth]{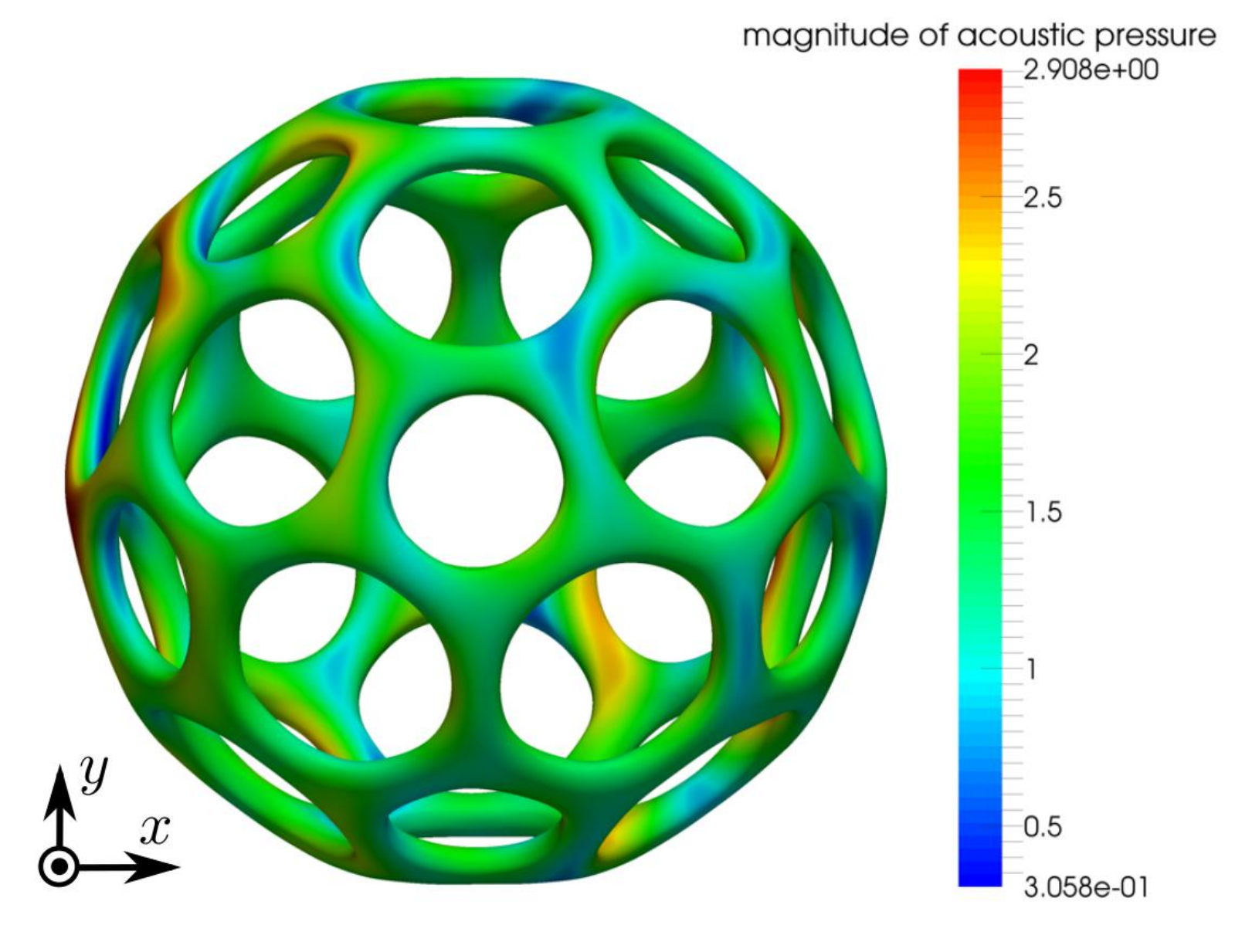}
        \caption{$|p|$}
        \label{fig:molecule-surface-results}
    \end{subfigure}
    
    \begin{subfigure}{0.7\textwidth}
        \includegraphics[width=\textwidth]{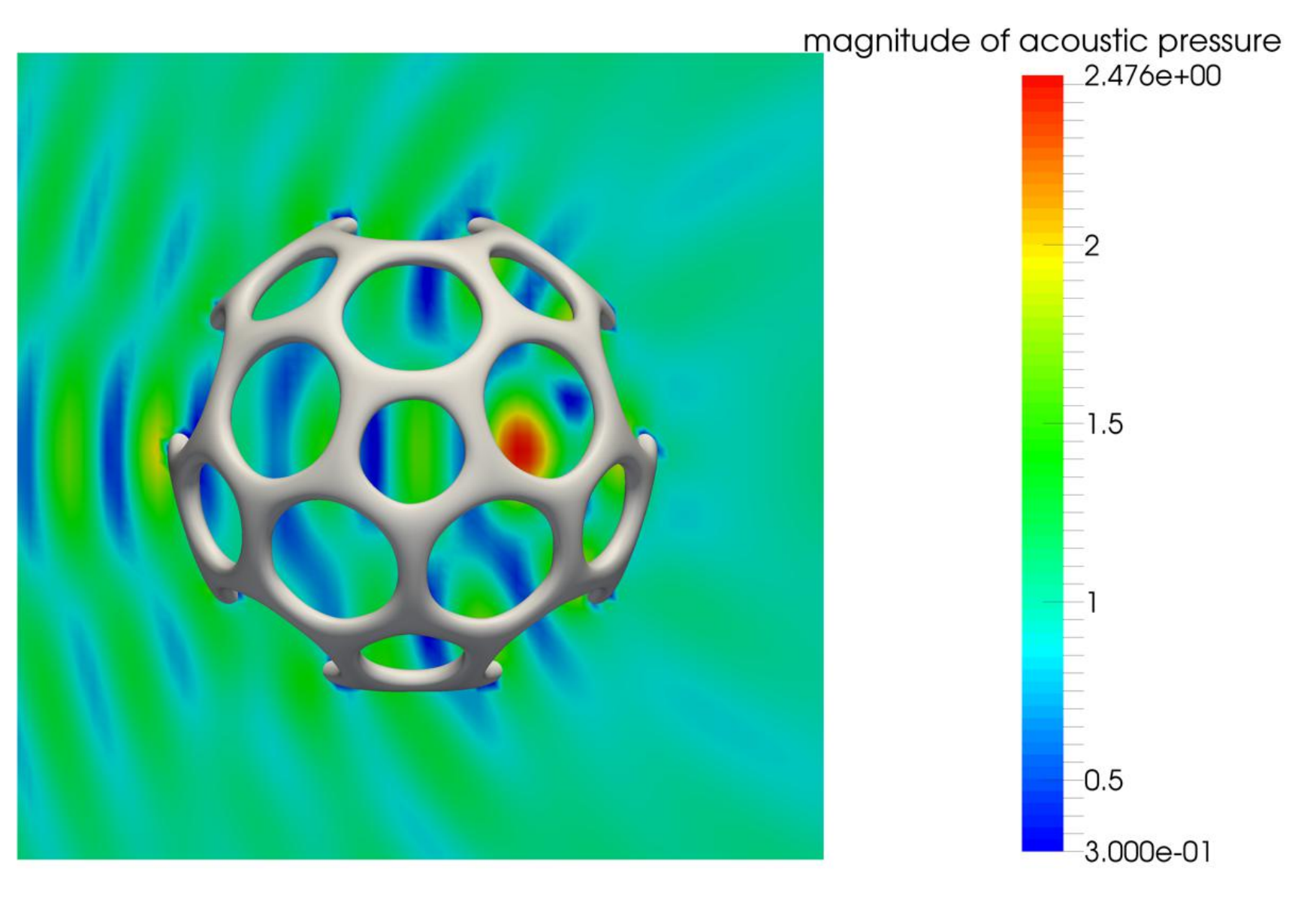}
        \caption{$|p|$ sampled over the $x\mbox{-}y$ plane.}
        \label{fig:molecule-xy-sampled}
    \end{subfigure}
    \caption{Coupled structural-acoustic analysis over a Loop subdivision surface with complex topology: acoustic pressure magnitudes for $ka=3.0$.}
    \label{fig:molecule-results}
\end{figure}

The final example we consider is that of a Loop subdivision surface with complex topology as illustrated in Figure~\ref{fig:molecule-discretisation}.   Such topologies are often challenging for parametric surfaces based on tensor product formulations but present no issues for subdivision surfaces.  The material properties and shell thickness as detailed in Table~\ref{tab:spherical-scattering-properties} are prescribed and a plane wave of unit magnitude travelling in the positive $x$ direction is specified with a normalised wavenumber of $ka=3$.  Plots of acoustic pressure magnitude sampled over the model surface and $x\mbox{-}y$ plane are illustrated in Figures~\ref{fig:molecule-surface-results} and \ref{fig:molecule-xy-sampled} respectively where a localised region of high acoustic pressure is created in the model interior. We remark that models with complex topology such as the present problem are often encountered in electromagnetic and acoustic scattering over metamaterial structures that exhibit non-intuitive scattered profiles and we envisage that our approach will have key benefits for such applications, particularly when applied to topology and shape optimisation.

% control grid and limit surface

% acoustic pressure plots

% *****************************************************
% ********** Conclusion *****************************
% *****************************************************
\section{Conclusion}

We have presented a novel BEM/FEM coupled method for structural acoustic analysis of shell geometries using Loop subdivision discretisations. Our approach utilises a collocation approach to generate a system of equations of the fluid domain using a boundary element formulation and a traditional Galerkin finite element method with Kirchhoff-Love shell theory to discretise the structural domain.   $\mathscr{H}$-matrices are employed to construct efficient low-rank approximations of dense matrices allowing for boundary element models generated through Loop subdivision discretisations with over 10,000 vertices (40,000 degrees of freedom).  We  verify our method through a closed-form solution of acoustic scattering over an elastic spherical shell geometry and demonstrate solutions of normalised wavenumbers up to $ka=80$.  Finally, the ability to model arbitrarily complex geometries with smooth surfaces generated through the Loop subdivision scheme is demonstrated thus highlighting the benefits of the present method for industrial design scenarios involving structural-acoustic interaction and optimisation of geometries with complex topologies.

\bibliographystyle{wileyj}
\bibliography{bibliography}

\begin{thebibliography}{10}
\providecommand{\url}[1]{\texttt{#1}}
\providecommand{\urlprefix}{URL }
\expandafter\ifx\csname urlstyle\endcsname\relax
  \providecommand{\doi}[1]{doi:\discretionary{}{}{}#1}\else
  \providecommand{\doi}{doi:\discretionary{}{}{}\begingroup
  \urlstyle{rm}\Url}\fi

\bibitem{cummer2016controlling}
Cummer SA, Christensen J, Al{\`u} A. Controlling sound with acoustic
  metamaterials. \emph{Nature Reviews Materials}  2016; \textbf{1}.

\bibitem{cirakscott2002}
Cirak F, Scott MJ, Antonsson EK, Ortiz M, Schr\"{o}der P. {Integrated modeling,
  finite-element analysis, and engineering design for thin-shell structures
  using subdivision}. \emph{Computer Aided Design}  2002;
  \textbf{34}(2):137--148.

\bibitem{cottrell2006isogeometric}
Cottrell JA, Reali A, Bazilevs Y, Hughes TJR. Isogeometric analysis of
  structural vibrations. \emph{Computer Methods in Applied Mechanics and
  Engineering}  2006; \textbf{195}(41):5257--5296.

\bibitem{SchDeScEvBoRaHu12}
Schillinger D, Ded\'{e} L, Scott MA, Evans JA, Borden MJ, Rank E, Hughes TJR.
  An isogeometric design-through-analysis methodology based on adaptive
  hierarchical refinement of {NURBS}, immersed boundary methods, and {T}-spline
  {CAD} surfaces. \emph{{Computer Methods in Applied Mechanics and
  Engineering}}  2012; \textbf{249-252}:116--150.

\bibitem{hughes2005isogeometric}
Hughes TJR, Cottrell JA, Bazilevs Y. Isogeometric analysis: {CAD}, finite
  elements, {NURBS}, exact geometry and mesh refinement. \emph{Computer Methods
  in Applied Mechanics and Engineering}  2005; \textbf{194}(39):4135--4195.

\bibitem{cirakortiz2000}
Cirak F, Ortiz M, Schr\"{o}der P. {Subdivision surfaces: a new paradigm for
  thin-shell finite-element analysis}. \emph{{International Journal for
  Numerical Methods in Engineering}}  2000; \textbf{47}(12):2039--2072.

\bibitem{nguyen2016c}
Nguyen T, Kar{\v{c}}iauskas K, Peters J. {$C^1$} finite elements on
  non-tensor-product 2d and 3d manifolds. \emph{Applied Mathematics and
  Computation}  2016; \textbf{272}:148--158.

\bibitem{majeedCirak:2016}
Majeed M, Cirak F. Isogeometric analysis using manifold-based smooth basis
  functions. \emph{Computer Methods in Applied Mechanics and Engineering}
  2017; \textbf{316}:547--567.

\bibitem{bazilevs2010isogeometric}
Bazilevs Y, Calo VM, Cottrell JA, Evans JA, Hughes TJR, Lipton S, Scott MA,
  Sederberg TW. Isogeometric analysis using {T}-splines. \emph{Computer Methods
  in Applied Mechanics and Engineering}  2010; \textbf{199}(5):229--263.

\bibitem{nguyen2011rotation}
Nguyen-Thanh N, Kiendl J, Nguyen-Xuan H, W{\"u}chner R, Bletzinger KU, Bazilevs
  Y, Rabczuk T. Rotation free isogeometric thin shell analysis using
  {PHT-splines}. \emph{Computer Methods in Applied Mechanics and Engineering}
  2011; \textbf{200}(47):3410--3424.

\bibitem{wei2015truncated}
Wei X, Zhang Y, Hughes TJR, Scott MA. Truncated hierarchical {Catmull--Clark}
  subdivision with local refinement. \emph{Computer Methods in Applied
  Mechanics and Engineering}  2015; \textbf{291}:1--20.

\bibitem{Grinspun:2002aa}
Grinspun E, Krysl P, Schr{\"o}der P. {CHARMS}: {A} simple framework for
  adaptive simulation. \emph{SIGGRAPH 2002 Conference Proceedings}, San
  Antonio, TX, 2002; 281--290.

\bibitem{zhang2016geometric}
Zhang YJ. \emph{Geometric Modeling and Mesh Generation from Scanned Images},
  vol.~6. CRC Press, 2016.

\bibitem{hu2016centroidal}
Hu K, Zhang YJ. Centroidal voronoi tessellation based polycube construction for
  adaptive all-hexahedral mesh generation. \emph{Computer Methods in Applied
  Mechanics and Engineering}  2016; \textbf{305}:405--421.

\bibitem{Sanches:2011aa}
Sanches R, Bornemann P, Cirak F. Immersed b-spline (i-spline) finite element
  method for geometrically complex domains. \emph{Computer Methods in Applied
  Mechanics and Engineering}  2011; \textbf{200}:1432--1445.

\bibitem{fritze2005fem}
Fritze D, Marburg S, Hardtke HJ. {FEM--BEM-coupling} and structural--acoustic
  sensitivity analysis for shell geometries. \emph{Computers \& structures}
  2005; \textbf{83}(2):143--154.

\bibitem{chen2014fem}
Chen L, Zheng C, Chen H. {FEM/wideband FMBEM} coupling for structural--acoustic
  design sensitivity analysis. \emph{Computer Methods in Applied Mechanics and
  Engineering}  2014; \textbf{276}:1--19.

\bibitem{webb1989absorbing}
Webb J, Kanellopoulos V. Absorbing boundary conditions for the finite element
  solution of the vector wave equation. \emph{Microwave and Optical Technology
  Letters}  1989; \textbf{2}(10):370--372.

\bibitem{safjan1998highly}
Safjan AJ. Highly accurate non-reflecting boundary conditions for finite
  element simulations of transient acoustics problems. \emph{Computer Methods
  in Applied Mechanics and Engineering}  1998; \textbf{152}(1):175--193.

\bibitem{djellouli2000finite}
Djellouli R, Farhat C, Macedo A, Tezaur R. Finite element solution of
  two-dimensional acoustic scattering problems using arbitrarily shaped convex
  artificial boundaries. \emph{Journal of Computational Acoustics}  2000;
  \textbf{8}(01):81--99.

\bibitem{harari2004analytical}
Harari I, Djellouli R. Analytical study of the effect of wave number on the
  performance of local absorbing boundary conditions for acoustic scattering.
  \emph{Applied Numerical Mathematics}  2004; \textbf{50}(1):15--47.

\bibitem{marburg2002six}
Marburg S. Six boundary elements per wavelength: Is that enough? \emph{Journal
  of Computational Acoustics}  2002; \textbf{10}(01):25--51.

\bibitem{marburg2003influence}
Marburg S, Schneider S. Influence of element types on numeric error for
  acoustic boundary elements. \emph{Journal of Computational Acoustics}  2003;
  \textbf{11}(03):363--386.

\bibitem{simpsonbordas2012}
Simpson RN, Bordas SPA, Trevelyan J, Rabczuk T. A two-dimensional isogeometric
  boundary element method for elastostatic analysis. \emph{{Computer Methods in
  Applied Mechanics and Engineering}}  2012; \textbf{209-212}:87--100.

\bibitem{simpsonscott2014}
Simpson RN, Scott MA, Taus M, Thomas DC, Lian H. Acoustic isogeometric boundary
  element analysis. \emph{{Computer Methods in Applied Mechanics and
  Engineering}}  2014; \textbf{269}:265--290.

\bibitem{peaketrevelyan2013}
Peake MJ, Trevelyan J, Coates G. Extended isogeometric boundary element method
  {(XIBEM)} for two-dimensional helmholtz problems. \emph{{Computer Methods in
  Applied Mechanics and Engineering}}  2013; \textbf{259}:93--102.

\bibitem{scottsimpson2012}
Scott MA, Simpson RN, Evans JA, Lipton S, Bordas SPA, Hughes TJR, Sederberg TW.
  Isogeometric boundary element analysis using unstructured {T}-splines.
  \emph{{Computer Methods in Applied Mechanics and Engineering}}  2013;
  \textbf{254}:197--221.

\bibitem{ginnis2014isogeometric}
Ginnis AI, Kostas K, Politis CG, Kaklis PD, Belibassakis KA, Gerostathis TP,
  Scott MA, Hughes TJR. Isogeometric boundary-element analysis for the
  wave-resistance problem using {T}-splines. \emph{Computer Methods in Applied
  Mechanics and Engineering}  2014; \textbf{279}:425--439.

\bibitem{bandara2015boundary}
Bandara K, Cirak F, Of G, Steinbach O, Zapletal J. Boundary element based
  multiresolution shape optimisation in electrostatics. \emph{Journal of
  Computational Physics}  2015; \textbf{297}:584--598.

\bibitem{kiendl2009isogeometric}
Kiendl J, Bletzinger KU, Linhard J, W{\"u}chner R. Isogeometric shell analysis
  with {Kirchhoff--Love} elements. \emph{Computer Methods in Applied Mechanics
  and Engineering}  2009; \textbf{198}(49):3902--3914.

\bibitem{benson2010isogeometric}
Benson D, Bazilevs Y, Hsu MC, Hughes TJR. Isogeometric shell analysis: the
  {Reissner--Mindlin} shell. \emph{Computer Methods in Applied Mechanics and
  Engineering}  2010; \textbf{199}(5):276--289.

\bibitem{benson2011large}
Benson D, Bazilevs Y, Hsu MC, Hughes TJR. A large deformation, rotation-free,
  isogeometric shell. \emph{Computer Methods in Applied Mechanics and
  Engineering}  2011; \textbf{200}(13):1367--1378.

\bibitem{Cirak:2001aa}
Cirak F, Ortiz M. Fully ${C}^1$-conforming subdivision elements for finite
  deformation thin-shell analysis. \emph{International Journal for Numerical
  Methods in Engineering}  2001; \textbf{51}:813--833.

\bibitem{cirak2011subdivision}
Cirak F, Long Q. Subdivision shells with exact boundary control and
  non-manifold geometry. \emph{International Journal for Numerical Methods in
  Engineering}  2011; \textbf{88}(9):897--923.

\bibitem{Long:2012aa}
Long Q, Bornemann PB, Cirak F. Shear-flexible subdivision shells.
  \emph{International Journal for Numerical Methods in Engineering}  2012;
  \textbf{90}:1549--1577.

\bibitem{loop1987smooth}
Loop C. Smooth subdivision surfaces based on triangles. Master's {T}hesis,
  Department of Mathematics, The University of Utah 1987.

\bibitem{heltai2016interaction}
Heltai L, Kiendl J, DeSimone A, Reali A. A natural framework for isogeometric
  fluid--structure interaction based on {BEM}--shell coupling. \emph{Computer
  Methods in Applied Mechanics and Engineering}  2016; (in press).

\bibitem{Zorin:2000aa}
Zorin D, Schr{\"o}der P. Subdivision for modeling and animation. SIGGRAPH 2000
  Course Notes 2000.

\bibitem{Peters:2008aa}
Peters J, Reif U. \emph{Subdivision Surfaces}. Springer Series in Geometry and
  Computing, Springer, 2008.

\bibitem{Stam1998a}
Stam J. {Exact evaluation of Loop subdivision surfaces at arbitrary parameter
  values}. \emph{SIGGRAPH Course Note}  1998; :111--124.

\bibitem{schillinger2013isogeometric}
Schillinger D, Evans JA, Reali A, Scott MA, Hughes TJR. Isogeometric
  collocation: Cost comparison with galerkin methods and extension to adaptive
  hierarchical {NURBS} discretizations. \emph{Computer Methods in Applied
  Mechanics and Engineering}  2013; \textbf{267}:170--232.

\bibitem{liu1991some}
Liu Y, Rudolphi T. Some identities for fundamental solutions and their
  applications to weakly-singular boundary element formulations.
  \emph{Engineering analysis with boundary elements}  1991;
  \textbf{8}(6):301--311.

\bibitem{Ciarlet:2005aa}
Ciarlet PG. \emph{An Introduction to Differential Geometry with Applications to
  Elasticity}. Springer, 2005.

\bibitem{gough2009gnu}
Gough B. \emph{GNU scientific library reference manual}. Network Theory Ltd.,
  2009.

\bibitem{eigenweb}
Guennebaud G, Jacob B, \emph{et~al.}. Eigen v3. http://eigen.tuxfamily.org
  2010.

\bibitem{trilinos-overview}
Heroux M, Bartlett R, Hoekstra VHR, Hu J, Kolda T, Lehoucq R, Long K, Pawlowski
  R, Phipps E, Salinger A, \emph{et~al.}. {An Overview of Trilinos}.
  \emph{Technical {R}eport SAND2003-2927}, Sandia National Laboratories 2003.

\bibitem{kriemann2008hlibpro}
Kriemann R. {HLIBpro} user manual. \emph{Max-Planck-Institute for Mathematics
  in the Sciences, Leipzig}  2008; .

\bibitem{hackbusch2004hierarchical}
Hackbusch W, Khoromskij BN, Kriemann R. Hierarchical matrices based on a weak
  admissibility criterion. \emph{Computing}  2004; \textbf{73}(3):207--243.

\bibitem{bebendorf2008hierarchical}
Bebendorf M. \emph{Hierarchical matrices}. Springer, 2008.

\bibitem{junger1986sound}
Junger MC, Feit D. \emph{Sound, structures, and their interaction}, vol. 225.
  MIT press Cambridge, MA, 1986.

\end{thebibliography}

\end{document}